# Multiple orthogonal polynomials: Pearson equations and Christoffel formulas

Amílcar Branquinho[1] · Ana Foulquié-Moreno[2] · Manuel Mañas[3,4]

## Abstract
Multiple orthogonal polynomials with respect to two weights on the step-line are considered. A connection between different dual spectral matrices, one banded (recursion matrix) and one Hessenberg, respectively, and the Gauss–Borel factorization of the moment matrix is given. It is shown a *hidden* freedom exhibited by the spectral system related to the multiple orthogonal polynomials. Pearson equations are discussed, a Laguerre–Freud matrix is considered, and differential equations for type I and II multiple orthogonal polynomials, as well as for the corresponding linear forms are given. The Jacobi–Piñeiro multiple orthogonal polynomials of type I and type II are used as an illustrating case and the corresponding differential relations are presented. A permuting Christoffel transformation is discussed, finding the connection between the different families of multiple orthogonal polynomials. The Jacobi–Piñeiro case provides a convenient illustration of these symmetries, giving linear relations between different polynomials with shifted and permuted parameters. We also present the general theory for the perturbation of each weight by a different polynomial or rational function aka called Christoffel and Geronimus transformations. The connections formulas between the type II multiple orthogonal polynomials, the type I linear forms,

✉ Manuel Mañas
manuel.manas@ucm.es

Amílcar Branquinho
ajplb@mat.uc.pt

Ana Foulquié-Moreno
foulquie@ua.pt

1 Departamento de Matemática, Universidade de Coimbra, 3001-454 Coimbra, Portugal

2 Departamento de Matemática, Universidade de Aveiro, 3810-193 Aveiro, Portugal

3 Departamento de Física Teórica, Universidad Complutense de Madrid, Plaza Ciencias 1, 28040 Madrid, Spain

4 Instituto de Ciencias Matematicas (ICMAT), Campus de Cantoblanco UAM, 28049 Madrid, Spain

as well as the vector Stieltjes–Markov vector functions is also presented. We illustrate these findings by analyzing the special case of modification by an even polynomial.

## Contents

## 1 Introduction

Multiple orthogonality is a very close topic to that of simultaneous rational approximation (simultaneous Padé aproximants) of systems of Cauchy transforms of measures. The history of simultaneous rational approximation starts in 1873 with the well known article [28] in where Charles Hermite proved the transcendence of Euler's constant $e$. Later on, along the years 1934-35, Kurt Mahler delivered in the University of Groningen several lectures [35] where he settled down the foundations of this theory. In the mean time, two Mahler's students, Coates and Jager, made important contributions in this respect (see [19, 30]).

There are two formulations of multiple orthogonality, the so called type I and type II ones, see [42]. Both are equivalent in the sense of the duality. If we set a problem involving orthogonality conditions regarding several measures (or several weights) we call it a type II problem. Under this view, the fundamental objects are polynomials. The dual objects for these polynomials are linear forms considered in the type I version of multiple orthogonality. Due to the fact that perhaps type II version is more *natural*, the research interest has been centered in the characterization of this type of systems (cf. [5]).

The Angelescu and the Nikishin [4, 38, 43] are the two main systems of multiple orthogonal polynomials that have been discussed in the literature. Contributions regarding zero distribution, interlacing property and confinement of the zeros, are studied in [21, 25, 27, 31–34]. Multiple orthogonal polynomials also have differ-ent expressions for the Christoffel–Darboux kernel, coming from the different ways of defining them from the recurrence relation [20, 22, 41]. Studying the asymptotic behavior of this kernel is essential in the study of universality classes in random matrix theory due to its connection with eigenvalue distribution. Recent trends of multiple orthogonal polynomials can be seen in [29]. For multiple Gaussian quadratures see [15, 22] and for Laguerre–Freud equations, Darboux transformations and related perturbations see [11–14, 37]. Laguerre–Freud equations play a key role in the study of the coefficients of the three term recurrence relation satisfied by monic standard orthogonal polynomials [14] and see also [37] for the so called structure relation.

Multiple orthogonal polynomials can be built recursively, once we choose a chain of indexes. In [1] the chain of indexes is given through the generalized Euclidean division, and using the Gauss–Borel factorization of the related matrix of moments we find that the type I and type II orthogonality conditions appear now clearly related through bi-orthogonality conditions.

In [6, 44] multiple orthogonal polynomials with respect to $p$ weights satisfying Pearson's differential equation are studied, giving a classification of type II multiple orthogonal polynomials that can be represented by a product of commuting Rodrigues type operators.

This paper is devoted to the study of multiple orthogonal polynomials with two weights on the step-line, also known in the literature as 2-orthogonal polynomials in the case of type II (see [36], among others). In the research on the relation between Markov chains beyond birth and death processes and multiple orthogonal polynomials [17, 18] we get, as a byproduct, some results that hold for families of these type I and II simple cases. Moreover, in order to better understand some of these phenomena we turned back to the definition of multiple orthogonality and find, to the best of our knowledge, some unknown facts so far. Given their general character and interest, we decided to collect them in a separate publication. Let us now give an account of the contents of this paper.

First we apply the Gauss–Borel factorization described in [1]. Within this description we introduce, associated to the lower and upper triangular factors of the moment matrix, two matrices $J$ and $\tilde{J}$ having the sequences of orthogonal polynomials of type II and type I, respectively, as eigenvectors. The first one is the recursion matrix, which happens to be a banded matrix with four non-vanishing diagonals, the first superdiagonal, the main diagonal and the two first subdiagonals. However, $\tilde{J}^\top$ is not a band matrix but an upper Hessenberg matrix so that all its superdiagonals, but for the first one, are zero. Nevertheless, its square $\tilde{J}^2$ happens to be equal to the transpose of $J$ and, consequently, is a banded matrix. These relations lead to expressions of the matrix $\tilde{J}$ in terms of the four nonzero diagonals of the recursion matrix $J$. This is described in Sect. 1.2.

In Sect. 1.3 it is shown that for multiple orthogonal polynomials the recursion matrix $J$, or equivalently the sequence of type II multiple orthogonal polynomials $\{B^{(n)}\}_{n=0}^{\infty}$, related to the spectral system $(w_1, w_2, \mathrm{d}\,\mu)$, determines uniquely $w_1$. However, the second weight $w_2$ is not uniquely determined. We describe this phenomena as a *gauge freedom*. It is also shown that the type I linear form associated to the sequence of type II multiple orthogonal polynomials $\{B^{(n)}\}_{n=0}^{\infty}$ is also uniquely determined. Examples of this gauge freedom with almost uniform Jacobi matrices are obtained in [18].

In Sect. 2 we assume that the two weights $(w_1, w_2)$ satisfy first order scalar Pearson equations, as it does for instance the Jacobi–Piñeiro. We are able to find in the matrix setting a symmetry for the matrix of moments and derived a Laguerre–Freud matrix, a finite band matrix with only four nonzero diagonals, that accounts for the consequences. These ideas lead to differential equations satisfied by the multiple orthogonal polynomials of type II as well as for the corresponding linear forms of type I.

Then, in Sect. 3 we deal with a nice symmetry easily described in this two-component case. There is a relation between the moment matrices with weights $(w_1, w_2)$ and $(xw_2, w_1)$, that we call permuting Christoffel symmetry. The connection matrix for both multiple orthogonal polynomials of type II and linear forms of type I is remarkably simple in this case. We will apply it to the case of Jacobi–Piñeiro multiple orthogonal polynomials [39, 42], that we take as *case study*. This provides explicit formulas connecting the polynomials and linear forms with parameters $(\alpha, \beta, \gamma)$ and $(\beta, \alpha + 1, \gamma)$. We also present the Christoffel formula when the transformation is the multiplication of both weights by $x$, with no permutation. Both types of Christoffel transformations are relevant in the understanding of the mentioned uniform Jacobi matrices related to hypergeometric weights see [1]. Then, motivated by the previous examples we analyze the general Christoffel and Geronimus transformations for a system of weights $(w_1, w_2)$, and we end the paper by taking the modification by an even polynomial as a case study of these two type of transformations.

### 1.1 Two component multiple orthogonal polynomials on the step-line

We present here well basic results regarding multiple orthogonal polynomials, see for example [42]. We follow the Gauss–Borel factorization problem approach of [1]. We consider

$$X := \begin{bmatrix} 1 & x & x^2 \cdots \end{bmatrix}^\top,$$

and, in terms of a system of a couple of weights $\vec{w} = (w_1, w_2)$, the vector of monomials

$$X_1(x) = \begin{bmatrix} 1 & 0 & x & 0 & x^2 & 0 \cdots \end{bmatrix}^\top, \qquad X_2(x) = \begin{bmatrix} 0 & 1 & 0 & x & 0 & x^2 \cdots \end{bmatrix}^\top,$$

we define the following *vector of undressed linear forms*

$$\xi := X_1 w_1 + X_2 w_2 = \begin{bmatrix} w_1 & w_2 & x w_1 & x w_2 & x^2 w_1 & x^2 w_2 \cdots \end{bmatrix}^\top. \tag{1}$$

**Definition 1** For a given measure $\mu$ with support on $\Delta$, an interval in $\mathbb{R}$, and weights $w_1$ and $w_2$ as above, the moment matrix is given by

$$g := \int_\Delta X(x)(\xi(x))^\top \, \mathrm{d}\, \mu(x). \tag{2}$$

**Definition 2** The Gauss–Borel factorization of the moment matrix $g$ is the problem of finding the solution of

$$g = S^{-1} H \tilde{S}^{-\top}, \tag{3}$$

with $S$, $\tilde{S}$ lower unitriangular semi-infinite matrices

$$S = \begin{bmatrix} 1 & 0 & \cdots \\ S_{1,0} & 1 & \ddots \\ S_{2,0} & S_{2,1} & \ddots \\ \vdots & \ddots & \ddots \end{bmatrix}, \qquad \tilde{S} = \begin{bmatrix} 1 & 0 & \cdots \\ \tilde{S}_{1,0} & 1 & \ddots \\ \tilde{S}_{2,0} & \tilde{S}_{2,1} & \ddots \\ \vdots & \ddots & \ddots \end{bmatrix},$$

and $H$ a semi-infinite diagonal matrix with diagonal entries $H_l \neq 0$, $l \in \mathbb{N}_0$.

Let us assume that the Gauss–Borel factorization exists, that is, the system $(\vec{w}, \mathrm{d}\,\mu)$ is perfect. In terms of $S$ and $\tilde{S}$ we construct the type II multiple orthogonal polynomials

$$B^{(m)} := x^m + \sum_{i=0}^{m-1} S_{m,i} x^i, \tag{4}$$

with $m \in \mathbb{N}_0$, as well as the type I multiple orthogonal polynomials

$$A_1^{(2m)} = \frac{1}{H_{2m}} \left( x^m + \sum_{i=0}^{m-1} \tilde{S}_{2m,2i} x^i \right),$$

$$A_1^{(2m+1)} = \frac{1}{H_{2m+1}} \left( \sum_{i=0}^{m} \tilde{S}_{2m+1,2i} x^i \right), \quad m \in \mathbb{N}_0, \tag{5}$$

and

$$A_2^{(0)} = 0, \qquad A_2^{(1)} = 1,$$

$$A_2^{(2m)} = \frac{1}{H_{2m}} \left( \sum_{i=0}^{m-1} \tilde{S}_{2m,2i+1} x^i \right), \quad A_2^{(2m+1)} = \frac{1}{H_{2m+1}} \left( x^m + \sum_{i=0}^{m-1} \tilde{S}_{2m+1,2i+1} x^i \right), \quad m \in \mathbb{N}. \tag{6}$$

For $m \in \mathbb{N}_0$, the linear forms are

$$Q^{(m)} := w_1 A_1^{(m)} + w_2 A_2^{(m)}.$$

**Proposition 1** (Type I and II multiorthogonality relations) *In terms of the weight vectors, $\vec{\nu}(2m) = (m+1, m)$ and $\vec{\nu}(2m+1) = (m+1, m+1)$, $m \in \mathbb{N}_0$, and corresponding multiple orthogonal polynomials*

$$B^{(2m)} = B_{(m,m)}, \qquad B^{(2m+1)} = B_{(m+1,m)},$$

$$A_1^{(2m)} = A_{(m+1,m),1} \qquad A_1^{(2m+1)} = A_{(m+1,m+1),1},$$
$$A_2^{(2m)} = A_{(m+1,m),2}, \qquad A_2^{(2m+1)} = A_{(m+1,m+1),2},$$

*the following type I orthogonality relations*

$$\int_\Delta x^j (A_{\vec{\nu},1}(x)w_1(x) + A_{\vec{\nu},2}(x)w_2(x))\,\mathrm{d}\,\mu(x) = 0,$$
$$\deg A_{\vec{\nu},1} \le \nu_1 - 1, \quad \deg A_{\vec{\nu},2} \le \nu_2 - 1,$$

*for* $j \in \{0, \dots, |\vec{\nu}| - 2\}$*, and type II orthogonality relations*

$$\int_\Delta B_{\vec{\nu}}(x)w_a(x)x^j\,\mathrm{d}\,\mu(x) = 0, \quad \deg B_{\vec{\nu}} \le |\vec{\nu}|, \quad j = 0, \dots, \nu_a - 1, \quad a = 1, 2, \tag{7}$$

*are fulfilled.*

Defining

$$\Lambda = \begin{bmatrix} 0 & 1 & 0 & 0 & \cdots \\ 0 & 0 & 1 & 0 & \ddots \\ 0 & 0 & 0 & 1 & \ddots \\ \vdots & \ddots & \ddots & \ddots & \ddots \end{bmatrix},$$

the lower unitriangular factors $S$, $\tilde{S}$ can be written in terms of its subdiagonals as follows

$$S = I + \Lambda^\top S^{[1]} + \left(\Lambda^\top\right)^2 S^{[2]} + \cdots, \qquad \tilde{S} = I + \Lambda^\top \tilde{S}^{[1]} + \left(\Lambda^\top\right)^2 \tilde{S}^{[2]} + \cdots,$$

with $S^{[k]}$, $\tilde{S}^{[k]}$ diagonal matrices with entries $S_l^{[k]}$, $\tilde{S}_l^{[k]}$, $l \in \mathbb{N}_0$, respectively. Hence,

$$S_m^{[k]} = S_{k+m,m}, \qquad \tilde{S}_m^{[k]} = \tilde{S}_{k+m,m},$$

so that $B^{(m)} = x^m + \sum_{i=0}^{m-1} S_i^{[m-i]} x^i$ and

$$A_1^{(2m)} = \frac{1}{H_{2m}} \left( x^m + \sum_{i=0}^{m-1} \tilde{S}_{2i}^{[2(m-i)]} x^i \right),$$
$$A_1^{(2m+1)} = \frac{1}{H_{2m+1}} \left( \sum_{i=0}^{m} \tilde{S}_{2i}^{[2(m-i)+1]} x^i \right), \quad m \in \mathbb{N}_0,$$

$$A_2^{(2m)} = \frac{1}{H_{2m}}\left(\sum_{i=0}^{m-1} \tilde{S}_{2i+1}^{[2(m-i)-1]} x^i\right),$$
$$A_2^{(2m+1)} = \frac{1}{H_{2m+1}}\left(x^m + \sum_{i=0}^{m-1} \tilde{S}_{2i+1}^{[2(m-i)]} x^i\right), \quad m \in \mathbb{N}.$$

**Definition 3** Vector of type II multiple orthogonal polynomials and of type I linear forms associated with $(w_1, w_2, \mathrm{d}\,\mu)$ are defined, respectively, by

$$B := \begin{bmatrix} B^{(0)} \\ B^{(1)} \\ \vdots \end{bmatrix} = S\,X, \tag{8}$$

$$A_1 := \begin{bmatrix} A_1^{(0)} \\ A_1^{(1)} \\ \vdots \end{bmatrix} = H^{-1}\,\tilde{S}\,X_1, \qquad A_2 := \begin{bmatrix} A_2^{(0)} \\ A_2^{(1)} \\ \vdots \end{bmatrix} = H^{-1}\,\tilde{S}\,X_2, \qquad \textit{and} \qquad Q := \begin{bmatrix} Q^{(0)} \\ Q^{(1)} \\ \vdots \end{bmatrix} = H^{-1}\,\tilde{S}\,\xi. \tag{9}$$

**Proposition 2** (Biorthogonality) *The following multiple biorthogonality relations*

$$\int_\Delta B^{(m)}(x) Q^{(k)}(x)\,\mathrm{d}\,\mu(x) = \delta_{m,k}, \qquad m, k \in \mathbb{N}_0, \tag{10}$$

*hold.*

### 1.2 Recurrence relations and the shift matrices

**Definition 4** *(Shifted matrices)* We introduce the following shift matrices

$$\Lambda_1 = \begin{bmatrix} 0 & 0 & 1 & 0 & \cdots \\ \vdots & \ddots & \ddots & 0 & \ddots \\ \vdots & & \ddots & \ddots & 1 & \ddots \\ \vdots & & & \ddots & \ddots & 0 & \ddots \\ \vdots & & & & \ddots & \ddots & 1 & \ddots \end{bmatrix}, \qquad \Lambda_2 = \begin{bmatrix} 0 & 0 & 0 & 0 & \cdots \\ \vdots & \ddots & \ddots & 1 & \ddots \\ \vdots & & \ddots & \ddots & 0 & \ddots \\ \vdots & & & \ddots & \ddots & 1 & \ddots \\ \vdots & & & & \ddots & \ddots & 0 & \ddots \end{bmatrix},$$

with

$$\Lambda_1 + \Lambda_2 = \Lambda^2 = \begin{bmatrix} 0 & 0 & 1 & 0 & \cdots \\ \vdots & \ddots & \ddots & 1 & \ddots \\ \vdots & & \ddots & \ddots & 1 & \ddots \\ \vdots & & & \ddots & \ddots & 1 & \ddots \\ \vdots & & & & \ddots & \ddots & 1 & \ddots \end{bmatrix}. \tag{11}$$

From this definition we get the following technical lemma that determines the algebra of the multiple orthogonality studied in this work.

**Lemma 1** *The shift matrices satisfy*

$$\Lambda\, X(x) = x\, X(x), \qquad \Lambda_1\, X_1(x) = x\, X_1(x), \qquad \Lambda_2\, X_2(x) = x\, X_2(x), \tag{12}$$

*as well as*

$$\Lambda_1\, X_2(x) = \mathbf{0}, \qquad\qquad \Lambda_2\, X_1(x) = \mathbf{0}.$$

*Moreover, the projection matrices*

$$\Pi_1 := \begin{bmatrix} 1 & 0 & \cdots & & \\ 0 & 0 & \ddots & & \\ \vdots & \ddots & 1 & \ddots & \\ \vdots & & \ddots & 0 & \ddots \end{bmatrix}, \qquad \Pi_2 := \begin{bmatrix} 0 & 0 & \cdots & & \\ 0 & 1 & \ddots & & \\ \vdots & \ddots & 0 & \ddots & \\ \vdots & & \ddots & 1 & \ddots \end{bmatrix},$$

*satisfy*

$$I = \Pi_1 + \Pi_2, \qquad \Pi_1^2 = \Pi_1, \qquad \Pi_2^2 = \Pi_2, \qquad \Pi_1\, \Pi_2 = \Pi_2\, \Pi_1 = \mathbf{0}.$$

**Proposition 3** (Bi-Hankel structure of the moment matrix) *The moment matrix fulfills the symmetry condition*

$$\Lambda\, g = g\, (\Lambda^\top)^2. \tag{13}$$

***Remark 1*** Therefore, it satisfies a bi-Hankel condition

$$g_{n+1,m} = g_{n,m+2}, \qquad\qquad n, m \in \mathbb{N}_0,$$

and we can write the moment matrix in terms of

$$\vec{g}_n = (g_{n,n}, g_{n,n+1})$$

as the following *vectorial Hankel* type matrix

$$g = \begin{bmatrix} \vec{g}_0 & \vec{g}_1 & \vec{g}_2 & \cdots\cdots & \vec{g}_n & \cdots\cdots \\ \vec{g}_1 & \vec{g}_2 & & & & \\ \vec{g}_2 & & & \iddots & & \\ \vdots & & \iddots & & & \\ \vdots & \iddots & & & & \\ \vec{g}_n & & & & & \\ \vdots & & & & & \\ \vdots & & & & & \end{bmatrix}.$$

**Proposition 4** (Recursion matrix) *The matrix*

$$T := S\Lambda S^{-1} = \left(H^{-1}\tilde{S}\Upsilon\tilde{S}^{-1}H\right)^{\top}.$$

*is a Hessenberg matrix of the form*

$$T = \begin{bmatrix} b_{0,0} & 1 & 0 & \cdots & & & \\ c_{1,0} & b_{1,0} & 1 & \ddots & & & \\ d_{1,1} & c_{1,1} & b_{1,1} & 1 & \ddots & & \\ 0 & d_{2,1} & c_{2,1} & b_{2,1} & 1 & \ddots & \\ \vdots & \ddots & \ddots & \ddots & \ddots & \ddots & \ddots \end{bmatrix}, \tag{14}$$

*and $T$ is a tetradiagonal matrix.*

**Theorem 1** (Recursion relations)

*(i) The recursion matrix, type II multiple orthogonal polynomials, and corresponding type I multiple orthogonal polynomials and linear forms of type I fulfill the eigenvalue property*

$$T\,B = x\,B, \qquad T^{\top}A_1 = x\,A_1, \qquad T^{\top}A_2 = x\,A_2, \qquad T^{\top}Q = x\,Q. \tag{15}$$

*Componentwise, this eigenvalue property leads to the recursion relations.*

*(ii) The lower Hessenberg matrix*

$$\hat{T} := H^{-1}\tilde{T}H, \qquad \text{with} \qquad \tilde{T} := \tilde{S}\Lambda\tilde{S}^{-1},$$

*is such that*

$$T^{\top} = \hat{T}^2, \tag{16}$$

*and has the following eigenvalue property*

$$\hat{T}\,A_1 = x\,A_2, \qquad \hat{T}\,A_2 = A_1. \tag{17}$$

*Moreover, if $T_1 := H^{-1}\tilde{S}\Lambda_1\tilde{S}^{-1}H$ and $T_2 := H^{-1}\tilde{S}\Lambda_2\tilde{S}^{-1}H$, then*

$$T_1\,A_1 = x\,A_1, \qquad T_2\,A_2 = x\,A_2, \tag{18}$$

*as well as,*

$$T_1\,A_2 = \mathbf{0}, \qquad T_2\,A_1 = \mathbf{0}.$$

***Proof*** Equation (16) follows from the representation of $T$ in Proposition 4 and taking into account (11). From

$$\Lambda\, X_1(x) = x\, X_2(x), \qquad \Lambda\, X_2(x) = X_1(x). \tag{19}$$

we deduce (17).

To deduce (18) use (12) and the definition of $A_1$ and $A_2$ in (9). □

***Remark 2*** (i) Despite $\hat{T}$ is not a banded matrix, its square $\hat{T}^2$ is a tetradiagonal matrix.

(ii) Consequently, we obtain the spectral equations

$$\hat{T}^2 A_1 = x\, A_1, \qquad \hat{T}^2 A_2 = x\, A_2.$$

(iii) Notice that $\hat{T}^2 = T^\top = T_1 + T_2$. Moreover $T_1 T_2 = T_2 T_1 = \mathbf{0}$ so that for any polynomial $p(x)$ we have $p(T^\top) = p(T_1 + T_2) = p(T_1) + p(T_2)$.

**Definition 5** *(Shift operators)* The *shift operators* $\mathfrak{a}_\pm$ acts over the diagonal matrices as follows

$$\mathfrak{a}_-\operatorname{diag}(d_0, d_1, \dots) := \operatorname{diag}(d_1, d_2, \dots) \qquad \mathfrak{a}_+\operatorname{diag}(d_0, d_1, \dots) := \operatorname{diag}(0, d_0, d_1, \dots).$$

**Lemma 2** *Shift operators have the following properties,*

$$\Lambda\, D = (\mathfrak{a}_-\, D)\Lambda, \quad A\,\Lambda = \Lambda(\mathfrak{a}_+\, D), \quad A\,\Lambda^\top = \Lambda^\top(\mathfrak{a}_-\, D), \quad \Lambda^\top A = (\mathfrak{a}_+\, D)\Lambda^\top, \tag{20}$$

*for any diagonal matrix D.*

A simple computation leads to:

**Proposition 5** *The inverse matrix $S^{-1}$ of a lower unitriangular matrix $S$ expands in terms of subdiagonals as follows*

$$S^{-1} = I + \Lambda^\top S^{[-1]} + \left(\Lambda^\top\right)^2 S^{[-2]} + \cdots,$$

*with first few subdiagonals $S^{[-k]}$ given by*

$$\begin{aligned}
S^{[-1]} &= -S^{[1]},\\
S^{[-2]} &= -S^{[2]} + (\mathfrak{a}_-\, S^{[1]})S^{[1]},\\
S^{[-3]} &= -S^{[3]} + (\mathfrak{a}_-\, S^{[2]})S^{[1]} + (\mathfrak{a}_-^2\, S^{[1]})S^{[2]} - (\mathfrak{a}_-^2\, S^{[1]})(\mathfrak{a}_-\, S^{[1]})S^{[1]},\\
S^{[-4]} &= -S^{[4]} + (\mathfrak{a}_-\, S^{[3]})S^{[1]} + (\mathfrak{a}_-^2\, S^{[2]})S^{[2]}\\
&\quad - (\mathfrak{a}_-^2\, S^{[2]})(\mathfrak{a}_-\, S^{[1]})S^{[1]} + (\mathfrak{a}_-^3\, S^{[1]})S^{[3]}\\
&\quad - (\mathfrak{a}_-^3\, S^{[1]})(\mathfrak{a}_-\, S^{[2]})S^{[1]} - (\mathfrak{a}_-^3\, S^{[1]})(\mathfrak{a}_-^2\, S^{[1]})S^{[2]}\\
&\quad + (\mathfrak{a}_-^3\, S^{[1]})(\mathfrak{a}_-^2\, S^{[1]})(\mathfrak{a}_-\, S^{[1]})S^{[1]}.
\end{aligned}$$

**Lemma 3** *For the recursion matrix*

$$T = S\,\Lambda\,S^{-1} = \Lambda + T^{[0]} + \Lambda^\top T^{[1]} + \big(\Lambda^\top\big)^2 T^{[2]} + \cdots,$$

*we have*

$$\begin{aligned}
T^{[0]} &= \mathfrak{a}_+ S^{[1]} - S^{[1]},\\
T^{[1]} &= \mathfrak{a}_+ S^{[2]} - S^{[2]} + \big(\mathfrak{a}_-\,S^{[1]} - S^{[1]}\big)S^{[1]},\\
T^{[2]} &= \mathfrak{a}_+\,S^{[3]} - S^{[3]} + \big(\mathfrak{a}_-\,S^{[2]}\big)S^{[1]} + \big(\mathfrak{a}_-^2\,S^{[1]}\big)S^{[2]}\\
&\quad - \big(\mathfrak{a}_-^2\,S^{[1]}\big)\big(\mathfrak{a}_-\,S^{[1]}\big)S^{[1]} - S^{[2]}S^{[1]}\\
&\quad + \big(\mathfrak{a}_+\,S^{[1]}\big)\big(-S^{[2]} + \big(\mathfrak{a}_-\,S^{[1]}\big)S^{[1]}\big).
\end{aligned}$$

***Proof*** We consider

$$\begin{aligned}
T &= \Big(I + \Lambda^\top S^{[1]} + \big(\Lambda^\top\big)^2 S^{[2]} + \cdots\Big)\,\Lambda\,\Big(I + \Lambda^\top S^{[-1]} + \big(\Lambda^\top\big)^2 S^{[-2]} + \cdots\Big)\\
&= \Lambda + \mathfrak{a}_+\,S^{[1]} + S^{[-1]} + \Lambda^\top\Big(\mathfrak{a}_+\,S^{[2]} + S^{[-2]} + S^{[1]}S^{[-1]}\Big)\\
&\quad + \big(\Lambda^\top\big)^2\Big(\mathfrak{a}_+\,S^{[3]} + S^{[-3]} + S^{[2]}S^{[-1]} + \big(\mathfrak{a}_+\,S^{[1]}\big)S^{[-2]}\Big) + \cdots.
\end{aligned}$$

Now, from Proposition 5 we get the desired representation. □

**Lemma 4** *We have the following diagonal expansion*

$$\begin{aligned}
\hat{T}^2 &= H(\mathfrak{a}_-^2 H)^{-1}\Lambda^2 + H(\mathfrak{a}_- H)^{-1}\big(\mathfrak{a}_-\tilde{T}^{[1]} + \tilde{T}^{[0]}\big)\Lambda + \sum_{k=0}^{\infty}(\Lambda^\top)^k\hat{T}^{[2,k]},\\
\hat{T}^{[2,n]} &:= T^{[2,n]}(\mathfrak{a}_-^n H)H^{-1}.
\end{aligned}$$

***Proof*** Observe that

$$\tilde{T} = \Lambda + \tilde{T}^{[0]} + \Lambda^\top\tilde{T}^{[1]} + (\Lambda^\top)^2\tilde{T}^{[2]} + \cdots,$$

so that

$$\hat{T} = H(\mathfrak{a}_- H)^{-1}\Lambda + \hat{T}^{[0]} + \Lambda^\top\hat{T}^{[1]} + (\Lambda^\top)^2\hat{T}^{[2]} + \cdots, \quad \hat{T}^{[k]} = \tilde{T}^{[k]}(\mathfrak{a}_-^k H)H^{-1}.$$

Hence,

$$\begin{aligned}
\tilde{T}^2 &= (\Lambda + \tilde{T}^{[0]} + \Lambda^\top\tilde{T}^{[1]} + (\Lambda^\top)^2\tilde{T}^{[2]} + \cdots)^2 = \Lambda^2 + \big(\mathfrak{a}_-\tilde{T}^{[0]} + \tilde{T}^{[0]}\big)\Lambda\\
&\quad + \sum_{k=0}^{\infty}(\Lambda^\top)^k\tilde{T}^{[2,k]},\\
\tilde{T}^{[2,n]} &:= \tilde{T}^{[n+1]} + \mathfrak{a}_+\tilde{T}^{[n+1]} + \sum_{k=0}^{n}(\mathfrak{a}_-^{n-k}\tilde{T}^{[k]})\tilde{T}^{[n-k]},
\end{aligned}$$

and, as $\hat{T}^2 = H\tilde{T}^2H^{-1}$, we get the announced result. □

**Theorem 2** *For the recursion matrix given in* (14) *written in terms of its diagonals* $T = \Lambda + Ib + \Lambda^\top c + (\Lambda^\top)^2 d$, *with* $b$, $c$, $d$ *diagonal matrices from* (16), *we get*

$$\begin{aligned} d &= H(\mathfrak{a}_-^2 H)^{-1}, \\ c &= H(\mathfrak{a}_- H)^{-1}\big(\mathfrak{a}_-\tilde{T}^{[0]} + \tilde{T}^{[0]}\big), \\ b &= \tilde{T}^{[1]} + \mathfrak{a}_+\tilde{T}^{[1]} + (\tilde{T}^{[0]})^2, \\ (\mathfrak{a}_- H)H^{-1} &= \tilde{T}^{[2]} + \mathfrak{a}_+\tilde{T}^{[2]} + (\mathfrak{a}_-\tilde{T}^{[0]} + \tilde{T}^{[0]})\tilde{T}^{[1]}, \end{aligned}$$

*as well as*

$$\tilde{T}^{[n+1]} + \mathfrak{a}_+\tilde{T}^{[n+1]} + \sum_{k=0}^{n}(\mathfrak{a}_-^{n-k}\tilde{T}^{[k]})\tilde{T}^{[n-k]} = \mathbf{0}, \qquad n \in \{2, 3, \dots\}. \tag{21}$$

***Proof*** Use the previous two lemmas. □

**Corollary 1** *We have*

$$H_{2n+2} = \frac{H_0}{\prod_{k=0}^{n} d_{2(n-k)}}, \qquad H_{2n+3} = \frac{H_1}{\prod_{k=0}^{n} d_{2(n-k)+1}}.$$

***Proof*** It follows form $(\mathfrak{a}_-^2 H)d = H$. □

Finally, to end this section, we discuss on the Christoffel–Darboux (CD) kernels, see [40] for the standard orthogonality, and corresponding CD formulas in this multiple context. The two partial and complete CD kernels are given by

$$\begin{aligned} K_1^{(n)}(x, y) &:= \sum_{m=0}^{n-1} B^{(m)}(x)A_1^{(m)}(y), \\ K_2^{(n)}(x, y) &:= \sum_{m=0}^{n-1} B^{(m)}(x)A_2^{(m)}(y), \quad K^{(n)}(x, y) := \sum_{m=0}^{n-1} B^{(m)}(x)Q^{(m)}(y). \end{aligned} \tag{22}$$

These kernels satisfy Christoffel–Darboux formulas. Sorokin and Van Iseghem [41] derived a CD formula that can be applied to multiple orthogonal polynomials, see also [20, 22]. In Theorem 1 of [9] a CD formula for the mixed case was proven. Daems–Kuijlaars' CD formula, that is not in sequence, was derived in [23, 24]. The extension to partials kernels follow the ideas of [9]. The CD formula reads as follows

$$\begin{aligned} (x - y)K^{(n)}(x, y) = B^{(n)}(x)Q^{(n-1)}(y) & \\ -\big(T_{n,n-2}B^{(n-2)}(x) + T_{n,n-1}B^{(n-1)}(x)\big)Q^{(n)}(y) & \\ -B^{(n-1)}(x)Q^{(n+1)}(y)T_{n+1,n-1}, & \end{aligned} \tag{23}$$

as well as the partial CD formulas, for $a \in \{1, 2\}$, are

$$\begin{aligned}(x-y)K_a^{(n)}(x,y) = {} & B^{(n)}(x)A_a^{(n-1)}(y) \\ & -\big(T_{n,n-2}B^{(n-2)}(x) + T_{n,n-1}B^{(n-1)}(x)\big)A_a^{(n)}(y) \\ & -B^{(n-1)}(x)A_a^{(n+1)}(y)T_{n+1,n-1}. \end{aligned} \tag{24}$$

### 1.3 Hidden symmetry

In this section we explore further the connection of the system of weights $(w_1, w_2, \mathrm{d}\,\mu)$ and the recursion matrix and its sequences of multiple orthogonal polynomials of type II and linear forms of type I.

In particular we discuss on a hidden symmetry already remarked in [22]. For multiple orthogonal polynomials the recursion matrix $T$, the sequence of type II multiple orthogonal polynomials $\{B^{(n)}\}_{n=0}^{\infty}$, and even the sequence of type I linear forms $\{Q^{(n)}\}_{n=0}^{\infty}$ do not determine uniquely the spectral system $(w_1, w_2, \mathrm{d}\,\mu)$, as it determines uniquely $w_1$, for the second weight one has the hidden freedom described below.

**Proposition 6** *Let $\Delta \subset \mathbb{R}$ be the compact support of two perfect systems, $(w_1, w_2, \mathrm{d}\,\mu)$, with $\int_\Delta w_1(x)\,\mathrm{d}\,\mu(x) = 1$ and $(\hat{w}_1, \hat{w}_2, \mathrm{d}\,\mu)$, with $\int_\Delta \hat{w}_1(x)\,\mathrm{d}\,\mu(x) = 1$, that have the same sequence of type II multiple orthogonal polynomials $\{B^{(n)}(x)\}_{n=0}^{\infty}$. Then $\hat{w}_1 = w_1$ and there exists $\alpha, \beta, \gamma \in \mathbb{R}$, $\alpha, \beta \neq 0$ such that*

$$\gamma w_1 + \alpha w_2 + \beta \hat{w}_2 = 0. \tag{25}$$

*If $A_1^{(m)}$, $A_2^{(m)}$ are the type I multiple orthogonal polynomials, associated with the system $(w_1, w_2, \mathrm{d}\,\mu)$, then the type I multiple orthogonal polynomials, $\hat{A}_1^{(m)}$, $\hat{A}_2^{(m)}$ associated with the system $(w_1, \hat{w}_2, \mathrm{d}\,\mu)$ are given by*

$$\hat{A}_1^{(m)} = A_1^{(m)} - \frac{\gamma}{\alpha}A_2^{(m)}, \qquad \hat{A}_2^{(m)} = -\frac{\beta}{\alpha}A_2^{(m)},$$

*and both systems have the same linear associated forms*

$$\hat{Q}^{(m)} = Q^{(m)}, \quad Q^{(m)} := w_1 A_1^{(m)} + w_2 A_2^{(m)}, \quad \hat{Q}^{(m)} := w_1 \hat{A}_1^{(m)} + \hat{w}_2 \hat{A}_2^{(m)}.$$

***Proof*** For $B = SX$, the Gauss–Borel factorization (3) leads to

$$Sg = \int_\Delta SX(x)(X_1(x)w_1(x) + X_2(x)w_2(x))^\top \,\mathrm{d}\,\mu(x). \tag{26}$$

In terms of the projection matrices $\Pi_1$ and $\Pi_2$, we split (26) as follows

$$Sg\Pi_1 = \int_\Delta B(x)(X_1(x))^\top w_1(x)\,\mathrm{d}\,\mu(x), \quad Sg\Pi_2 = \int_\Delta B(x)(X_2(x))^\top w_2(x)\,\mathrm{d}\,\mu(x).$$

Assuming $\int_\Delta w_1(x)\,\mathrm{d}\,\mu(x) = 1$ and recalling that $Sg$ is an upper triangular matrix we get

$$\begin{bmatrix}1\\0\\\vdots\end{bmatrix} = \int_\Delta SX(x)w_1(x)\,\mathrm{d}\,\mu(x) \quad\Longrightarrow\quad S^{-1}\begin{bmatrix}1\\0\\\vdots\end{bmatrix} = \int_\Delta X(x)w_1(x)\,\mathrm{d}\,\mu(x),$$

and we conclude that given $S$, the moments $\int_\Delta x^n w_1\,\mathrm{d}\,\mu(x)$ are uniquely determined, and as the Hausdorff moment problem when solvable is determined, we obtain that $\hat{w}_1 = w_1$.

Similarly, we get

$$\begin{bmatrix}a\\b\\0\\\vdots\end{bmatrix} = \int_\Delta SX(x)w_2(x)\,\mathrm{d}\,\mu(x), \qquad \text{and} \qquad \begin{bmatrix}\hat{a}\\\hat{b}\\0\\\vdots\end{bmatrix} = \int_\Delta SX(x)\hat{w}_2(x)\,\mathrm{d}\,\mu(x).$$

We can assure that $b$ and $\hat{b}$ are different from zero, otherwise $(w_1, w_2, \mathrm{d}\,\mu)$ and $(w_1, \hat{w}_2, \mathrm{d}\,\mu)$ would not be perfect systems. Then, we can find $\alpha, \beta, \gamma$, where $\alpha, \beta \neq 0$ such that

$$\alpha\begin{pmatrix}a & b & 0\cdots\cdots\end{pmatrix} + \beta\begin{pmatrix}\hat{a} & \hat{b} & 0\cdots\cdots\end{pmatrix} + \gamma\begin{pmatrix}1 & 0 & 0\cdots\cdots\end{pmatrix} = \begin{pmatrix}0 & 0\cdots\cdots\end{pmatrix}.$$

Therefore,

$$\begin{bmatrix}0\\0\\\vdots\end{bmatrix} = \int_\Delta SX(x)(\alpha w_2 + \beta\hat{w}_2 + \gamma w_1)\,\mathrm{d}\,\mu(x)$$

which implies

$$\int_\Delta X(x)(\alpha w_2 + \beta\hat{w}_2 + \gamma w_1)\;\mathrm{d}\,\mu(x) = S^{-1}\begin{bmatrix}0\\0\\\vdots\end{bmatrix} = \begin{bmatrix}0\\0\\\vdots\end{bmatrix},$$

and as the Hausdorff moment problem is determined we get (25).

Now, we consider the sequence of type I linear form associated to the system $(w_1, w_2, \mathrm{d}\,\mu)$

$$Q^{(m)} := w_1 A_1^{(m)} + w_2 A_2^{(m)},$$

then it holds, that

$$Q^{(m)} := w_1 A_1^{(m)} + \Big(-\frac{\gamma}{\alpha}w_1 - \frac{\beta}{\alpha}\hat{w}_2\Big)A_2^{(m)} = w_1\Big(A_1^{(m)} - \frac{\gamma}{\alpha}A_2^{(m)}\Big) + \hat{w}_2\Big(-\frac{\beta}{\alpha}A_2^{(m)}\Big).$$

Using the biorthogonality of the linear form, that the degree of the polynomial $A_2^{(m)}$ is less than or equal to the degree of the polynomial $A_1^{(m)}$, and the uniqueness of the type $I$ polynomials associated to the perfect system $(w_1, \hat{w}_2, \mathrm{d}\,\mu)$ we recover that $\hat{A}_1^{(m)} = A_1^{(m)} - \frac{\gamma}{\alpha} A_2^{(m)}$ and $\hat{A}_2^{(m)} = -\frac{\beta}{\alpha} A_2^{(m)}$, and also that both systems have the same sequence of type I linear form associated $\hat{Q}^{(m)} = Q^{(m)}$ as we wanted to prove. □

## 2 Pearson equation and differential equations

Notice that the Jacobi–Piñeiro weights, $\tilde{w}_a(x) = w_a(x)(1-x)^\gamma$, $a \in \{1, 2\}$, i.e.

$$\tilde{w}_1 = x^\alpha (1-x)^\gamma, \qquad \tilde{w}_2 = x^\beta (1-x)^\gamma,$$

fulfill the following Pearson equations

$$x(1-x)\tilde{w}_1'(x) = (\alpha(1-x) - \gamma x)\tilde{w}_1(x), \quad x(1-x)\tilde{w}_2'(x) = (\beta(1-x) - \gamma x)\tilde{w}_2(x), \tag{27}$$

as well as

$$x^\beta \tilde{w}_1 = x^\alpha \tilde{w}_2.$$

Let us assume that $\mathrm{d}\,\mu = v(x)\,\mathrm{d}\,x$, denote $\tilde{w}_a = w_a v$ with $v(x) = (1-x)^\gamma$, and consider the following Pearson relations

$$\sigma(x)\tilde{w}_1'(x) = q_1(x)\tilde{w}_1(x), \qquad \sigma(x)\tilde{w}_2'(x) = q_2(x)\tilde{w}_2(x), \tag{28}$$

where for the Jacobi–Piñeiro case we have $\sigma(x) = x(1-x)$, $q_1(x) = (\alpha(1-x) - \gamma x)$, $q_2(x) = (\beta(1-x) - \gamma x)$, with $\sigma$ a polynomial such that $\sigma\tilde{w}_1 = \sigma\tilde{w}_2 = 0$ at the boundary $\partial\Delta$ of the support. If we introduce the notation

$$\tau_1(x) := q_1(x) + \sigma'(x), \qquad \tau_2(x) := q_2(x) + \sigma'(x),$$

the Pearson equations are

$$(\sigma\tilde{w}_1)' = \tau_1\tilde{w}_1, \qquad (\sigma\tilde{w}_2)' = \tau_2\tilde{w}_2.$$

On the other hand, to handle derivatives of type II multiple orthogonal polynomials we introduce the matrices

$$\mathscr{D} := \begin{bmatrix} 0 & 0 & 0 & \cdots & \\ 1 & 0 & 0 & \ddots & \\ 0 & 2 & 0 & \ddots & \ddots \\ 0 & 0 & 3 & \ddots & \ddots \\ \vdots & \ddots & \ddots & \ddots & \ddots \end{bmatrix} = \Lambda^\top \mathscr{N}, \qquad \mathscr{N} := \operatorname{diag}(1, 2, 3, \dots).$$

**Proposition 7** *The derivative of the vector of monomials is*

$$X' = \mathscr{D} X.$$

*The following identity*

$$[\Lambda, \mathscr{D}] = I$$

*is satisfied.*

***Proof*** This relation can be checked directly, from the explicit for of the matrices or by observing that

$$[\Lambda, \mathscr{D}]X = (\Lambda\mathscr{D} - \mathscr{D}\Lambda)X = \Lambda X' - x X' = \begin{bmatrix} 1 \\ 2x \\ 3x^2 \\ \vdots \end{bmatrix} - x \begin{bmatrix} 0 \\ 1 \\ 2x \\ \vdots \end{bmatrix} = \begin{bmatrix} 1 \\ x \\ x^2 \\ \vdots \end{bmatrix},$$

and the result follows. □

On the other hand, to deal with derivatives of type I multiple orthogonal polynomials, we require of the matrices

$$\mathscr{D}_1 := \begin{bmatrix} 0 & 0 & 0 & 0 & 0 & \cdots \\ 0 & 0 & 0 & 0 & 0 & \ddots \\ 1 & 0 & 0 & 0 & 0 & \ddots \\ 0 & 0 & 0 & 0 & 0 & \ddots \\ 0 & 0 & 2 & 0 & 0 & \ddots \\ 0 & 0 & 0 & 0 & 0 & \ddots \\ 0 & 0 & 0 & 0 & 3 & \ddots \\ \vdots & \ddots & \ddots & \ddots & \ddots & \ddots \end{bmatrix}, \qquad \mathscr{D}_2 := \begin{bmatrix} 0 & 0 & 0 & 0 & 0 & \cdots \\ 0 & 0 & 0 & 0 & 0 & \ddots \\ 0 & 0 & 0 & 0 & 0 & \ddots \\ 0 & 1 & 0 & 0 & 0 & \ddots \\ 0 & 0 & 0 & 0 & 0 & \ddots \\ 0 & 0 & 0 & 2 & 0 & \ddots \\ 0 & 0 & 0 & 0 & 0 & \ddots \\ \vdots & \ddots & \ddots & \ddots & \ddots & \ddots \end{bmatrix}.$$

**Lemma 5** *We have*

$$X_1' = \mathscr{D}_1 X_1, \qquad \mathbf{0} = \mathscr{D}_1 X_2, \qquad \mathbf{0} = \mathscr{D}_2 X_1, \qquad X_2' = \mathscr{D}_2 X_2.$$

**Definition 6** Let us introduce

$$\begin{aligned}\mathscr{N}_1 &:= \operatorname{diag}(1,0,2,0,3,\dots), \quad \mathscr{N}_2 := \operatorname{diag}(0,1,0,2,0,3,\dots)),\\ \left\lceil\frac{\mathscr{N}}{2}\right\rceil &:= \mathscr{N}_1+\mathscr{N}_2 = \operatorname{diag}(1,1,2,2,3,3,\dots).\end{aligned}$$

**Lemma 6** *The following relations are fulfilled*

$$\mathscr{D}_1 = \Lambda_1^\top \mathscr{N}_1 = (\Lambda^\top)^2\mathscr{N}_1, \quad \mathscr{D}_2 = \Lambda_2^\top \mathscr{N}_1 = (\Lambda^\top)^2\mathscr{N}_2, \quad \mathscr{D}_1+\mathscr{D}_2 = (\Lambda^\top)^2\left\lceil\frac{\mathscr{N}}{2}\right\rceil. \tag{29}$$

*We also have*

$$[\Lambda_1,\mathscr{D}_1] = \Pi_1, \qquad [\Lambda_2,\mathscr{D}_2] = \Pi_2.$$

**Theorem 3** *If the Pearson equations* (28) *hold, then the following symmetry for the moment matrix*

$$\sigma(\Lambda)\mathscr{D}g = -g\big(\mathscr{D}_1\sigma(\Lambda_1)+q_1(\Lambda_1)+\mathscr{D}_2\sigma(\Lambda_2)+q_2(\Lambda_2)\big)^\top \tag{30}$$

*is fulfilled.*

***Proof*** Observing

$$\big(\sigma(x)X(x)\big)' = \big(\sigma(\Lambda)X(x)\big)' = \sigma(\Lambda)X'(x) = \sigma(\Lambda)\mathscr{D}X(x)$$

we deduce that

$$\begin{aligned}\sigma(\Lambda)\mathscr{D}g &= \int_\Delta \sigma(\Lambda)\mathscr{D}X(x)\big(\tilde w_1(x)X_1(x)+\tilde w_2(x)X_2(x)\big)^\top \mathrm{d}\,x\\ &= \int_\Delta \big(\sigma(x)X(x)\big)'\big(\tilde w_1(x)X_1(x)+\tilde w_2(x)X_2(x)\big)^\top \mathrm{d}\,x.\end{aligned}$$

Integrating by parts we get

$$\begin{aligned}\sigma(\Lambda)\mathscr{D}g =& \big(\sigma(x)X(x)(\tilde w_1(x)X_1(x)+\tilde w_2X_2(x))^\top\big)|_{\partial\Delta}\\ &-\int_\Delta \sigma(x)X(x)(\tilde w_1(x)X_1^\top(x)+\tilde w_2X_2^\top(x))'\,\mathrm{d}\,x.\end{aligned}$$

Hence, recalling that $\sigma(x)\tilde{w}_1(x) = 0$ and $\sigma(x)\tilde{w}_2(x) = 0$ at the boundaries of the support, we obtain

$$\begin{aligned}\sigma(\Lambda)\mathscr{D}g &= -\int_\Delta \sigma(x)X(x)\big(\tilde{w}_1(x)X_1^\top(x) + \tilde{w}_2(x)X_2^\top(x)\big)'\,\mathrm{d}\,x\\ &= -\int_\Delta \sigma(x)X(x)\big(\tilde{w}_1'(x)X_1^\top(x) + \tilde{w}_2'(x)X_2^\top(x)\big)\,\mathrm{d}\,x\\ &\quad - \int_\Delta \sigma(x)X(x)\big(\tilde{w}_1(x)(X_1'(x))^\top + \tilde{w}_2(x)(X_2'(x))^\top\big)\,\mathrm{d}\,x.\end{aligned}$$

Then, the Pearson equations (27) lead to

$$\begin{aligned}\sigma(\Lambda)\mathscr{D}g &= -\int_\Delta X(x)\big(q_1(x)\tilde{w}_1(x)X_1^\top(x) + q_2(x)\tilde{w}_2(x)X_2^\top(x)\big)\,\mathrm{d}\,x\\ &\quad - \int_\Delta \sigma(x)X(x)\big(\tilde{w}_1(x)(X_1'(x))^\top + \tilde{w}_2(x)(X_2'(x))^\top\big)\,\mathrm{d}\,x\\ &= -\Big(\int_\Delta X(x)\big(\tilde{w}_1(x)X_1^\top(x) + \tilde{w}_2(x)X_2^\top(x)\big)\,\mathrm{d}\,x\Big)(q_1(\Lambda_1) + q_2(\Lambda_2))^\top\\ &\quad - \Big(\int_\Delta \sigma(x)X(x)\big(\tilde{w}_1(x)(X_1(x))^\top + \tilde{w}_2(x)(X_2(x))^\top\big)\,\mathrm{d}\,x\Big)(\mathscr{D}_1 + \mathscr{D}_2)^\top,\end{aligned}$$

and we conclude that

$$\sigma(\Lambda)\big(\mathscr{D}g + g(\mathscr{D}_1 + \mathscr{D}_2)^\top\big) = -g(q_1(\Lambda_1) + q_2(\Lambda_2))^\top.$$

Finally, recalling that

$$\Lambda g = g(\Lambda_1 + \Lambda_2)^\top, \qquad \text{so that} \qquad \sigma(\Lambda)g = g(\sigma(\Lambda_1) + \sigma(\Lambda_2))^\top,$$

and that

$$(\mathscr{D}_1 + \mathscr{D}_2)(\sigma(\Lambda_1) + \sigma(\Lambda_2)) = \mathscr{D}_1\sigma(\Lambda_1) + \mathscr{D}_2\sigma(\Lambda_2)$$

we obtain (30). □

**Definition 7** Let us define

$$\Phi := S\mathscr{D}S^{-1}, \qquad \Phi_1 := H^{-1}\tilde{S}\mathscr{D}_1\tilde{S}^{-1}H, \qquad \Phi_2 := H^{-1}\tilde{S}\mathscr{D}_2\tilde{S}^{-1}H.$$

***Remark 3*** (i) The matrices $\Phi$, $\Phi_1$ and $\Phi_2$ model the derivatives of the multiple orthogonal polynomials as follows

$$B' = \Phi B, \qquad A_1' = \Phi_1 A_1, \qquad A_2' = \Phi_2 A_2.$$

Moreover, we have

$$\Phi_1 A_2 = \mathbf{0}, \qquad \Phi_2 A_1 = \mathbf{0}.$$

(ii) They also satisfy

$$\Phi_1 T_2 = T_2 \Phi_1 = \Phi_2 T_1 = T_1 \Phi_2 = \mathbf{0}.$$

(iii) The matrices $\Phi$, $\Phi_1$ and $\Phi_2$ are strictly lower triangular. The first possibly non zero subdiagonal of $\Phi$ is the first one and of $\Phi_1$ and $\Phi_2$ the second one.

(iv) We also have

$$[T, \Phi] = I.$$

(v) Introducing the lower triangular matrices $C_1 = H^{-1}\tilde{S}\Pi_1\tilde{S}^{-1}H$ and $C_2 = H^{-1}\tilde{S}\Pi_2\tilde{S}^{-1}H$, that are projections $C_1^2 = C_1$, $C_2^2 = C_2$, $C_1C_2 = C_2C_1$ and $C_1 + C_2 = I$, we have

$$[T_1, \Phi_1] = C_1, \qquad [T_1, \Phi_2] = C_2, \qquad [T_1 + T_2, \Phi_1 + \Phi_2] = I.$$

Now, using the Gauss–Borel factorization, we are ready to express the symmetry for the moment matrix (30) in the following terms

**Theorem 4** *We have*

$$\sigma(T)\Phi = -\big(\Phi_1\sigma(T_1) + \Phi_2\sigma(T_2) + q_1(T_1) + q_2(T_2)\big)^\top. \tag{31}$$

***Proof*** Using the Gauss–Borel factorization $g = S^{-1}H\tilde{S}^{-\top}$ of the moment matrix $g$ we can write (30) as follows

$$\sigma(\Lambda)\mathscr{D}S^{-1}H\tilde{S}^{-\top} = -S^{-1}H\tilde{S}^{-\top}(\mathscr{D}_1\sigma(\Lambda_1) + q_1(\Lambda_1) + \mathscr{D}_2\sigma(\Lambda_2) + q_2(\Lambda_2))^\top,$$

so that

$$S\sigma(\Lambda)S^{-1}S\mathscr{D}S^{-1}H = -S\sigma(\Lambda)S^{-1}H\tilde{S}^{-\top}(\mathscr{D}_1 + \mathscr{D}_2)^\top\tilde{S}^\top - H\tilde{S}^{-\top}(q_1(\Lambda_1) + q_2(\Lambda_2))^\top\tilde{S}^\top,$$

from where we get (31). □

**Definition 8** We introduce the Laguerre–Freud matrix

$$\Psi := \sigma(T)\Phi = -\big(\Phi_1\sigma(T_1) + \Phi_2\sigma(T_2) + q_1(T_1) + q_2(T_2)\big)^\top.$$

**Proposition 8** *The Laguerre–Freud matrix $\Psi$ is a banded matrix with* $\deg\sigma - 1$ *possible non zero superdiagonals and* $2\max(\deg\sigma - 1, \deg q_1, \deg q_2)$ *possibly non zero subdiagonals.*

***Proof*** This follows from (31). Indeed, the highest non zero superdiagonal of $\Psi = \sigma(T)\Phi$ is the $(\deg\sigma-1)$-th one, and the lowest subdiagonal of $\Psi = -\big(\Phi_1\sigma(T_1)+\Phi_2\sigma(T_2)+q_1(T_1)+q_2(T_2)\big)^\top$ is $2\max(\deg\sigma-1,\deg q_1,\deg q_2)$-th one. □

**Theorem 5** *In terms of the Laguerre–Freud matrix we have:*

*(i) The type I multiple orthogonal polynomials fulfill*

$$\begin{aligned}\sigma(x)A_1'(x)+q_1(x)A_1(x)+\Psi^\top A_1(x)&=\mathbf{0},\\ \sigma(x)A_2'(x)+q_2(x)A_2(x)+\Psi^\top A_2(x)&=\mathbf{0}.\end{aligned}\tag{32}$$

*(ii) The linear form* $\tilde{Q}(x)=A_1(x)\tilde{w}_1(x)+A_2(x)\tilde{w}_2(x)=Q(x)v(x)$ *is subject to*

$$\sigma\tilde{Q}'(x)+\Psi^\top\tilde{Q}(x)=\mathbf{0}.\tag{33}$$

*(iii) The type II orthogonal polynomials satisfy*

$$(\sigma(x)B(x))'=\Psi B(x).\tag{34}$$

***Proof*** From (31) we deduce that $A_1$ fulfills

$$-(\sigma(T)\Phi)^\top A_1=\Phi_1\sigma(T_1)A_1+q_1(T_1)A_1=\sigma(x)A_1'(x)+q_1(x)A_1(x).$$

For $A_2$ we proceed analogously. For the linear form we have

$$\begin{aligned}\sigma\tilde{Q}'&=\sigma(A_1\tilde{w}_1+A_2\tilde{w}_2)'=\sigma(A_1'\tilde{w}_1+A_1\tilde{w}_1'+A_2'\tilde{w}_2+A_2\tilde{w}_2')\\&=(\sigma A_1'+q_1A_1)\tilde{w}_1+(\sigma A_2'+q_2A_2)\tilde{w}_2,\end{aligned}$$

where we have used the Pearson equations (28). Then, from (32) the relation (33) for the linear form $\tilde{Q}$ follows. Moreover, for the multiple orthogonal polynomials of type II we get that $B$ is subject to

$$\sigma(T)\Phi B=\sigma(T)B'=(\sigma(T)B))'=(\sigma(x)B(x))',$$

which completes the proof. □

In terms of the determinants

$$\delta^{(n,m)}:=\begin{vmatrix}B^{(n)}(a)&B^{(n)}(b)\\B^{(m)}(a)&B^{(m)}(b)\end{vmatrix},$$

we have:

**Proposition 9** *For a support* $\Delta=[a,b]$*, a quadratic polynomial* $\sigma(x)=-(x-a)(x-b)$*,* $\deg q_1=\deg q_2=1$*, the Laguerre–Freud matrix has the following nonzero coefficients*

$$\Psi_{n,n+1}=-n-2,$$

$$\Psi_{n+2,n} = \left\lfloor \frac{n}{2} \right\rfloor,$$
$$\Psi_{n+2,n+1} = -\frac{-2(b-a)B^{(n+2)}(a)B^{(n+2)}(b) + \lfloor \frac{n}{2} \rfloor \delta^{(n,n+2)} - (n+4)\delta^{(n+3,n+2)}}{\delta^{(n+1,n+2)}},$$
$$\Psi_{n+2,n+2} = -\frac{(b-a)(B^{(n+2)}(a)B^{(n+1)}(b) + B^{(n+1)}(a)B^{(n+2)}(b)) + \lfloor \frac{n}{2} \rfloor \delta^{(n+1,n)} - (n+4)\delta^{(n+1,n+3)}}{\delta^{(n+1,n+2)}}.$$

***Proof*** From $\Psi = \sigma(T)\Phi$ we have that the first superdiagonal of $\Psi$ is

$$-\Lambda^2\mathscr{D} = -\Lambda^2\Lambda^\top\mathscr{N} = -\Lambda\mathscr{N} = -(\mathfrak{a}_-\mathscr{N})\Lambda$$

and from $\Psi = -\big(\Phi_1\sigma(T_1) + \Phi_2\sigma(T_2) + q_1(T_1) + q_2(T_2)\big)^\top$ using (29) we get that the second subdiagonal is the transpose of

$$-\mathscr{D}_1\Lambda_1^2 - \mathscr{D}_2\Lambda_2^2 = -(\mathscr{D}_1 + \mathscr{D}_2)\Lambda^4 = (\Lambda^\top)^2\left\lceil \frac{\mathscr{N}}{2} \right\rceil \Lambda^4 = \left(\mathfrak{a}_+^2\left\lceil \frac{\mathscr{N}}{2} \right\rceil\right)\Lambda^2$$
$$= \left\lceil \frac{\mathscr{N}}{2} - I \right\rceil \Lambda^2 = \left\lfloor \frac{\mathscr{N} - I}{2} \right\rfloor \Lambda^2.$$

From (34) we deduce that

$$\sigma'(x)B(x) + \sigma(x)B'(x) = \Psi B(x)$$

so that, as $a, b$ are zeros of $\sigma$

$$\sigma'(a)B(a) = \Psi B(a), \qquad \sigma'(b)B(b) = \Psi B(b).$$

Thus, as $\sigma'(a) = -\sigma'(b) = (b-a)$, we get

$$(b-a)B(a) = \Psi B(a), \qquad -(b-a)B(b) = \Psi B(b).$$

These two equations can be written, for $n \geq 2$ as

$$(\Psi_{n+2,n+1}, \Psi_{n+2,n+2}) \begin{bmatrix} B^{(n+1)}(a) & B^{(n+1)}(b) \\ B^{(n+2)}(a) & B^{(n+2)}(b) \end{bmatrix}$$
$$= (b-a)(B^{(n+2)}(a), -B^{(n+2)}(b)) - \left(\left\lfloor \frac{n}{2} \right\rfloor, -(n+4)\right) \begin{bmatrix} B^{(n)}(a) & B^{(n)}(b) \\ B^{(n+3)}(a) & B^{(n+3)}(b) \end{bmatrix}.$$

Therefore, we get

$$
\begin{aligned}
(\Psi_{n+2,n+1},\Psi_{n+2,n+2}) &= \Bigg((b-a)(B^{(n+2)}(a),-B^{(n+2)}(b))\\
&\quad -\left(\left\lfloor\frac{n}{2}\right\rfloor,-(n+4)\right)\begin{bmatrix} B^{(n)}(a) & B^{(n)}(b)\\ B^{(n+3)}(a) & B^{(n+3)}(b)\end{bmatrix}\Bigg)\begin{bmatrix} B^{(n+1)}(a) & B^{(n+1)}(b)\\ B^{(n+2)}(a) & B^{(n+2)}(b)\end{bmatrix}^{-1}\\
&= \frac{1}{\delta^{(n+1,n+2)}}\Bigg((b-a)\big(B^{(n+2)}(a),-B^{(n+2)}(b)\big)\\
&\quad -\left(\left\lfloor\frac{n}{2}\right\rfloor,-(n+4)\right)\begin{bmatrix} B^{(n)}(a) & B^{(n)}(b)\\ B^{(n+3)}(a) & B^{(n+3)}(b)\end{bmatrix}\Bigg)\begin{bmatrix} B^{(n+2)}(b) & -B^{(n+1)}(b)\\ -B^{(n+2)}(a) & B^{(n+1)}(a)\end{bmatrix}\\
&= \frac{1}{\delta^{(n+1,n+2)}}\Bigg((b-a)\big(2B^{(n+2)}(a)B^{(n+2)}(b),-B^{(n+2)}(a)B^{(n+1)}(b)\\
&\quad -B^{(n+1)}(a)B^{(n+2)}(b)\big)-\left(\left\lfloor\frac{n}{2}\right\rfloor,-(n+4)\right)\begin{bmatrix} \delta^{(n,n+2)} & \delta^{(n+1,n)}\\ \delta^{(n+3,n+2)} & \delta^{(n+1,n+3)}\end{bmatrix}\Bigg),
\end{aligned}
$$

and the result is proven. □

***Remark 4*** Therefore, the Laguerre–Freud matrix has the following structure

$$
\Psi=\begin{bmatrix}
\Psi_{0,0} & -2 & 0 & \cdots & & & & & \\
\Psi_{1,0} & \Psi_{1,1} & -3 & \ddots & & & & & \\
0 & \Psi_{2,1} & \Psi_{2,2} & -4 & \ddots & & & & \\
0 & 0 & \Psi_{3,2} & \Psi_{3,3} & -5 & \ddots & & & \\
\vdots & \ddots & 1 & \Psi_{4,3} & \Psi_{4,4} & -6 & \ddots & & \\
\vdots & & \ddots & 1 & \Psi_{5,4} & \Psi_{5,5} & -7 & \ddots & \\
\vdots & & & \ddots & 2 & \Psi_{6,5} & \Psi_{6,6} & -8 & \ddots \\
\vdots & & & & 2 & \Psi_{7,6} & \Psi_{7,7} & -9 & \ddots \\
\vdots & & & & & \ddots & \ddots & \ddots & \ddots
\end{bmatrix}.
$$

***Remark 5*** Componentwise Eq. (33) for the linear forms $\{\tilde{Q}^{(k)}(x)\}_{k=0}^{\infty}$ is

$$
\begin{aligned}
\sigma\frac{\mathrm{d}\,\tilde{Q}^{(n+1)}}{\mathrm{d}\,x}-(n+2)\tilde{Q}^{(n)}+\Psi_{n+1,n+1}\tilde{Q}^{(n+1)}+\Psi_{n+2,n+1}\tilde{Q}^{(n+2)}\\
+\left\lfloor\frac{n+1}{2}\right\rfloor\tilde{Q}^{(n+3)}=0.
\end{aligned}
$$

Equation (34) for the type II orthogonal polynomials $\{B^{(k)}(x)\}_{k=0}^{\infty}$ reads componentwise as follows

$$
\frac{\mathrm{d}\,\sigma B^{(n+2)}}{\mathrm{d}\,x}=\left\lfloor\frac{n}{2}\right\rfloor B^{(n)}+\Psi_{n+2,n+1}B^{(n+1)}+\Psi_{n+2,n+2}B^{(n+2)}-(n+4)B^{(n+3)}.
$$

***Remark 6*** The Jacobi–Piñeiro case corresponds to the particular choice $a = 0$ and $b = 1$ and $\sigma = x(1-x)$. Moreover, according to [17]

$$
\begin{aligned}
B_{(n,m)}(1) &= \frac{(\gamma+1)_{m+n}}{(\alpha+\gamma+m+n+1)_n(\beta+\gamma+m+n+1)_m},\\
B_{(n,m)}(0) &= \frac{(-1)^{m+n}(\alpha+1)_n(\beta+1)_m}{(\alpha+\gamma+m+n+1)_n(\beta+\gamma+m+n+1)_m},
\end{aligned}
$$

and, consequently,the determinants involved are

$$
\begin{aligned}
&\delta^{((n,m),(k,l))}\\
&= \frac{(-1)^{m+n}(\alpha+1)_n(\beta+1)_m(\gamma+1)_{k+l} - (-1)^{k+l}(\alpha+1)_k(\beta+1)_l(\gamma+1)_{m+n}}{(\alpha+\gamma+m+n+1)_n(\beta+\gamma+m+n+1)_m(\alpha+\gamma+k+l+1)_k(\beta+\gamma+k+l+1)_l}.
\end{aligned}
$$

***Remark 7*** In [5, Corollary 2] the authors prove that the type II, multiple orthogonal polynomials that admits a Rodrigues formula representation, satisfies a $p+1$ order differential equations of type

$$
\sum_{k=0}^{p+1} q_{k,p+1}\big(B^{(n_1,n_2)}\big)^{(p+1-k)} = 0
$$

with explicitly given polynomial coefficients $q_{k,p+1}$ that could depend on $n$. In the case of Jacobi–Piñeiro, and for a system of two weight functions $\tilde{w}_1, \tilde{w}_2$ we have, as $p = 2$, a third order differential equation with polynomial coefficients given by,

$$
\begin{aligned}
q_{0,3} &= z^2(1-z)^2,\\
q_{1,3} &= z(z-1)\big(2(\gamma+1)z + (\alpha+\beta+3)(z-1)\big),\\
q_{2,3} &= \gamma(\gamma+1)z^2 + (\alpha+1)(\beta+1)(z-1)^2\\
&\quad - z(z-1)\big((\gamma+\alpha+n_1+1)n_1 + (\gamma+\beta+n_2+1)n_2\\
&\quad + n_1n_2 - (\gamma+1)(\alpha+\beta+3)\big),\\
q_{3,3} &= \gamma z\big((\gamma+\alpha+n_1+1)n_1 + (\gamma+\beta+n_2+1)n_2\big)\\
&\quad - (z-1)\big(n_1(\gamma+\alpha+n_1+1)(\beta+n_2+1)\\
&\quad + n_2(\gamma+\beta+n_2+1)(\alpha+n_1+1)\big) - \gamma n_1n_2.
\end{aligned}
$$

In the proof the author's heavily used the Rodrigues type formula representation for these polynomials.

We think that this technique cannot be used to derive the third order differential equations for the vector of type I linear forms. Nevertheless, by duality of the one derived in [5] for the type II multiple orthogonal polynomials of Jacobi–Piñeiro, we think that this can be achieved.

# 3 Christoffel and Geronimus perturbations: Christoffel formulas

We will start by considering basic simple cases and then we move to more general polynomial perturbations.

**Definition 9** *(Permuting Christoffel transformation)* Let us consider $\vec{w} = (w_1, w_2)$ and the transformed vector of weights $\underline{\vec{w}} = (w_2, x\, w_1)$, that is a simple Christoffel transformation of $w_1$ followed by a permutation of the two weights.

**Proposition 10** (Permuting Christoffel transformation and the moment matrix) *The moment matrix satisfies*

$$g_{\vec{w}}\, \Lambda^\top = g_{\underline{\vec{w}}}. \tag{35}$$

***Proof*** In fact, we successively have

$$\begin{aligned} g_{(w_1,w_2)}\, \Lambda^\top &= \int_\Delta X(x)\big(X_1(x)w_1(x) + X_2(x)w_2(x)\big)^\top \,\mathrm{d}\,\mu(x)\Lambda^\top \\ &= \int_\Delta X(x)\big(\Lambda\, X_1(x)w_1(x) + \Lambda\, X_2(x)w_2(x)\big)^\top \,\mathrm{d}\,\mu(x) \\ &= \int_\Delta X(x)\big(x X_2(x)w_1(x) + X_1(x)w_2(x)\big)^\top \,\mathrm{d}\,\mu(x) \qquad \text{by (19)} \\ &= g_{(w_2, x\, w_1)}, \end{aligned}$$

which completes the proof. □

***Remark 8*** Observe that iterating (35) we get

$$g_{(w_1,w_2)}\, (\Lambda^\top)^2 = g_{(w_2,x\, w_1)}\, \Lambda^\top = g_{(x\, w_1, x\, w_2)} = \Lambda\, g_{(w_1,w_2)}, \tag{36}$$

and bi-Hankel property (13) of the moment matrix is recovered. In this sense the transformation $\vec{w} \to \underline{\vec{w}}$ can be understood as a square root of the Christoffel transformation $\vec{w} \to x\vec{w}$.

Let us assume that $g_{\underline{\vec{w}}}$ has a Gauss–Borel factorization

$$g_{\underline{\vec{w}}} = S_{\underline{\vec{w}}}^{-1} H_{\underline{\vec{w}}}\, \tilde{S}_{\underline{\vec{w}}}^{-\top}.$$

**Definition 10** We introduce the connection matrix

$$\Omega := S_{\vec{w}}\, S_{\underline{\vec{w}}}^{-1}.$$

Then,

**Theorem 6** (Permuting Christoffel transformation) *The connection matrix can be written as follows*

$$\Omega = \begin{bmatrix} 1 & 0 & \cdots & \cdots \\ \frac{H_{\vec{w},1}}{H_{\underline{\vec{w}},0}} & 1 & \ddots & \\ 0 & \frac{H_{\vec{w},2}}{H_{\underline{\vec{w}},1}} & 1 & \ddots \\ \vdots & \ddots & \ddots & \ddots \end{bmatrix}.$$

*The following connection formulas hold true*

$$\Omega B_{\underline{\vec{w}}} = B_{\vec{w}}, \tag{37}$$

$$\Omega^\top A_{1,\vec{w}} = x A_{2,\underline{\vec{w}}}, \quad \Omega^\top A_{2,\vec{w}} = A_{1,\underline{\vec{w}}}, \tag{38}$$

$$\Omega^\top Q_{\vec{w}} = Q_{\underline{\vec{w}}}. \tag{39}$$

***Proof*** The Gauss–Borel factorization of (35) leads to

$$S_{\vec{w}}^{-1} H_{\vec{w}} \tilde{S}_{\vec{w}}^{-\top} \Lambda^\top = S_{\underline{\vec{w}}}^{-1} H_{\underline{\vec{w}}} \tilde{S}_{\underline{\vec{w}}}^{-\top},$$

so that

$$\Omega = S_{\vec{w}} S_{\underline{\vec{w}}}^{-1} = \left(H_{\underline{\vec{w}}}^{-1} \tilde{S}_{\underline{\vec{w}}} \Lambda \tilde{S}_{\vec{w}}^{-1} H_{\vec{w}}\right)^\top.$$

Thus, we deduce that the matrix $\Omega$ is an unitriangular matrix with only its first subdiagonal different from zero as well as the corresponding subdiagonal coefficients.

Moreover, from definition we get

$$\Omega B_{\underline{\vec{w}}} = S_{\vec{w}} S_{\underline{\vec{w}}}^{-1} S_{\underline{\vec{w}}} X = S_{\vec{w}} X = B_{\vec{w}}.$$

Now,

$$\Omega^\top A_{1,\vec{w}} = H_{\underline{\vec{w}}}^{-1} \tilde{S}_{\underline{\vec{w}}} \Lambda \tilde{S}_{\vec{w}}^{-1} H_{\vec{w}} H_{\vec{w}}^{-1} \tilde{S}_{\vec{w}} X_1 = H_{\underline{\vec{w}}}^{-1} \tilde{S}_{\underline{\vec{w}}} \Lambda X_1 = x H_{\underline{\vec{w}}}^{-1} \tilde{S}_{\underline{\vec{w}}} X_2 = x A_{2,\underline{\vec{w}}},$$

$$\Omega^\top A_{2,\vec{w}} = H_{\underline{\vec{w}}}^{-1} \tilde{S}_{\underline{\vec{w}}} \Lambda \tilde{S}_{\vec{w}}^{-1} H_{\vec{w}} H_{\vec{w}}^{-1} \tilde{S}_{\vec{w}} X_2 = H_{\underline{\vec{w}}}^{-1} \tilde{S}_{\underline{\vec{w}}} \Lambda X_2 = H_{\underline{\vec{w}}}^{-1} \tilde{S}_{\underline{\vec{w}}} X_1 = A_{1,\underline{\vec{w}}}.$$

Finally,

$$\Omega^\top Q_{\vec{w}} = \Omega^\top (A_{1,\vec{w}} w_1 + A_{2,\vec{w}} w_2) = A_{2,\underline{\vec{w}}}(x w_1) + A_{1,\underline{\vec{w}}} w_2 = Q_{\underline{\vec{w}}},$$

which ends the proof. □

We now discuss similar connection formulas but for the CD kernels.

**Proposition 11** *For $n \in \mathbb{N}$, the CD kernels (22) satisfy*

$$
\begin{aligned}
K^{(n)}_{\underline{\vec{w}},1}(x,y) &= K^{(n)}_{\vec{w},2}(x,y) + B^{(n-1)}_{\underline{\vec{w}}}(x)\Omega_{n,n-1}A^{(n)}_{\vec{w},2}(y),\\
yK^{(n)}_{\underline{\vec{w}},2}(x,y) &= K^{(n)}_{\vec{w},1}(x,y) + B^{(n-1)}_{\underline{\vec{w}}}(x)\Omega_{n,n-1}A^{(n)}_{\vec{w},1}(y),\\
K^{(n)}_{\underline{\vec{w}}}(x,y) &= K^{(n)}_{\vec{w}}(x,y) + B^{(n-1)}_{\underline{\vec{w}}}(x)\Omega_{n,n-1}Q^{(n)}_{\vec{w}}(y).
\end{aligned} \tag{40}
$$

***Proof*** Given any semi-infinite vector $C$ or matrix $M$ we will denote by $C_{[n]}$ or matrix $M_{[n]}$ the truncations, where we keep only the first $n$ rows or $n$ rows and columns (the indices will run from 0 up to $n-1$), respectively. Then, given the band form of $\Omega$ we find

$$
(\Omega^\top C)_{[n]} = \Omega^\top_{[n]}C_{[n]} + \begin{bmatrix} 0 \\ \vdots \\ 0 \\ \Omega_{n,n-1}C^{(n)} \end{bmatrix}, \tag{41}
$$

with $C^{(n)}$ the corresponding entry of the vector $C$.

For $a \in \{1,2\}$, we can write (22) using he following vector notation:

$$
K^{(n)}_a(x,y) := (B(x))^\top_{[n]}(A(y))_{[n]}, \qquad K^{(n)}(x,y) := (B(x))^\top_{[n]}(Q(x))_{[n]}. \tag{42}
$$

Then, using (37), (38) and (39) we deduce

$$
\begin{aligned}
K^{(n)}_{\underline{\vec{w}},1}(x,y) &= (B_{\underline{\vec{w}}}(x))^\top_{[n]}(A_{\underline{\vec{w}},1})_{[n]}(y) = (B_{\underline{\vec{w}}}(x))^\top_{[n]}(\Omega^\top A_{\vec{w},2}(x))_{[n]},\\
yK^{(n)}_{\underline{\vec{w}},2}(x,y) &= y(B_{\underline{\vec{w}}}(x))^\top_{[n]}(A_{\underline{\vec{w}},2}(y))_{[n]} = (B_{\underline{\vec{w}}}(x))^\top_{[n]}(\Omega^\top A_{\vec{w},1}(x))_{[n]},\\
K^{(n)}_{\underline{\vec{w}}}(x,y) &= (B_{\underline{\vec{w}}}(x))^\top_{[n]}(Q_{\underline{\vec{w}}}(y))_{[n]} = (B_{\underline{\vec{w}}}(x))^\top_{[n]}(\Omega^\top Q_{\vec{w}}(x))_{[n]},
\end{aligned}
$$

so that using (41) and back again (42) we obtain the desired result. □

**Lemma 7** *For $n \in \mathbb{N}_0$, we have $A^{(n)}_1(0) \neq 0$.*

***Proof*** If $A^{(n)}_1(0) = 0$, from (40) we deduce that $K^{(n)}_{\vec{w},1}(x,0) = 0$. Attending to (22) we have $K^{(n)}_1(x,0) = \sum\limits_{m=0}^{n-1} B^{(m)}(x)A^{(m)}_1(0)$, as the sequence $\{B^{(n)}(x)\}_{n=0}^\infty$ is linearly independent we deduce that $A^{(m)}_1(0) = 0$, for all $m \in \{0,1,\dots,n\}$, therefore $A^{(0)}_1(0) = 0$, which is impossible. □

**Lemma 8** *For $n \in \mathbb{N}_0$, the matrix coefficient $\Omega_{n+1,n}$ of $\Omega$ is*

$$
\Omega_{n+1,n} = -\frac{A^{(n)}_{1,\vec{w}}(0)}{A^{(n+1)}_{1,\vec{w}}(0)}. \tag{43}
$$

***Proof*** In the first relation in Eq. (38) put $x = 0$ and clean up to get the form of the unknown $\Omega_{n+1,n}$. □

**Theorem 7** (Christoffel formulas) *For $n \in \mathbb{N}_0$, the type I orthogonal polynomials and linear forms fulfill*

$$
\begin{aligned}
Q^{(n)}_{\underline{\vec{w}}}(y) &= Q^{(n)}_{\vec{w}}(y) - \frac{A^{(n)}_{1,\vec{w}}(0)}{A^{(n+1)}_{1,\vec{w}}(0)} Q^{(n+1)}_{\vec{w}}(y),\\
A^{(n)}_{1,\underline{\vec{w}}}(y) &= A^{(n)}_{2,\vec{w}}(y) - \frac{A^{(n)}_{1,\vec{w}}(0)}{A^{(n+1)}_{1,\vec{w}}(0)} A^{(n+1)}_{2,\vec{w}}(y),\\
A^{(n)}_{2,\underline{\vec{w}}}(y) &= \frac{1}{y}\left(A^{(n)}_{1,\vec{w}}(y) - \frac{A^{(n)}_{1,\vec{w}}(0)}{A^{(n+1)}_{1,\vec{w}}(0)} A^{(n+1)}_{1,\vec{w}}(y)\right).
\end{aligned}
$$

*For $n \in \mathbb{N}_0$, the type II orthogonal polynomials satisfy*

$$
\begin{aligned}
B^{(n)}_{\underline{\vec{w}}}(x) = \frac{K^{(n+1)}_1(x,0)}{A^{(n)}_{1,\vec{w}}(0)} &= \frac{1}{x}\left(B^{(n+1)}_{\vec{w}}(x) + \left(\frac{A^{(n-1)}_{1,\vec{w}}(0)}{A^{(n)}_{1,\vec{w}}(0)} + T_{n,n}\right) B^{(n)}_{\vec{w}}(x)\right.\\
&\left. - \frac{A^{(n+1)}_{1,\vec{w}}(0)}{A^{(n)}_{1,\vec{w}}(0)} T_{n+1,n-1} B^{(n-1)}_{\vec{w}}(x)\right).
\end{aligned}
$$

*For $n \in \mathbb{N}_0$, we find*

$$
H_{\underline{\vec{w}},n} = -\frac{A^{(n+1)}_{1,\vec{w}}(0)}{A^{(n)}_{1,\vec{w}}(0)} H_{\vec{w},n}.
$$

***Proof*** We just evaluate (40) and (24) at $x = 0$, use (43) and do some clearing. □

**Definition 11** The Jacobi–Piñeiro multiple orthogonal polynomials, correspond to the choice $w_1 = x^\alpha$, $w_2 := x^\beta$ and $\mathrm{d}\,\mu = (1-x)^\gamma$ with $\alpha, \beta, \gamma > -1$ and $\alpha - \beta \notin \mathbb{Z}$, $\delta = [0, 1]$.

**Theorem 8** (Permuting Christoffel transformation for Jacobi–Piñeiro I) *The transformation $(\alpha, \beta) \to (\beta, \alpha + 1)$ in the Jacobi–Piñeiro case has as connection matrix*

$$
\Omega = \begin{bmatrix}
1 & 0 & \cdots & \cdots \\
\frac{H_1^{\alpha,\beta}}{H_0^{\beta,\alpha+1}} & 1 & \ddots & \\
0 & \frac{H_2^{\alpha,\beta}}{H_1^{\beta,\alpha+1}} & 1 & \ddots \\
\vdots & \ddots & \ddots & \ddots
\end{bmatrix}.
$$

*The corresponding connection formulas*

$$
\begin{gathered}
\Omega\, B^{\beta,\alpha+1} = B^{\alpha,\beta},\\
\Omega^\top A_1^{\alpha,\beta} = x\, A_2^{\beta,\alpha+1},\ \Omega^\top A_2^{\alpha,\beta} = A_1^{\beta,\alpha+1},\\
\Omega^\top Q^{\alpha,\beta} = Q^{\beta,\alpha+1},
\end{gathered}
$$

*hold.*

In [17] the explicit expressions for the Jacobi–Piñeiro's $H_n$ were given as

$$
\begin{aligned}
H_{2n}^{\alpha,\beta} &= \frac{n!\Gamma(2n+1+\gamma)\Gamma(n+1+\alpha)}{\Gamma(3n+2+\alpha+\gamma)} \frac{(\alpha-\beta+1)_n}{(2n+1+\alpha+\gamma)_n(2n+1+\beta+\gamma)_n},\\
H_{2n+1}^{\alpha,\beta} &= \frac{n!\Gamma(2n+2+\gamma)\Gamma(n+1+\beta)}{\Gamma(3n+2+\beta+\gamma)} \frac{(\beta-\alpha)_{n+1}}{(2n+2+\beta+\gamma)_{n+1}(2n+2+\alpha+\gamma)_{n+1}}.
\end{aligned}
$$

Then,

**Corollary 2** (Permuting Christoffel transformation for Jacobi–Piñeiro II) *The connection coefficients are explicitly given by*

$$
\begin{aligned}
\Omega_{2n+1,2n} &= \frac{H_{2n+1}^{\alpha,\beta}}{H_{2n}^{\beta,\alpha+1}} = \frac{(2n+1+\gamma)(2n+1+\beta+\gamma)(\beta-\alpha+n)}{(3n+2+\alpha+\gamma)(3n+1+\beta+\gamma)(3n+2+\beta+\gamma)},\\
\Omega_{2n+2,2n+1} &= \frac{H_{2n+2}^{\alpha,\beta}}{H_{2n+1}^{\beta,\alpha+1}} = \frac{(n+1)(2n+2+\gamma)}{(3n+3+\alpha+\gamma)(3n+4+\alpha+\gamma)}.
\end{aligned}
$$

*In terms of which we have the following connection formulas between Jacobi–Piñeiro multiple orthogonal polynomials with permuted parameters*

$$
\begin{aligned}
\Omega_{2n+1,2n} B^{\beta,\alpha+1}_{(n,n)} + B^{\beta,\alpha+1}_{(n+1,n)} &= B^{\alpha,\beta}_{(n+1,n)},\\
\Omega_{2n+2,2n+1} B^{\beta,\alpha+1}_{(n+1,n)} + B^{\beta,\alpha+1}_{(n+1,n+1)} &= B^{\alpha,\beta}_{(n+1,n+1)},\\
\Omega_{2n+1,2n} A^{\alpha,\beta}_{(n+1,n+1),1} + A^{\alpha,\beta}_{(n+1,n),1} &= x A^{\beta,\alpha+1}_{(n+1,n),2},\\
\Omega_{2n+2,2n+1} A^{\alpha,\beta}_{(n+2,n+1),1} + A^{\beta,\alpha+1}_{(n+1,n+1),1} &= x A^{\beta,\alpha+1}_{(n+1,n+1),2},\\
\Omega_{2n+1,2n} A^{\alpha,\beta}_{(n+1,n+1),2} + A^{\alpha,\beta}_{(n+1,n),2} &= A^{\beta,\alpha+1}_{(n+1,n),1},\\
\Omega_{2n+2,2n+1} A^{\alpha,\beta}_{(n+2,n+1),2} + A^{\alpha,\beta}_{(n+1,n+1),2} &= A^{\beta,\alpha+1}_{(n+1,n+1),1},\\
\Omega_{2n+1,2n} Q^{\alpha,\beta}_{(n+1,n+1)} + Q^{\alpha,\beta}_{(n+1,n)} &= Q^{\beta,\alpha+1}_{(n+1,n)},\\
\Omega_{2n+2,2n+1} Q^{\alpha,\beta}_{(n+2,n+1)} + Q^{\alpha,\beta}_{(n+1,n+1)} &= Q^{\beta,\alpha+1}_{(n+1,n+1)}.
\end{aligned}
$$

We go back to Remark 8 and consider the basic Christoffel transformation

$$
\underline{\underline{\vec{w}}} := x\vec{w}.
$$

The following Lemma was first presented in [22, 7, Lemma 2.4], see also [16, Theorem 9]

**Lemma 9** (Coussement-Van Assche) *For* $n \in \mathbb{N}_0$*, there are nonzero constants* $C_n$ *such that*

$$B^{(n)}(x) = C_n \begin{vmatrix} A_1^{(n)}(x) & A_2^{(n)}(x) \\ A_1^{(n+1)}(x) & A_2^{(n+1)}(x) \end{vmatrix}.$$

***Proof*** Despite this was proven elsewhere [22], for the reader convenience, we give a proof of it. In the one hand, we have that $\begin{vmatrix} A_1^{(n)}(x) & A_2^{(n)}(x) \\ A_1^{(n+1)}(x) & A_2^{(n+1)}(x) \end{vmatrix} = n.$. On the other hand for $k \in \{0, \dots, \lfloor \frac{n-1}{2} \rfloor\}$ we have

$$\begin{aligned}
&\int_\Delta x^k \begin{vmatrix} A_1^{(n)}(x) & A_2^{(n)}(x) \\ A_1^{(n+1)}(x) & A_2^{(n+1)}(x) \end{vmatrix} w_1(x) \,\mathrm{d}\,\mu(x) \\
&= \int_\Delta x^k \begin{vmatrix} w_1(x)A_1^{(n)}(x) & A_2^{(n)}(x) \\ w_1(x)A_1^{(n+1)}(x) & A_2^{(n+1)}(x) \end{vmatrix} \mathrm{d}\,\mu(x) \\
&= \int_\Delta x^k \begin{vmatrix} w_1(x)A_1^{(n)}(x) + w_2(x)A_2^{(n)}(x) & A_2^{(n)}(x) \\ w_1(x)A_1^{(n+1)}(x) + w_2 A_2^{(n+1)}(x) & A_2^{(n+1)}(x) \end{vmatrix} \mathrm{d}\,\mu(x) \\
&= \int_\Delta x^k \begin{vmatrix} Q^{(n)}(x) & A_2^{(n)}(x) \\ Q^{(n+1)}(x) & A_2^{(n+1)}(x) \end{vmatrix} \mathrm{d}\,\mu(x) = 0,
\end{aligned}$$

and for $k \in \{0, \dots, \lfloor \frac{n}{2} \rfloor - 1\}$ we find

$$\begin{aligned}
&\int_\Delta x^k \begin{vmatrix} A_1^{(n)}(x) & A_2^{(n)}(x) \\ A_1^{(n+1)}(x) & A_2^{(n+1)}(x) \end{vmatrix} w_2(x) \,\mathrm{d}\,\mu(x) \\
&= \int_\Delta x^k \begin{vmatrix} A_1^{(n)}(x) & w_2(x)A_2^{(n)}(x) \\ A_1^{(n+1)}(x) & w_2(x)A_2^{(n+1)}(x) \end{vmatrix} \mathrm{d}\,\mu(x) \\
&= \int_\Delta x^k \begin{vmatrix} A_1^{(n)}(x) & w_1(x)A_1^{(n)}(x) + w_2(x)A_2^{(n)}(x) \\ A_1^{(n+1)}(x) & w_1(x)A_1^{(n+1)}(x) + w_2 A_2^{(n+1)}(x) \end{vmatrix} \mathrm{d}\,\mu(x) \\
&= \int_\Delta x^k \begin{vmatrix} A_1^{(n)}(x) & Q^{(n)}(x) \\ A_1^{(n+1)}(x) & Q^{(n+1)}(x) \end{vmatrix} \mathrm{d}\,\mu(x) = 0.
\end{aligned}$$

These two orthogonality relations are precisely those satisfied by the type II multiple orthogonal polynomials (7), and the result follows. □

For the remaining of this section we also assume that the zeros of the type II orthogonal polynomials belong to $\mathring{\Delta}$, and thtat $0 \notin \mathring{\Delta}$ as for an AT system, see [42]. Then, we find

**Lemma 10** *For $n\in\mathbb{N}_0$, the matrix* $\begin{bmatrix}A^{(n+1)}_{\vec{w},1}(0) & A^{(n+1)}_{\vec{w},2}(0)\\ A^{(n+2)}_{\vec{w},1}(0) & A^{(n+2)}_{\vec{w},2}(0)\end{bmatrix}$ *is nonsingular.*

***Proof*** From previous Lemma 9 we have $\begin{vmatrix}A^{(n)}_{\vec{w},1}(0) & A^{(n)}_{\vec{w},2}(0)\\ A^{(n+1)}_{\vec{w},1}(0) & A^{(n+1)}_{\vec{w},2}(0)\end{vmatrix}=C_n^{-1}B^{(n)}(0)\neq 0,$ , and the matrix $\begin{bmatrix}A^{(n+1)}_{\vec{w},1}(0) & A^{(n+1)}_{\vec{w},2}(0)\\ A^{(n+2)}_{\vec{w},1}(0) & A^{(n+2)}_{\vec{w},2}(0)\end{bmatrix}^{-1}$ exists. □

**Theorem 9** (Christoffel formulas) *For $n\in\mathbb{N}_0$, we have the following relations*

$$B^{(n)}_{\underline{\underline{\vec{w}}}}(x)=\frac{1}{x}\frac{\begin{vmatrix}B^{(n)}_{\vec{w}}(0) & B^{(n)}_{\vec{w}}(x)\\ B^{(n+1)}_{\vec{w}}(0) & B^{(n+1)}_{\vec{w}}(x)\end{vmatrix}}{B^{(n)}_{\vec{w}}(0)},\qquad Q^{(n)}_{\underline{\underline{\vec{w}}}}(x)=\frac{1}{x}\frac{\begin{vmatrix}A^{(n)}_{\vec{w},1}(0) & A^{(n)}_{\vec{w},2}(0) & Q^{(n)}_{\vec{w}}(x)\\ A^{(n+1)}_{\vec{w},1}(0) & A^{(n+1)}_{\vec{w},2}(0) & Q^{(n+1)}_{\vec{w}}(x)\\ A^{(n+2)}_{\vec{w},1}(0) & A^{(n+2)}_{\vec{w},2}(0) & Q^{(n+2)}_{\vec{w}}(x)\end{vmatrix}}{\begin{vmatrix}A^{(n+1)}_{\vec{w},1}(0) & A^{(n+1)}_{\vec{w},2}(0)\\ A^{(n+2)}_{\vec{w},1}(0) & A^{(n+2)}_{\vec{w},2}(0)\end{vmatrix}}.$$

*Moreover, if $T_{\vec{w}}=\Omega\,\omega$ is the LU factorization of the original system then $T_{\underline{\underline{\vec{w}}}}=\omega\,\Omega$.*

***Proof*** Notice that now we have (36) , i.e.

$$g_{\underline{\underline{\vec{w}}}}=\Lambda\, g_{\vec{w}}=g_{\vec{w}}\,(\Lambda^\top)^2,$$

in where the bi-Hankel structure is evident. Therefore, assuming the Gauss–Borel factorization of the perturbed moment matrix, we find that

$$S^{-1}_{\underline{\underline{\vec{w}}}}H_{\underline{\underline{\vec{w}}}}\tilde{S}^{-1}_{\underline{\underline{\vec{w}}}}=\Lambda S^{-1}_{\vec{w}}H_{\vec{w}}\tilde{S}^{-\top}_{\vec{w}}=S^{-1}_{\vec{w}}H_{\vec{w}}\tilde{S}^{-\top}_{\vec{w}}(\Lambda^\top)^2.$$

Hence, for each of the relations, we get that

$$H_{\underline{\underline{\vec{w}}}}\big(\tilde{S}_{\vec{w}}\tilde{S}^{-1}_{\underline{\underline{\vec{w}}}}\big)^\top=S_{\underline{\underline{\vec{w}}}}\Lambda S^{-1}_{\vec{w}}H_{\vec{w}},\qquad S_{\vec{w}}S^{-1}_{\underline{\underline{\vec{w}}}}H_{\underline{\underline{\vec{w}}}}=H_{\vec{w}}\Big(\tilde{S}_{\underline{\underline{\vec{w}}}}\Lambda^2\tilde{S}^{-1}_{\vec{w}}\Big)^\top.$$

Then, we see that the matrices

$$\omega:=S_{\underline{\underline{\vec{w}}}}\Lambda S^{-1}_{\vec{w}},\qquad \Omega:=S_{\vec{w}}S^{-1}_{\underline{\underline{\vec{w}}}},\tag{44}$$

are banded matrices. Indeed, from definition $\omega$ is Hessenberg, but as $\omega=H_{\underline{\underline{\vec{w}}}}\big(\tilde{S}_{\vec{w}}\tilde{S}^{-1}_{\underline{\underline{\vec{w}}}}\big)^\top H^{-1}_{\vec{w}}$ we conclude that only its main diagonal and first superdiagonal are possibly nonzero,

$$\omega=\begin{bmatrix}\omega_0 & 1 & 0 & \cdots \\ 0 & \omega_1 & 1 & \ddots \\ \vdots & \ddots & \ddots & \ddots\end{bmatrix},\qquad \omega_n=\frac{H_{\underline{\underline{\vec{w}}},n}}{H_{\vec{w},n}}.$$

We also deduce, as $\Omega$ is by definition lower triangular and as $\Omega = H_{\vec{w}}\left(\tilde{S}_{\underline{\underline{\vec{w}}}}\Lambda^2\tilde{S}_{\vec{w}}^{-1}\right)^\top H_{\underline{\underline{\vec{w}}}}^{-1}$, that only the main diagonal, the first an second subdiagonals are possibly nonzero,

$$\Omega=\begin{bmatrix}1 & 0 & \cdots & & \\ \Omega_{1,0} & 1 & \ddots & & \\ \Omega_{2,0} & \Omega_{2,1} & 1 & \ddots & \\ 0 & \Omega_{3,1} & \Omega_{3,2} & 1 & \ddots \\ \vdots & \ddots & \ddots & \ddots & \ddots\end{bmatrix}, \qquad \Omega_{n+2,n}=\frac{H_{\vec{w},n+2}}{H_{\underline{\underline{\vec{w}}},n}}.$$

From (44) we deduce the following connection formulas

$$xB_{\underline{\underline{\vec{w}}}}=\omega B_{\vec{w}},\quad xA_{\underline{\underline{\vec{w}}},1}=\Omega^\top A_{\vec{w},1},\quad xA_{\underline{\underline{\vec{w}}},2}=\Omega^\top A_{\vec{w},2},\quad xQ_{\underline{\underline{\vec{w}}}}=\Omega^\top Q_{\vec{w}},\tag{45}$$

and evaluating at 0 and using standard techniques [2, 7, 8] and Lemma 10, we get

$$\omega_n=-\frac{B^{(n+1)}_{\vec{w}}(0)}{B^{(n)}_{\vec{w}}(0)},\qquad \begin{bmatrix}\Omega_{n+1,n} & \Omega_{n+2,n}\end{bmatrix}=-\begin{bmatrix}A^{(n)}_{\vec{w},1}(0) & A^{(n)}_{\vec{w},2}(0)\end{bmatrix}\begin{bmatrix}A^{(n+1)}_{\vec{w},1}(0) & A^{(n+1)}_{\vec{w},2}(0)\\ A^{(n+2)}_{\vec{w},1}(0) & A^{(n+2)}_{\vec{w},2}(0)\end{bmatrix}^{-1},$$

that introduced back in(45) leads to

$$xB^{(n)}_{\underline{\underline{\vec{w}}}}(x)=B^{(n+1)}_{\vec{w}}(x)-\frac{B^{(n+1)}_{\vec{w}}(0)}{B^{(n)}_{\vec{w}}(0)}B^{(n)}_{\vec{w}}(x),$$

$$xQ^{(n)}_{\underline{\underline{\vec{w}}}}(x)=Q^{(n)}_{\vec{w}}(x)-\begin{bmatrix}A^{(n)}_{\vec{w},1}(0) & A^{(n)}_{\vec{w},2}(0)\end{bmatrix}\begin{bmatrix}A^{(n+1)}_{\vec{w},1}(0) & A^{(n+1)}_{\vec{w},2}(0)\\ A^{(n+2)}_{\vec{w},1}(0) & A^{(n+2)}_{\vec{w},2}(0)\end{bmatrix}^{-1}\begin{bmatrix}Q^{(n+1)}(x)\\ Q^{(n+2)}(x)\end{bmatrix}.$$

Alternatively, in terms of determinants we get the announced result. □

***Remark 9*** Hence, the transformed recursion matrix is obtained from the $LU$ factorization of the initial recursion matrix by flipping the factors, that is considering the $UL$ factorization.

***Remark 10*** For the type I orthogonal polynomials the Christoffel formulas read as follows

$$A^{(n)}_{\underline{\underline{\vec{w}}},1}(x)=\frac{1}{x}\frac{\begin{vmatrix}A^{(n)}_{\vec{w},1}(0) & A^{(n)}_{\vec{w},2}(0) & A^{(n)}_{\vec{w},1}(x)\\ A^{(n+1)}_{\vec{w},1}(0) & A^{(n+1)}_{\vec{w},2}(0) & A^{(n+1)}_{\vec{w},1}(x)\\ A^{(n+2)}_{\vec{w},1}(0) & A^{(n+2)}_{\vec{w},2}(0) & A^{(n+2)}_{\vec{w},1}(x)\end{vmatrix}}{\begin{vmatrix}A^{(n+1)}_{\vec{w},1}(0) & A^{(n+1)}_{\vec{w},2}(0)\\ A^{(n+2)}_{\vec{w},1}(0) & A^{(n+2)}_{\vec{w},2}(0)\end{vmatrix}},$$

$$A^{(n)}_{\underline{\underline{\vec{w}}},2}(x) = \frac{1}{x}\frac{\begin{vmatrix} A^{(n)}_{\vec{w},2}(0) & A^{(n)}_{\vec{w},2}(0) & A^{(n)}_{\vec{w},2}(x) \\ A^{(n+1)}_{\vec{w},1}(0) & A^{(n+1)}_{\vec{w},2}(0) & A^{(n+1)}_{\vec{w},2}(x) \\ A^{(n+2)}_{\vec{w},1}(0) & A^{(n+2)}_{\vec{w},2}(0) & A^{(n+2)}_{\vec{w},2}(x) \end{vmatrix}}{\begin{vmatrix} A^{(n+1)}_{\vec{w},1}(0) & A^{(n+1)}_{\vec{w},2}(0) \\ A^{(n+2)}_{\vec{w},1}(0) & A^{(n+2)}_{\vec{w},2}(0) \end{vmatrix}}.$$

### 3.1 General Christoffel transformations

Motivated by the previous examples we now construct a more general Christoffel transformation. Here we follow the ideas in [2, 26]. Let us consider a perturbing monic polynomial $P \in \mathbb{R}[x]$ with $\deg P = N$, and let consider its decomposition into even and odd parts

$$P(x) = P_{\mathrm{e}}(x^2) + xP_{\mathrm{o}}(x^2), \qquad P_{\mathrm{e}}, P_{\mathrm{o}} \in \mathbb{R}[x], \tag{46}$$

so that

$$2P_{\mathrm{e}}(x^2) = P(x) + P(-x), \qquad 2xP_{\mathrm{o}}(x^2) = P(x) - P(-x), \tag{47}$$

and let us define the polynomial

$$\pi(x) := P_{\mathrm{e}}^2(x) - xP_{\mathrm{o}}^2(x).$$

**Definition 12** We say that a polynomial is non-symmetric if whenever $x_0 \neq 0$ is a root then $-x_0$ is not.

**Lemma 11** *(i) Any monic polynomial $\tilde{P}$, can be factorized as $\tilde{P}(x) = p(x^2)P(x)$, where $P = (x - x_1)\cdots(x - x_N)$ is a non-symmetric monic polynomial and $p(x) = (x - r_1^2)\cdots(x - r_M^2)$. Thus, $\{\pm r_1, \dots, \pm r_M, x_1, \dots x_N\}$ is the zero set for $P$, and if $\tilde{N} = \deg \tilde{P}$, we have $\tilde{N} = 2M + N$.*

*(ii) We have that $P_{\mathrm{e}}, P_{\mathrm{o}} \neq 0$ and*

$$\tilde{P}_{\mathrm{e}}(x) = p(x)P_{\mathrm{e}}(x), \qquad \tilde{P}_{\mathrm{o}}(x) = p(x)P_{\mathrm{o}}(x).$$

*(iii) If $x_0$ is a non-zero root of $P$, then $P_{\mathrm{e}}(x_0) \neq 0$ and $P_{\mathrm{o}}(x_0) \neq 0$.*

*(iv) $x_0$ is a root of $P(x)$ if and only if $x_0^2$ is a root of $\pi$.*

***Proof*** (i) It follows from the Fundamental Theorem of Algebra.

(ii) If $x_0 \neq 0$ is a root of $\tilde{P}$ such that $-x_0$ is not then Eq. (47) gives

$$2P_{\mathrm{e}}(x_0^2) = P(-x_0) \neq 0, \qquad 2x_0P_{\mathrm{o}}(x_0^2) = -P(-x_0) \neq 0$$

and we conclude that $P_{\mathrm{e}}, P_{\mathrm{o}} \neq 0$. The second statement follows from the even/odd decomposition (46).

(iii) It follows from (47) and the non-symmetric character of $P$.
(iv) From (47) we get

$$\begin{aligned}4P_{\mathrm{e}}^2(x^2) &= P^2(x)+P^2(-x)+2P(x)P(-x),\\ 4x^2P_{\mathrm{o}}^2(x^2) &= P^2(x)+P^2(-x)-2P(x)P(-x),\end{aligned}$$

so that $\pi(x^2)=P(x)P(-x)$.

□

Then, we define Christoffel perturbation of the vector of weights as follows

$$\hat{\vec{w}}=\vec{w}\begin{bmatrix}\tilde{P}_{\mathrm{e}}(x) & x\tilde{P}_{\mathrm{o}}(x)\\ \tilde{P}_{\mathrm{o}}(x) & \tilde{P}_{\mathrm{e}}(x)\end{bmatrix}=\vec{w}p(x)\begin{bmatrix}P_{\mathrm{e}}(x) & xP_{\mathrm{o}}(x)\\ P_{\mathrm{o}}(x) & P_{\mathrm{e}}(x)\end{bmatrix}. \tag{48}$$

Then,

**Proposition 12** *The moment matrices satisfy*

$$g_{\hat{\vec{w}}}=g_{\vec{w}}\tilde{P}(\Lambda^\top).$$

***Proof*** It follows from $\Lambda X_1=xX_2$ and $\Lambda X_2=X_1$ and the expressions of the moment matrices. □

**Proposition 13** *Let us assume that the moment matrices $g_{\hat{\vec{w}}}$, $g_{\vec{w}}$ have a Gauss–Borel factorization, i.e.*

$$g_{\hat{\vec{w}}}=\hat{S}^{-1}\hat{H}\hat{\tilde{S}}^{-\top}, \qquad g_{\vec{w}}=S^{-1}H\tilde{S}^{-\top}.$$

*Then, for the connection matrix*

$$\Omega:=S\hat{S}^{-1}$$

*one has the alternative expression*

$$\Omega=H\tilde{S}^{-\top}\tilde{P}(\Lambda^\top)\hat{\tilde{S}}^\top\hat{H}^{-1},$$

*so that is lower unitriangular matrix with only its first $\tilde{N}$ subdiagonals possibly nonzero.*

***Proof*** Direct consequence of the Gauss–Borel factorization. □

**Proposition 14** (Polynomial connection formulas) *The following formulas hold*

$$\begin{aligned}\Omega\hat{B}(x)&=B(x),\\ \Omega^\top A_1(x)&=p(x)\big(P_{\mathrm{e}}(x)\hat{A}_1(x)+xP_{\mathrm{o}}(x)\hat{A}_2(x)\big),\\ \Omega^\top A_2(x)&=p(x)\big(P_{\mathrm{o}}(x)\hat{A}_1(x)+P_{\mathrm{e}}(x)\hat{A}_2(x)\big).\end{aligned}$$

***Proof*** Use the definitions $B = SX$, $A_1 = H^{-1}\tilde{S}X_1$, $A_2 = H^{-1}\tilde{S}X_2$ , $\hat{B} = \hat{S}X$, $\hat{A}_1 = \hat{H}^{-1}\tilde{\hat{S}}X_1$, $\hat{A}_2 = \hat{H}^{-1}\tilde{\hat{S}}X_2$ and the two expressions for $\Omega$. The following

$$\begin{aligned}\tilde{P}(\Lambda)X_1(x) &= p(\Lambda^2)\big(P_{\rm e}(\Lambda^2) + \Lambda P_{\rm o}(\Lambda^2)\big)X_1 = p(x)\big(P_{\rm e}(x)X_1(x) + xP_{\rm o}(x)X_2(x)\big),\\ \tilde{P}(\Lambda)X_2(x) &= p(\Lambda^2)\big(P_{\rm e}(\Lambda^2) + \Lambda P_{\rm o}(\Lambda^2)\big)X_2 = p(x)\big(P_{\rm o}(x)X_1(x) + P_{\rm e}(x)X_2(x)\big),\end{aligned}$$

leads to the desired representations. □

**Corollary 3** *We have the relations*

$$\begin{aligned}\Omega^\top\big(P_{\rm e}(x)A_1(x) - xP_{\rm o}(x)A_2(x)\big) &= \pi(x)p(x)\hat{A}_1(x),\\ \Omega^\top\big(-P_{\rm o}(x)A_1(x) + P_{\rm e}(x)A_2(x)\big) &= \pi(x)p(x)\hat{A}_2(x).\end{aligned}$$

***Proof*** Solve for $p(x)\hat{A}_1(x)$ and $p(x)\hat{A}_2(x)$ the last two equations in Proposition 14. □

***Remark 11*** Notice that $\tilde{\pi} = p^2\pi$ is not $p\pi$.

**Definition 13** Let us consider the polynomials

$$\begin{aligned}v_1^{(n)}(x) &:= P_{\rm e}(x)A_1^{(n)}(x) - xP_{\rm o}(x)A_2^{(n)}(x),\\ v_2^{(n)}(x) &:= -P_{\rm o}(x)A_1^{(n)}(x) + P_{\rm e}(x)A_2^{(n)}(x).\end{aligned}$$

**Definition 14** For $n > \tilde{N}$, let us introduce the matrix

$$\Omega[n] = \begin{bmatrix}\Omega_{n,n-\tilde{N}} & 0 & \cdots & 0\\ \vdots & \ddots & \ddots & \vdots\\ \vdots & & \ddots & 0\\ \Omega_{n,n-1} & \cdots & \cdots & \Omega_{n+\tilde{N}-1,n-1}\end{bmatrix}.$$

**Lemma 12** *Given any two semi-infinite vectors $C$, $D$, we find*

$$D_{[n]}^\top(\Omega^\top C)_{[n]} = (\Omega D)_{[n]}^\top C_{[n]} + \begin{bmatrix}D^{(n-\tilde{N})} & \cdots & D^{(n-1)}\end{bmatrix}\Omega[n]\begin{bmatrix}C^{(n)}\\ \vdots\\ C^{(n+\tilde{N}-1)}\end{bmatrix}, \tag{49}$$

*with $C^{(n)}$, $D^{(n)}$, $n \in \mathbb{N}_0$, the corresponding entries of the vectors $C$ and $D$.*

***Proof*** Given the band form of $\Omega$ we find for any two semi-infinite vectors $C$, $D$ that

$$D_{[n]}^{\top}(\Omega^{\top}C)_{[n]} = (\Omega D)_{[n]}^{\top}C_{[n]} + D_{[n]}^{\top}\begin{bmatrix} 0 \\ \vdots \, (\tilde{N}-n)\text{times} \\ 0 \\ \Omega_{n,n-\tilde{N}}C^{(n)} \\ \Omega_{n,n-\tilde{N}+1}C^{(n)} + \Omega_{n+1,n-\tilde{N}+1}C^{(n+1)} \\ \vdots \\ \Omega_{n,n-1}C^{(n)} + \cdots + \Omega_{n+\tilde{N}-1,n-1}C^{(n+\tilde{N}-1)} \end{bmatrix},$$

and the result is proven. □

**Definition 15** Let us use the notation

$$\begin{aligned} k_1^{(n)}(x,y) &:= P_{\mathrm{e}}(y)K_1^{(n)}(x,y) - xP_{\mathrm{o}}(y)K_2^{(n)}(x,y), \\ k_2^{(n)}(x,y) &:= -P_{\mathrm{o}}(y)K_1^{(n)}(x,y) + P_{\mathrm{e}}(y)K_2^{(n)}(x,y). \end{aligned}$$

**Proposition 15** (CD kernels connection formulas) *The CD kernels satisfy*

$$\pi(y)p(y)\hat{K}_1^{(n)}(x,y) = k_1^{(n)}(x,y) + \left[\hat{B}^{(n-\tilde{N})}(x) \cdots\cdots \hat{B}^{(n-1)}(x)\right]\Omega[n]\begin{bmatrix} v_1^{(n)}(y) \\ \vdots \\ v_1^{(\tilde{N}+n-1)}(y) \end{bmatrix},$$

$$\pi(y)p(y)\hat{K}_2^{(n)}(x,y) = k_2^{(n)}(x,y) + \left[\hat{B}^{(n-\tilde{N})}(x) \cdots\cdots \hat{B}^{(n-1)}(x)\right]\Omega[n]\begin{bmatrix} v_2^{(n)}(y) \\ \vdots \\ v_2^{(\tilde{N}+n-1)}(y) \end{bmatrix}.$$

***Proof*** Follows from Corollary 3, Lemma 12 and the definition of the CD kernels. □

The set of zeros, $Z_{p\pi}$, of $p(z)\pi(z)$ is

$$Z_{p\pi} := \{r_1^2, \ldots, r_M^2, x_1^2, \ldots, x_N^2\}. \tag{50}$$

From heron we assume that all these zeros are simple. That is, we assume that

$$\begin{aligned} r_i^2 &\neq r_j^2, && \text{for } 1 \leq i < j \leq M, \\ r_i^2 &\neq x_j^2, && \text{for } 1 \leq i \leq M \text{ and } 1 \leq j \leq N, \\ x_i^2 &\neq x_j^2, && \text{for } 1 \leq i < j \leq N. \end{aligned}$$

In other words, the zeros of $p(x)$ and $P(x)$ are simple and the zeros of $p(x)$ are different from those of $P(x)$.

**Definition 16** We introduce the function

$$\tau_n := \begin{vmatrix} v_1^{(n)}(r_1^2) \cdots\cdots v_1^{(n)}(r_M^2) & v_2^{(n)}(r_1^2) \cdots\cdots v_2^{(n)}(r_M^2) & v_2^{(n)}(x_1^2) \cdots\cdots v_2^{(n)}(x_N^2) \\ \vdots \qquad \vdots & \vdots \qquad \vdots & \vdots \qquad \vdots \\ v_1^{(\tilde{N}+n-1)}(r_1^2) \cdots v_1^{(\tilde{N}+n-1)}(r_M^2) & v_2^{(\tilde{N}+n-1)}(r_1^2) \cdots v_2^{(\tilde{N}+n-1)}(r_M^2) & v_2^{(\tilde{N}+n-1)}(x_1^2) \cdots v_2^{(\tilde{N}+n-1)}(x_N^2) \end{vmatrix}.$$

**Theorem 10** *For* $n \geq \tilde{N}$, $\tau_n \neq 0$.

***Proof*** As we have

$$\begin{aligned} &\left[\hat{B}^{(n-\tilde{N})}(x) \cdots \hat{B}^{(n-1)}(x)\right] \Omega[n] \\ &\times \begin{bmatrix} v_1^{(n)}(r_1^2) \cdots\cdots v_1^{(n)}(r_M^2) & v_2^{(n)}(r_1^2) \cdots\cdots v_2^{(n)}(r_M^2) & v_2^{(n)}(x_1^2) \cdots\cdots v_2^{(n)}(x_N^2) \\ \vdots \qquad \vdots & \vdots \qquad \vdots & \vdots \qquad \vdots \\ v_1^{(\tilde{N}+n-1)}(r_1^2) \cdots v_1^{(\tilde{N}+n-1)}(r_M^2) & v_2^{(\tilde{N}+n-1)}(r_1^2) \cdots v_2^{(\tilde{N}+n-1)}(r_M^2) & v_2^{(\tilde{N}+n-1)}(x_1^2) \cdots v_2^{(\tilde{N}+n-1)}(x_N^2) \end{bmatrix} \\ &= -\left[k_1^{(n)}(x, r_1^2) \cdots k_1^{(n)}(x, r_M^2) \; k_2^{(n)}(x, r_1^2) \cdots k_2^{(n)}(x, r_M^2) \; k_2^{(n)}(x, x_1^2) \cdots k_2^{(n)}(x, x_N^2)\right], \end{aligned} \tag{51}$$

if we assume that $\tau_n = 0$ we conclude that there exists a non-zero vector $\left[c_1 \cdots c_{\tilde{N}}\right]^\top$ such that

$$\begin{aligned} &c_1 k_1^{(n)}(x, r_1^2) + \cdots + c_M k_1^{(n)}(x, r_M^2) + c_{M+1} k_2^{(n)}(x, r_1^2) + \cdots + c_{2M} k_2^{(n)}(x, r_M^2) \\ &\quad + c_{2M+1} k_2^{(n)}(x, x_1^2) + \cdots + c_{\tilde{N}} k_2^{(n)}(x, x_N^2) = 0. \end{aligned}$$

Hence, replacing the explicit expressions of the kernel polynomials, we conclude

$$\begin{aligned} &\sum_{l=0}^{n-1} \Big(c_1 v_1^{(l)}(r_1^2) + \cdots + c_M v_1^{(l)}(r_M^2) + c_{M+1} v_2^{(l)}(r_1^2) + \cdots + c_{2M} v_2^{(l)}(r_M^2) \\ &\quad + c_{2M+1} v_2^{(l)}(x_1^2) + \cdots + c_{\tilde{N}} v_2^{(l)}(x_N^2)\Big) B^{(l)}(x) = 0. \end{aligned}$$

As the type II orthogonal polynomials are linearly independent we get that

$$\begin{aligned} &c_1 v_1^{(l)}(r_1^2) + \cdots + c_M v_1^{(l)}(r_M^2) + c_{M+1} v_2^{(l)}(r_1^2) + \cdots + c_{2M} v_2^{(l)}(r_M^2) \\ &\quad + c_{2M+1} v_2^{(l)}(x_1^2) + \cdots + c_{\tilde{N}} v_2^{(l)}(x_N^2) = 0. \end{aligned}$$

Hence, discussing the linear system that appears by considering the equations for $l \in \{0, 1, \ldots, \tilde{N} - 1\}$ we conclude that if $\tau_0 \neq 0$ then $\tau_n \neq 0$ for $n \geq \tilde{N}$. Now, we observe that if $\det(v_1, \ldots, v_{\tilde{N}})$ denotes the determinant considered as a multi-linear function of the columns $v_j$ of the matrix, recalling that for type I polynomials we have $(A_a(x))_{[n]} = (H^{-1}\tilde{S})_{[n]}(X_a(x))_{[n]}$ and that $\det(H^{-1}\tilde{S})_{[n]} = \frac{1}{H_0H_1\cdots H_{n-1}}$ we get

$$\tau_0 = \frac{\det \Theta}{H_0H_1\cdots H_{2M+N-1}}$$

with

$$\Theta := \begin{bmatrix} P_{\mathrm{e}}(r_1^2)\cdots\cdots P_{\mathrm{e}}(r_M^2) & -P_{\mathrm{o}}(r_1^2)\cdots\cdots -P_{\mathrm{o}}(r_M^2) & -P_{\mathrm{o}}(x_1^2)\cdots\cdots -P_{\mathrm{o}}(x_N^2) \\ -r_1^2P_{\mathrm{o}}(r_1^2)\cdots -r_M^2P_{\mathrm{o}}(r_M^2) & P_{\mathrm{e}}(r_1^2)\cdots\cdots P_{\mathrm{e}}(r_M^2) & P_{\mathrm{e}}(x_1^2)\cdots\cdots P_{\mathrm{e}}(x_N^2) \\ r_1^2P_{\mathrm{e}}(r_1^2)\cdots r_M^2P_{\mathrm{e}}(r_M^2) & -r_1^2P_{\mathrm{o}}(r_1^2)\cdots -r_M^2P_{\mathrm{o}}(r_M^2) & -x_1^2P_{\mathrm{o}}(x_1^2)\cdots -x_N^2P_{\mathrm{o}}(x_N^2) \\ -r_1^4P_{\mathrm{o}}(r_1^2)\cdots -r_M^4P_{\mathrm{o}}(r_M^2) & r_1^2P_{\mathrm{e}}(r_1^2)\cdots r_M^2P_{\mathrm{e}}(r_M^2) & x_1^2P_{\mathrm{e}}(x_1^2)\cdots x_N^2P_{\mathrm{e}}(x_N^2) \\ \vdots \qquad \vdots & \vdots \qquad \vdots & \vdots \qquad \vdots \\ r_1^{\tilde{N}-1}P_{\mathrm{e}}(r_1^2)\cdots r_M^{\tilde{N}-1}P_{\mathrm{e}}(r_M^2) & -r_1^{\tilde{N}-1}P_{\mathrm{o}}(r_1^2)\cdots -r_M^{\tilde{N}-1}P_{\mathrm{o}}(r_M^2) & -x_1^{\tilde{N}-1}P_{\mathrm{o}}(x_1^2)\cdots -x_N^{\tilde{N}-1}P_{\mathrm{o}}(x_N^2) \end{bmatrix}$$

for $N$ odd and

$$\Theta := \begin{bmatrix} P_{\mathrm{e}}(r_1^2)\cdots\cdots P_{\mathrm{e}}(r_M^2) & -P_{\mathrm{o}}(r_1^2)\cdots\cdots -P_{\mathrm{o}}(r_M^2) & -P_{\mathrm{o}}(x_1^2)\cdots\cdots -P_{\mathrm{o}}(x_N^2) \\ -r_1^2P_{\mathrm{o}}(r_1^2)\cdots -r_M^2P_{\mathrm{o}}(r_M^2) & P_{\mathrm{e}}(r_1^2)\cdots\cdots P_{\mathrm{e}}(r_M^2) & P_{\mathrm{e}}(x_1^2)\cdots\cdots P_{\mathrm{e}}(x_N^2) \\ r_1^2P_{\mathrm{e}}(r_1^2)\cdots r_M^2P_{\mathrm{e}}(r_M^2) & -r_1^2P_{\mathrm{o}}(r_1^2)\cdots -r_M^2P_{\mathrm{o}}(r_M^2) & -x_1^2P_{\mathrm{o}}(x_1^2)\cdots -x_N^2P_{\mathrm{o}}(x_N^2) \\ -r_1^4P_{\mathrm{o}}(r_1^2)\cdots -r_M^4P_{\mathrm{o}}(r_M^2) & r_1^2P_{\mathrm{e}}(r_1^2)\cdots r_M^2P_{\mathrm{e}}(r_M^2) & x_1^2P_{\mathrm{e}}(x_1^2)\cdots x_N^2P_{\mathrm{e}}(x_N^2) \\ \vdots \qquad \vdots & \vdots \qquad \vdots & \vdots \qquad \vdots \\ r_1^{\tilde{N}-2}P_{\mathrm{e}}(r_1^2)\cdots r_M^{\tilde{N}-2}P_{\mathrm{e}}(r_M^2) & -r_1^{\tilde{N}-2}P_{\mathrm{o}}(r_1^2)\cdots -r_M^{\tilde{N}-2}P_{\mathrm{o}}(r_M^2) & -x_1^{\tilde{N}-2}P_{\mathrm{o}}(x_1^2)\cdots -x_N^{\tilde{N}-2}P_{\mathrm{o}}(x_N^2) \\ -r_1^{\tilde{N}}P_{\mathrm{o}}(r_1^2)\cdots -r_M^{\tilde{N}}P_{\mathrm{o}}(r_M^2) & r_1^{\tilde{N}-2}P_{\mathrm{e}}(r_1^2)\cdots r_M^{\tilde{N}-2}P_{\mathrm{e}}(r_M^2) & x_1^{\tilde{N}-2}P_{\mathrm{e}}(x_1^2)\cdots x_N^{\tilde{N}-2}P_{\mathrm{e}}(x_N^2) \end{bmatrix}.$$

for $N$ even.

To analyze $\det\Theta$ let us split the matrix $\Theta$ into two submatrices $\Theta_1$, with the first $2M$ columns, and the submatrix $\Theta_2$ with the final $N$ columns. Regarding $\Theta_2$, recall that $P(x_j) = 0$ so that $P_{\mathrm{e}}(x_j^2) = -x_jP_{\mathrm{o}}(x_j^2)$, therefore we can write

$$\Theta_2 = \begin{bmatrix} -P_{\mathrm{o}}(x_1^2)\cdots\cdots -P_{\mathrm{o}}(x_N^2) \\ -x_1P_{\mathrm{o}}(x_1^2)\cdots -x_NP_{\mathrm{o}}(x_N^2) \\ -x_1^2P_{\mathrm{o}}(x_1^2)\cdots -x_N^2P_{\mathrm{o}}(x_N^2) \\ -x_1^3P_{\mathrm{o}}(x_1^2)\cdots -x_N^3P_{\mathrm{o}}(x_N^2) \\ \vdots \qquad \vdots \\ -x_1^{\tilde{N}-1}P_{\mathrm{o}}(x_1^2)\cdots -x_N^{\tilde{N}-1}P_{\mathrm{o}}(x_N^2) \end{bmatrix}.$$

To study $\Theta_1$ use (47) as well as

$$P_{\mathrm{e}}(r_i^2) = \frac{P(r_i)+P(-r_i)}{2}, \qquad P_{\mathrm{o}}(r_i^2) = \frac{P(r_i)-P(-r_i)}{2r_i},$$

so that for $N$ odd

$$\Theta_1 := \begin{bmatrix} \frac{P(r_1)+P(-r_1)}{2} \cdots\cdots \frac{P(r_M)+P(-r_M)}{2} & -\frac{P(r_1)-P(-r_1)}{2r_1} \cdots\cdots -\frac{P(r_M)-P(-r_M)}{2r_M} \\ -r_1\frac{P(r_1)-P(-r_1)}{2} \cdots -r_M\frac{P(r_M)-P(-r_M)}{2} & \frac{P(r_1)+P(-r_1)}{2} \cdots\cdots \frac{P(r_M)+P(-r_M)}{2} \\ r_1^2\frac{P(r_1)+P(-r_1)}{2} \cdots r_M^2\frac{P(r_M)+P(-r_M)}{2} & -r_1\frac{P(r_1)-P(-r_1)}{2} \cdots -r_M\frac{P(r_M)-P(-r_M)}{2} \\ -r_1^3\frac{P(r_1)-P(-r_1)}{2} \cdots -r_M^3\frac{P(r_1)-P(-r_1)}{2} & r_1^2\frac{P(r_1)+P(-r_1)}{2} \cdots r_M^2\frac{P(r_M)+P(-r_M)}{2} \\ \vdots \qquad\qquad \vdots & \vdots \qquad\qquad \vdots \\ r_1^{\tilde N-1}\frac{P(r_1)+P(-r_1)}{2} \cdots r_M^{\tilde N-1}\frac{P(r_M)+P(-r_M)}{2} & -r_1^{\tilde N-2}\frac{P(r_1)-P(-r_1)}{2} \cdots -r_M^{\tilde N-2}\frac{P(r_M)-P(-r_M)}{2} \end{bmatrix}.$$

We make a first transformation by performing column operations, that is we add to the $i$-th column, $i \in \{1, \dots, M\}$, $r_i$ times the $(M+i)$-th column::

$$\Theta_1' := \begin{bmatrix} P(-r_1) \cdots\cdots P(-r_M) & -\frac{P(r_1)-P(-r_1)}{2r_1} \cdots\cdots -\frac{P(r_M)-P(-r_M)}{2r_M} \\ r_1P(-r_1) \cdots r_MP(-r_M) & \frac{P(r_1)+P(-r_1)}{2} \cdots\cdots \frac{P(r_M)+P(-r_M)}{2} \\ r_1^2P(-r_1) \cdots r_M^2P(-r_M) & -r_1\frac{P(r_1)-P(-r_1)}{2} \cdots -r_M\frac{P(r_M)-P(-r_M)}{2} \\ r_1^3P(-r_1) \cdots r_M^3P(-r_M) & r_1^2\frac{P(r_1)+P(-r_1)}{2} \cdots r_M^2\frac{P(r_M)+P(-r_M)}{2} \\ \vdots \qquad\qquad \vdots & \vdots \qquad\qquad \vdots \\ r_1^{\tilde N-1}P(-r_1) \cdots r_M^{\tilde N-1}P(-r_M) & -r_1^{\tilde N-2}\frac{P(r_1)-P(-r_1)}{2} \cdots -r_M^{\tilde N-2}\frac{P(r_M)-P(-r_M)}{2} \end{bmatrix},$$

then we proceed with a second column operation by subtracting in $\Theta_1'$ to the $(M+i)$-th column $\frac{1}{2r_i}$ times the $i$-th column:

$$\Theta_1'' := \begin{bmatrix} P(-r_1) \cdots\cdots P(-r_M) & -\frac{P(r_1)}{2r_1} \cdots\cdots -\frac{P(r_M)}{2r_M} \\ r_1P(-r_1) \cdots r_MP(-r_M) & \frac{P(r_1)}{2} \cdots\cdots \frac{P(r_M)}{2} \\ r_1^2P(-r_1) \cdots r_M^2P(-r_M) & -r_1\frac{P(r_1)}{2} \cdots -r_M\frac{P(r_M)}{2} \\ r_1^3P(-r_1) \cdots r_M^3P(-r_M) & r_1^2\frac{P(r_1)}{2} \cdots r_M^2\frac{P(r_M)}{2} \\ \vdots \qquad\qquad \vdots & \vdots \qquad\qquad \vdots \\ r_1^{\tilde N-1}P(-r_1) \cdots r_M^{\tilde N-1}P(-r_M) & -r_1^{\tilde N-2}\frac{P(r_1)}{2} \cdots -r_M^{\tilde N-2}\frac{P(r_M)}{2} \end{bmatrix}.$$

Hence, recalling that $\pi(x^2) = P(x)P(-x)$ we get

$$\det \Theta = (-1)^N \prod_{i=1}^{M} \frac{\pi(r_i^2)}{2r_i} \prod_{j=1}^{N} P_{\mathrm{o}}(x_j^2) \det \Theta''',$$

with

$$\Theta''' = \begin{bmatrix} 1 \cdots\cdots\cdots 1 & -1 \cdots\cdots\cdots -1 & 1 \cdots\cdots\cdots 1 \\ r_1 \cdots\cdots\cdots r_M & r_1 \cdots\cdots\cdots r_M & x_1 \cdots\cdots\cdots x_N \\ r_1^2 \cdots\cdots\cdots r_M^2 & -r_1^2 \cdots\cdots\cdots -r_M^2 & x_1^2 \cdots\cdots\cdots x_N^2 \\ r_1^3 \cdots\cdots\cdots r_M^3 & r_1^3 \cdots\cdots\cdots r_M^3 & x_1^3 \cdots\cdots\cdots x_N^3 \\ \vdots \qquad \vdots & \vdots \qquad \vdots & \vdots \qquad \vdots \\ r_1^{\tilde{N}-1} \cdots\cdots r_M^{\tilde{N}-1} & -r_1^{\tilde{N}-1} \cdots\cdots -r_M^{\tilde{N}-1} & x_1^{\tilde{N}-1} \cdots\cdots x_N^{\tilde{N}-1} \end{bmatrix}. \tag{52}$$

Finally, by adding to the $(M+i)$-th column the $i$-th column, normalizing, and then in the resulting matrix subtracting to the $i$-th column the $(M+i)$-column we get

$$\det \Theta = (-1)^N \prod_{i=1}^{M} \frac{\pi(r_i^2)}{2r_i} \prod_{j=1}^{N} P_{\mathrm{o}}(x_j^2) \det \Theta'''.$$

Now, we observe that

$$\pi(r_i^2) = \prod_{j=1}^{N} (r_i^2 - x_j^2), \quad P_{\mathrm{o}}(x_j^2) = \frac{-P(-x_j)}{2x_j} = \\ - \prod_{\substack{1 \leq i \leq N \\ i \neq j}} (-x_j - x_i) = (-1)^{N-1} \prod_{\substack{1 \leq i \leq N \\ i \neq j}} (x_j + x_i).$$

to find that

$$\det \Theta = (-1)^N \prod_{\substack{1 \leq i \leq M \\ 1 \leq j \leq N}} \frac{r_i^2 - x_j^2}{2r_i} \prod_{\substack{1 \leq i, j \leq N \\ i \neq j}} (x_j + x_i)^2 \det \Theta'''.$$

In order to compute $\det \Theta'''$ we replace in (52) in the $2M$-th column $r_M$ by $x$. The determinant of the resulting matrix is an $(\tilde{N}-1)$-th degree polynomial in the $x$ variable with zeros at $\{\pm r_1, \ldots, \pm r_{M-1}, r_M, -x_1, \ldots, -x_N\}$. The leading coefficient of this polynomial is given by $(-1)^{N+1}\det\Theta^{IV}$ where

$$\Theta^{IV}=\begin{bmatrix} 1\cdots\cdots 1 & -1\cdots\cdots -1 & 1\cdots\cdots 1\\ r_1\cdots\cdots r_M & r_1\cdots\cdots r_{M-1} & x_1\cdots\cdots x_N\\ r_1^2\cdots\cdots r_M^2 & -r_1^2\cdots -r_{M-1}^2 & x_1^2\cdots\cdots x_N^2\\ r_1^3\cdots\cdots r_M^3 & r_1^3\cdots\cdots r_{M-1}^3 & x_1^3\cdots\cdots x_N^3\\ \vdots\qquad\vdots & \vdots\qquad\vdots & \vdots\qquad\vdots\\ r_1^{\tilde N-2}\cdots r_M^{\tilde N-2} & r_1^{\tilde N-2}\cdots r_{M-1}^{\tilde N-2} & x_1^{\tilde N-2}\cdots x_N^{\tilde N-2}\end{bmatrix}.$$

To compute $\det\Theta^{IV}$, we repeat the previous idea and replace in the $M$-th column $r_M$ by $x$. This is a polynomial of degree $(\tilde N-2)$ with zeros located at $\{\pm r_1,\dots,\pm r_{M-1},x_1,\dots,x_N\}$ and leading coefficient $(-1)^{M+N}\det\Theta^V$, with

$$\Theta^{V}=\begin{bmatrix} 1\cdots\cdots 1 & -1\cdots\cdots -1 & 1\cdots\cdots 1\\ r_1\cdots\cdots r_{M-1} & r_1\cdots\cdots r_{M-1} & x_1\cdots\cdots x_N\\ r_1^2\cdots\cdots r_{M-1}^2 & -r_1^2\cdots -r_{M-1}^2 & x_1^2\cdots\cdots x_N^2\\ r_1^3\cdots\cdots r_M^3 & r_1^3\cdots\cdots r_{M-1}^3 & x_1^3\cdots\cdots x_N^3\\ \vdots\qquad\vdots & \vdots\qquad\vdots & \vdots\qquad\vdots\\ r_1^{\tilde N-3}\cdots r_{M-1}^{\tilde N-3} & -r_1^{\tilde N-3}\cdots -r_{M-1}^{\tilde N-3} & x_1^{\tilde N-3}\cdots x_N^{\tilde N-3}\end{bmatrix}.$$

Hence,

$$\det\Theta'''=(-1)^M\prod_{i=1}^{M-1}\big(r_M^2-r_i^2\big)^2\prod_{j=1}^{N}(r_M^2-x_i^2)2r_M\det\Theta^V,$$

and by induction we get

$$\det\Theta'''=(-1)^{\frac{M(M+1)}{2}}\prod_{i=1}^{M}(2r_i)\prod_{1\le i<j\le M}(r_i^2-r_j^2)^2\prod_{\substack{1\le i\le M\\1\le j\le N}}(r_i^2-x_j^2)\prod_{1\le i<j\le N}(x_i-x_j).$$

Then, collecting all these facts, we get

$$\begin{aligned}\det\Theta=(-1)^{\frac{M(M+1)}{2}+N}&\prod_{1\le i<j\le M}(r_i^2-r_j^2)^2\prod_{\substack{1\le i\le M\\1\le j\le N}}(r_i^2-x_j^2)^2\prod_{1\le i<j\le N}(x_i-x_j)\\ &\times\prod_{1\le i<j\le N}(x_i+x_j)^2.\end{aligned}\tag{53}$$

Hence $\det\Theta \neq 0$ and $\tau_n$ never cancels for $n \geq \tilde{N}$.

Similar considerations lead to the same expression for the $\det\Theta$ when $N$ is even. □

***Remark 12*** Notice that

$$\frac{1}{\prod_{i=1}^{M}(2r_i)}\det\Theta''' = \begin{vmatrix} 1\cdots\cdots 1 & 0\cdots\cdots 0 & 1\cdots\cdots 1 \\ 0\cdots\cdots 0 & 1\cdots\cdots 1 & x_1\cdots\cdots x_N \\ r_1^2\cdots\cdots r_M^2 & 0\cdots\cdots 0 & x_1^2\cdots\cdots x_N^2 \\ 0\cdots\cdots 0 & r_1^2\cdots\cdots r_M^2 & x_1^3\cdots\cdots x_N^3 \\ r_1^4\cdots\cdots r_M^4 & 0\cdots\cdots 0 & x_1^4\cdots\cdots x_N^4 \\ \vdots\quad\quad\vdots & \vdots\quad\quad\vdots & \vdots\quad\quad\vdots \\ r_1^{\tilde{N}-1}\cdots\cdots r_M^{\tilde{N}-1} & 0\cdots\cdots 0 & x_1^{\tilde{N}-1}\cdots\cdots x_N^{\tilde{N}-1} \end{vmatrix}$$

$$= \det((X_1(r_1))_{[\tilde{N}]}, \dots, (X_1(r_M))_{[\tilde{N}]}, (X_2(r_1))_{[\tilde{N}]}, \dots, (X_2(r_M))_{[\tilde{N}]}, (X(x_1))_{[\tilde{N}]}, \dots, (X(x_N))_{[\tilde{N}]}).$$

**Lemma 13** *The following relations hold*

$$\begin{aligned}&\left[\Omega_{n+1,n}\cdots\cdots\Omega_{n+\tilde{N},n}\right]\\ &= -\left[v_1^{(n)}(r_1^2)\cdots\cdots v_1^{(n)}(r_M^2)\ v_2^{(n)}(r_1^2)\cdots\cdots v_2^{(n)}(r_M^2)\ v_2^{(n)}(x_1^2)\cdots\cdots v_2^{(n)}(x_N^2)\right]\\ &\times\begin{bmatrix} v_1^{(n+1)}(r_1^2)\cdots\cdots v_1^{(n+1)}(r_M^2) & v_2^{(n+1)}(r_1^2)\cdots\cdots v_2^{(n+1)}(r_M^2) & v_2^{(n+1)}(x_1^2)\cdots\cdots v_2^{(n+1)}(x_N^2) \\ \vdots\quad\quad\vdots & \vdots\quad\quad\vdots & \vdots\quad\quad\vdots \\ v_1^{(\tilde{N}+n)}(r_1^2)\cdots\cdots v_1^{(\tilde{N}+n)}(r_M^2) & v_2^{(\tilde{N}+n-1)}(r_1^2)\cdots\cdots v_2^{(\tilde{N}+n)}(r_M^2) & v_2^{(\tilde{N}+n)}(x_1^2)\cdots\cdots v_2^{(\tilde{N}+n)}(x_N^2) \end{bmatrix}^{-1}.\end{aligned} \quad (54)$$

***Proof*** Direct consequence of Corollary 3. □

In what follows, for the reader convenience, we use the more condensed notation

$$\begin{aligned}&\begin{bmatrix} v_1^{(n)}(r_1^2)\cdots\cdots v_2^{(n)}(x_N^2) \\ \vdots\quad\quad\vdots \\ v_1^{(\tilde{N}+n)}(r_1^2)\cdots\cdots v_2^{(\tilde{N}+n)}(x_N^2) \end{bmatrix}\\ &:= \begin{bmatrix} v_1^{(n)}(r_1^2)\cdots\cdots v_1^{(n)}(r_M^2) & v_2^{(n)}(r_1^2)\cdots\cdots v_2^{(n)}(r_M^2) & v_2^{(n)}(x_1^2)\cdots\cdots v_2^{(n)}(x_N^2) & v_a^{(n)}(x) \\ \vdots\quad\quad\vdots & \vdots\quad\quad\vdots & \vdots\quad\quad\vdots & \vdots \\ v_1^{(\tilde{N}+n)}(r_1^2)\cdots\cdots v_1^{(\tilde{N}+n)}(r_M^2) & v_2^{(\tilde{N}+n)}(r_1^2)\cdots\cdots v_2^{(\tilde{N}+n)}(r_M^2) & v_2^{(\tilde{N}+n)}(x_1^2)\cdots\cdots v_2^{(\tilde{N}+n)}(x_N^2) & v_a^{(\tilde{N}+n)}(x) \end{bmatrix},\end{aligned}$$

**Theorem 11** (Christoffel formulas for Christoffel transformations) *For $n \geq \tilde{N}$ and $a \in \{1, 2\}$, we have the following Christoffel formulas for the perturbed types I and II* *multiple orthogonal polynomials*

$$\hat{A}_a^{(n)}(y) = \frac{1}{\tau_{n+1}\pi(y)p(y)} \begin{vmatrix} v_1^{(n)}(r_1^2) \cdots\cdots v_2^{(n)}(x_N^2) & v_a^{(n)}(y) \\ \vdots \quad\quad \vdots & \vdots \\ v_1^{(\tilde{N}+n)}(r_1^2) \cdots v_2^{(\tilde{N}+n)}(x_N^2) & v_a^{(\tilde{N}+n)}(y) \end{vmatrix},$$

$$\hat{B}^{(n-1)}(x) = \frac{1}{\tau_{n-1}} \begin{vmatrix} k_1(x, r_1^2) \cdots\cdots k_2(x, x_N^2) \\ v_1^{(n)}(r_1^2) \cdots\cdots v_2^{(n)}(x_N^2) \\ \vdots \quad\quad \vdots \\ v_1^{(\tilde{N}+n-2)}(r_1^2) \cdots v_2^{(\tilde{N}+n-2)}(x_N^2) \end{vmatrix}.$$

***Proof*** The type I situation follows from from Lemma 13 and Corollary 3. For the type II polynomials, recalling (51) we can write

$$\begin{aligned} &\left[\hat{B}^{(n-\tilde{N})}(x) \cdots \hat{B}^{(n-1)}(x)\right] \Omega[n] \\ &= -\left[k_1(x, r_1^2) \cdots k_2(x, x_N^2)\right] \begin{bmatrix} v_1^{(n)}(r_1^2) \cdots\cdots v_2^{(n)}(x_N^2) \\ \vdots \quad\quad \vdots \\ v_1^{(\tilde{N}+n-1)}(r_1^2) \cdots v_2^{(\tilde{N}+n-1)}(x_N^2) \end{bmatrix}^{-1}. \end{aligned}$$

In particular, we get

$$\begin{aligned} &\Omega_{n-1+\tilde{N},n-1} \hat{B}^{(n-1)}(x) \\ &= -\left[k_1(x, r_1^2) \cdots k_2(x, x_N^2)\right] \begin{bmatrix} v_1^{(n)}(r_1^2) \cdots\cdots v_2^{(n)}(x_N^2) \\ \vdots \quad\quad \vdots \\ v_1^{(\tilde{N}+n-1)}(r_1^2) \cdots v_2^{(\tilde{N}+n-1)}(x_N^2) \end{bmatrix}^{-1} \begin{bmatrix} 0 \\ \vdots \\ 0 \\ 1 \end{bmatrix} \\ &= -\frac{1}{\tau_n} \begin{vmatrix} v_1^{(n)}(r_1^2) \cdots\cdots v_2^{(n)}(x_N^2) \\ \vdots \quad\quad \vdots \\ v_1^{(\tilde{N}+n-2)}(r_1^2) \cdots v_2^{(\tilde{N}+n-2)}(x_N^2) \\ k_1(x, r_1^2) \cdots\cdots k_2(x, x_N^2) \end{vmatrix}. \end{aligned}$$

Finally, from (54) we obtain

$$\Omega_{n-1+\tilde{N},n-1} = -\left[v^{(n-1)}(x_1^2) \cdots v^{(n-1)}(x_N^2)\right] \begin{bmatrix} v_1^{(n)}(r_1^2) \cdots\cdots v_2^{(n)}(x_N^2) \\ \vdots \quad\quad \vdots \\ v_1^{(\tilde{N}+n-1)}(r_1^2) \cdots v_2^{(N+n-1)}(x_N^2) \end{bmatrix}^{-1} \begin{bmatrix} 0 \\ \vdots \\ 0 \\ 1 \end{bmatrix}$$

$$= (-1)^{\tilde{N}} \frac{\tau_{n-1}}{\tau_n},$$

and we get the announced Christoffel formula. □

## 3.2 General Geronimus transformations

Here we follow the ideas in [3, 10]. In this discussion we set $\mathrm{d}\mu(x) = \mathrm{d}x$. For a polynomial as in (46) we define the Geronimus perturbation $\check{\vec{w}}$ of the vector of weights $\vec{w}$ by the condition

$$p(x)\check{\vec{w}} \begin{bmatrix} P_{\mathrm{e}}(x) & xP_{\mathrm{o}}(x) \\ P_{\mathrm{o}}(x) & P_{\mathrm{e}}(x) \end{bmatrix} = \check{\vec{w}} \begin{bmatrix} \tilde{P}_{\mathrm{e}}(x) & x\tilde{P}_{\mathrm{o}}(x) \\ \tilde{P}_{\mathrm{o}}(x) & \tilde{P}_{\mathrm{e}}(x) \end{bmatrix} = \vec{w}, \tag{55}$$

which determines the following Geronimus perturbation of the vector of weights

$$\check{\vec{w}} = \frac{1}{\pi(x)p(x)}\vec{w} \begin{bmatrix} P_{\mathrm{e}}(x) & -xP_{\mathrm{o}}(x) \\ -P_{\mathrm{o}}(x) & P_{\mathrm{e}}(x) \end{bmatrix} + \vec{c}(x)\Bigg(\sum_{j=1}^{M} \delta(x - r_j^2) + \sum_{j=1}^{N} \delta(x - x_j^2)\Bigg), \tag{56}$$

for some vector function $\vec{c}(x) = \big[c_1(x)\; c_2(x)\big]$. Here $\delta(x-a)$ stands for the Dirac's delta distribution. Then,

**Proposition 16** *The moment matrices satisfy*

$$g_{\check{\vec{w}}} \tilde{P}(\Lambda^\top) = g_{\vec{w}}.$$

***Proof*** It follows from $\Lambda X_1 = xX_2$ and $\Lambda X_2 = X_1$ and the expressions of the moment matrices. □

**Proposition 17** *Let us assume that the moment matrices $g_{\hat{\vec{w}}}$, $g_{\vec{w}}$ have a Gauss–Borel factorization, i.e.*

$$g_{\check{\vec{w}}} = \check{S}^{-1}\check{H}\check{\tilde{S}}^{-\top}, \qquad g_{\vec{w}} = S^{-1}H\tilde{S}^{-\top}.$$

*Then, for the connection matrix*

$$\Omega := \check{S}S^{-1},$$

*one has the alternative expression*

$$\Omega = \check{H}\check{\tilde{S}}^{-\top}\tilde{P}(\Lambda^\top)\tilde{S}^\top H^{-1},$$

*so that is lower unitriangular matrix with only its first $N$ subdiagonals possibly nonzero.*

***Proof*** Direct consequence of the Gauss–Borel factorization. □

**Proposition 18** (Connection formulas) *The following formulas hold*

$$\begin{aligned}\Omega B(x)&=\check{B}(x),\\ \Omega^\top\check{A}_1(x)&=p(x)\big(P_{\mathrm{e}}(x)A_1(x)+xP_{\mathrm{o}}(x)A_2(x)\big),\\ \Omega^\top\check{A}_2(x)&=p(x)\big(P_{\mathrm{o}}(x)A_1(x)+P_{\mathrm{e}}(x)A_2(x)\big).\end{aligned}$$

***Proof*** Use the definitions $B=SX$, $A_1=H^{-1}\tilde{S}X_1$, $A_2=H^{-1}\tilde{S}X_2$ , $\check{B}=\hat{S}X$, $\check{A}_1=\check{H}^{-1}\tilde{\hat{S}}X_1$, $\check{A}_2=\check{H}^{-1}\tilde{\check{S}}X_2$ and the two expressions for $\Omega$. □

**Definition 17** *(Second kind functions)* Let us introduce the Cauchy transforms semi-infinite vectors

$$C_1(z):=\int_\Delta\frac{B(x)}{z-x}w_1(x)\,\mathrm{d}x,\quad C_2(z):=\int_\Delta\frac{B(x)}{z-x}w_2(x)\,\mathrm{d}x,\quad C(z):=\int_\Delta\frac{Q(x)}{z-x}\,\mathrm{d}x.$$

**Proposition 19** (Connection formulas for Cauchy transforms) *The second kind functions are subject to the following relations*

$$\begin{aligned}\Omega C_1(z)&=\tilde{P}_{\mathrm{e}}(z)\check{C}_1(z)+\tilde{P}_{\mathrm{o}}(z)\check{C}_2(z)\\ &\quad-\int_\Delta\check{B}(x)\Big(\check{w}_1(x)\frac{\tilde{P}_{\mathrm{e}}(z)-\tilde{P}_{\mathrm{e}}(x)}{z-x}+\check{w}_2(x)\frac{\tilde{P}_{\mathrm{o}}(z)-\tilde{P}_{\mathrm{o}}(x)}{z-x}\Big)\,\mathrm{d}x,\\ \Omega C_2(z)&=z\tilde{P}_{\mathrm{o}}(z)\check{C}_1(z)+\tilde{P}_{\mathrm{e}}(z)\check{C}_2(z)\\ &\quad-\int_\Delta\check{B}(x)\Big(\check{w}_1(x)\frac{z\tilde{P}_{\mathrm{o}}(z)-x\tilde{P}_{\mathrm{o}}(x)}{z-x}+\check{w}_2(x)\frac{\tilde{P}_{\mathrm{e}}(z)-\tilde{P}_{\mathrm{e}}(x)}{z-x}\Big)\,\mathrm{d}x,\\ C(z)&=\Omega^\top\check{C}(z).\end{aligned}$$

***Proof*** It follows from the definitions of the Cauchy transforms and the connection formulas. The last formula follows from $Q=\Omega^\top\check{Q}$. □

**Corollary 4** *For* $n\geq\max(\deg pP_{\mathrm{e}},\deg xpP_{\mathrm{o}})$ *we have*

$$(\Omega C_1(z))^{(n)}=\tilde{P}_{\mathrm{e}}(z)\check{C}_1^{(n)}(z)+\tilde{P}_{\mathrm{o}}(z)\check{C}_2^{(n)}(z),(\Omega C_2(z))^{(n)}=z\tilde{P}_{\mathrm{o}}(z)\check{C}_1^{(n)}(z)+\tilde{P}_{\mathrm{e}}(z)\check{C}_2^{(n)}(z).$$

***Proof*** It follows from the orthogonality relations satisfied by the Geronimus transformed type II polynomials. □

**Corollary 5** *For* $n\geq\max(\deg\tilde{P}_{\mathrm{e}},\deg x\tilde{P}_{\mathrm{o}})$ *we have*

$$\begin{bmatrix}(\Omega C_1(z))^{(n)} & (\Omega C_2(z))^{(n)}\end{bmatrix}\begin{bmatrix}P_{\mathrm{e}}(z) & -zP_{\mathrm{o}}(z)\\ -P_{\mathrm{o}}(z) & P_{\mathrm{e}}(z)\end{bmatrix}=\pi(z)p(z)\begin{bmatrix}\check{C}_1^{(n)}(z) & \check{C}_2^{(n)}(z)\end{bmatrix}.$$

**Lemma 14** *Given any polynomial* $P\in\mathbb{R}[x]$*, for* $n>\deg P$ *and* $a,b\in\{1,2\}$ *we have that*

$$\mathscr{V}^{P}_{a,b}:=\int_\Delta(\check{B}(t))^\top_{[n]}(\check{A}_a(y))_{[n]}\check{w}_b(t)\frac{P(x)-P(t)}{x-t}\,\mathrm{d}t,$$

*is a bivariate polynomial in* $x,y$*, not depending on* $n$*, of degree, as polynomial in* $x$*, less than* $\deg P$.

***Proof*** It follows from the orthogonality relations fulfilled by $\check{B}$. □

**Definition 18** *(Mixed CD kernels)* For $a,b\in\{1,2\}$, we introduce the mixed CD kernels defined by

$$\begin{aligned}K^{(n)}_{a,b}(x,y)&:=\sum_{l=0}^{n-1}C^{(l)}_b(x)A^{(l)}_a(y)=\int_\Delta\sum_{l=0}^{n-1}B^{(l)}(t)A^{(l)}_a(y)\frac{w_b(t)}{x-t}\,\mathrm{d}t\\&=\int_\Delta\frac{w_b(t)}{x-t}K^{(n)}_a(t,y)\,\mathrm{d}t,\end{aligned}$$

which are Cauchy transforms of the CD kernels $K^{(n)}_1,K^{(n)}_2$.

**Definition 19** Let us also introduce the notation

$$V^{(n)}_1(x):=P_{\mathrm{e}}(x)C^{(n)}_1(x)-P_{\mathrm{o}}(x)C^{(n)}_2(x),\quad V^{(n)}_2(x):=-xP_{\mathrm{o}}(x)C^{(n)}_1(x)+P_{\mathrm{e}}(x)C^{(n)}_2(x).$$

**Theorem 12** (Connection formulas for mixed CD kernels) *For* $n\geq\tilde{N}$*, the following relations connecting the mixed CD kernels are fulfilled*

$$\begin{aligned}&\pi(x)p(x)\begin{bmatrix}\check{K}^{(n)}_{1,1}(x,y)&\check{K}^{(n)}_{1,2}(x,y)\\\check{K}^{(n)}_{2,1}(x,y)&\check{K}^{(n)}_{2,2}(x,y)\end{bmatrix}\\&\quad=\begin{bmatrix}\tilde{P}_{\mathrm{e}}(y)&y\tilde{P}_{\mathrm{o}}(y)\\\tilde{P}_{\mathrm{o}}(y)&\tilde{P}_{\mathrm{e}}(y)\end{bmatrix}\begin{bmatrix}K^{(n)}_{1,1}(x,y)&K^{(n)}_{1,2}(x,y)\\K^{(n)}_{2,1}(x,y)&K^{(n)}_{2,2}(x,y)\end{bmatrix}\begin{bmatrix}P_{\mathrm{e}}(x)&-xP_{\mathrm{o}}(x)\\-P_{\mathrm{o}}(x)&P_{\mathrm{e}}(x)\end{bmatrix}\\&\qquad+\begin{bmatrix}\mathscr{V}^{\tilde{P}_{\mathrm{e}}}_{1,1}(x,y)+\mathscr{V}^{\tilde{P}_{\mathrm{o}}}_{1,2}(x,y)&\mathscr{V}^{x\tilde{P}_{\mathrm{o}}}_{1,1}(x,y)+\mathscr{V}^{\tilde{P}_{\mathrm{e}}}_{1,2}(x,y)\\\mathscr{V}^{\tilde{P}_{\mathrm{e}}}_{2,1}(x,y)+\mathscr{V}^{\tilde{P}_{\mathrm{o}}}_{2,2}(x,y)&\mathscr{V}^{x\tilde{P}_{\mathrm{o}}}_{2,1}(x,y)+\mathscr{V}^{\tilde{P}_{\mathrm{e}}}_{2,2}(x,y)\end{bmatrix}\begin{bmatrix}P_{\mathrm{e}}(x)&-xP_{\mathrm{o}}(x)\\-P_{\mathrm{o}}(x)&P_{\mathrm{e}}(x)\end{bmatrix}\\&\qquad-\begin{bmatrix}\check{A}^{(n)}_1(y)\cdots\cdots\check{A}^{(n+N-1)}_1(y)\\\check{A}^{(n)}_2(y)\cdots\cdots\check{A}^{(n+N-1)}_2(y)\end{bmatrix}\Omega[n]^\top\begin{bmatrix}V^{(n-N)}_1(x)&V^{(n-N)}_2(x)\\\vdots&\vdots\\V^{(n-1)}_1(x)&V^{(n-1)}_2(x)\end{bmatrix}.\end{aligned}$$

***Proof*** Let us consider the kernels

$$\mathscr{K}^{(n)}_{a,b}(x,y):=(C_b(x))^\top_{[n]}\Omega^\top_{[n]}(\check{A}_a(y))_{[n]}.$$

In the one hand, from Proposition 19 we get

$$\begin{aligned}(C_1(x))^\top_{[n]}\Omega^\top_{[n]} &= \tilde P_{\rm e}(x)(\check C_1(x))^\top_{[n]} + \tilde P_{\rm o}(x)(\check C_2(x))^\top_{[n]}\\ &\quad - \int_\Delta (\check B(t))^\top_{[n]}\Big(\check w_1(t)\frac{\tilde P_{\rm e}(x)-\tilde P_{\rm e}(t)}{x-t} + \check w_2(t)\frac{\tilde P_{\rm o}(x)-\tilde P_{\rm o}(t)}{x-t}\Big)\,\mathrm d t,\\ (C_2(x))^\top_{[n]}\Omega^\top_{[n]} &= x\tilde P_{\rm o}(x)(\check C_1(x))^\top_{[n]} + \tilde P_{\rm e}(x)(\check C_2(x))^\top_{[n]}\\ &\quad - \int_\Delta (\check B(t))^\top_{[n]}\Big(\check w_1(t)\frac{x\tilde P_{\rm o}(x)-t\tilde P_{\rm o}(t)}{z-x} + \check w_2(t)\frac{\tilde P_{\rm e}(x)-\tilde P_{\rm e}(t)}{x-t}\Big)\,\mathrm d t.\end{aligned}$$

Hence, according to Lemma 14, for $n > \max(\deg \tilde P_{\rm e}, \deg \tilde P_{\rm o})$ we get

$$\begin{aligned}\mathscr K^{(n)}_{b,1}(x,y) &:= \tilde P_{\rm e}(x)\check K^{(n)}_{b,1}(x,y) + \tilde P_{\rm o}(x)\check K^{(n)}_{b,2}(x,y) - \mathscr V^{\tilde P_{\rm e}}_{b,1}(x,y) - \mathscr V^{\tilde P_{\rm o}}_{b,2}(x,y),\\ \mathscr K^{(n)}_{b,2}(x,y) &:= x\tilde P_{\rm o}(x)\check K^{(n)}_{b,1}(x,y) + \tilde P_{\rm e}(x)\check K^{(n)}_{b,2}(x,y) - \mathscr V^{x\tilde P_{\rm o}}_{b,1}(x,y) - \mathscr V^{\tilde P_{\rm e}}_{b,2}(x,y).\end{aligned}$$

On the other hand, attending to Lemma 12 for $a, b \in \{1, 2\}$ we have

$$\begin{aligned}\mathscr K^{(n)}_{b,a}(x,y) =& (C_a(x))^\top_{[n]}(\Omega^\top(\check A_b(y))_{[n]}\\ &- \begin{bmatrix}\check A^{(n)}_b(y)\cdots\cdots\check A^{(n+N-1)}_b(y)\end{bmatrix}\Omega[n]^\top\begin{bmatrix}C^{(n-N)}_a(x)\\ \vdots\\ C^{(n-1)}_a(x)\end{bmatrix},\end{aligned}$$

and Proposition 18 leads to

$$\begin{aligned}\mathscr K^{(n)}_{1,a}(x,y) &= K^{(n)}_{1,a}(x,y)\tilde P_{\rm e}(y) + K^{(n)}_{2,a}(x,y)y\tilde P_{\rm o}(y)\\ &\quad - \begin{bmatrix}\check A^{(n)}_1(y)\cdots\cdots\check A^{(n+N-1)}_1(y)\end{bmatrix}\Omega[n]^\top\begin{bmatrix}C^{(n-N)}_a(x)\\ \vdots\\ C^{(n-1)}_a(x)\end{bmatrix},\\ \mathscr K^{(n)}_{2,a}(x,y) &= K^{(n)}_{1,a}(x)\tilde P_{\rm o}(y) + K^{(n)}_{2,a}(x,y)\tilde P_{\rm e}(y)\\ &\quad - \begin{bmatrix}\check A^{(n)}_2(y)\cdots\cdots\check A^{(n+N-1)}_2(y)\end{bmatrix}\Omega[n]^\top\begin{bmatrix}C^{(n-N)}_a(x)\\ \vdots\\ C^{(n-1)}_a(x)\end{bmatrix}.\end{aligned}$$

Hence, we get the following equations

$$\begin{aligned}&\tilde P_{\rm e}(x)\check K^{(n)}_{1,1}(x,y) + \tilde P_{\rm o}(x)\check K^{(n)}_{1,2}(x,y) - \mathscr V^{\tilde P_{\rm e}}_{1,1}(x,y) - \mathscr V^{\tilde P_{\rm o}}_{1,2}(x,y)\\ &\quad = K^{(n)}_{1,1}(x,y)\tilde P_{\rm e}(y) + K^{(n)}_{2,1}(x,y)y\tilde P_{\rm o}(y)\end{aligned}$$

$$
-\left[\check{A}_1^{(n)}(y)\cdots\cdots\check{A}_1^{(n+\tilde{N}-1)}(y)\right]\Omega[n]^\top\begin{bmatrix}C_1^{(n-\tilde{N})}(x)\\ \vdots\\ C_1^{(n-1)}(x)\end{bmatrix},
$$

$$
x\tilde{P}_{\mathrm{o}}(x)\check{K}_{1,1}^{(n)}(x,y)+\tilde{P}_{\mathrm{e}}(x)\check{K}_{1,2}^{(n)}(x,y)-\mathscr{V}_{1,1}^{x\tilde{P}_{\mathrm{o}}}(x,y)-\mathscr{V}_{1,2}^{\tilde{P}_{\mathrm{e}}}(x,y)
$$
$$
=K_{1,2}^{(n)}(x,y)\tilde{P}_{\mathrm{e}}(y)+K_{2,2}^{(n)}(x,y)y\tilde{P}_{\mathrm{o}}(y)
$$
$$
-\left[\check{A}_1^{(n)}(y)\cdots\cdots\check{A}_1^{(n+\tilde{N}-1)}(y)\right]\Omega[n]^\top\begin{bmatrix}C_2^{(n-\tilde{N})}(x)\\ \vdots\\ C_2^{(n-1)}(x)\end{bmatrix},
$$

$$
\tilde{P}_{\mathrm{e}}(x)\check{K}_{2,1}^{(n)}(x,y)+\tilde{P}_{\mathrm{o}}(x)\check{K}_{2,2}^{(n)}(x,y)-\mathscr{V}_{2,1}^{\tilde{P}_{\mathrm{e}}}(x,y)-\mathscr{V}_{2,2}^{\tilde{P}_{\mathrm{o}}}(x,y)
$$
$$
=K_{1,1}^{(n)}(x,y)\tilde{P}_{\mathrm{o}}(y)+K_{2,1}^{(n)}(x,y)\tilde{P}_{\mathrm{e}}(y)
$$
$$
-\left[\check{A}_2^{(n)}(y)\cdots\cdots\check{A}_2^{(n+\tilde{N}-1)}(y)\right]\Omega[n]^\top\begin{bmatrix}C_1^{(n-\tilde{N})}(x)\\ \vdots\\ C_1^{(n-1)}(x)\end{bmatrix},
$$

$$
x\tilde{P}_{\mathrm{o}}(x)\check{K}_{2,1}^{(n)}(x,y)+\tilde{P}_{\mathrm{e}}(x)\check{K}_{2,2}^{(n)}(x,y)-\mathscr{V}_{2,1}^{x\tilde{P}_{\mathrm{o}}}(x,y)-\mathscr{V}_{2,2}^{\tilde{P}_{\mathrm{e}}}(x,y)
$$
$$
=K_{1,2}^{(n)}(x,y)\tilde{P}_{\mathrm{o}}(y)+K_{2,2}^{(n)}(x,y)\tilde{P}_{\mathrm{e}}(y)
$$
$$
-\left[\check{A}_2^{(n)}(y)\cdots\cdots\check{A}_2^{(n+\tilde{N}-1)}(y)\right]\Omega[n]^\top\begin{bmatrix}C_2^{(n-\tilde{N})}(x)\\ \vdots\\ C_2^{(n-1)}(x)\end{bmatrix}.
$$

These equations are equivalent to the stated result. □

**Definition 20** For $a\in\{1,2\}$, let us define

$$
W_a(x):=V_a(x)-c_a(x)(\pi p)'(x)B(x),
$$

$$
\begin{bmatrix}\mathscr{P}_{1,1}(x,y) & \mathscr{P}_{1,2}(x,y)\\ \mathscr{P}_{2,1}(x,y) & \mathscr{P}_{2,2}(x,y)\end{bmatrix}:=\begin{bmatrix}\mathscr{V}_{1,1}^{\tilde{P}_{\mathrm{e}}}(x,y)+\mathscr{V}_{1,2}^{\tilde{P}_{\mathrm{o}}}(x,y) & \mathscr{V}_{1,1}^{x\tilde{P}_{\mathrm{o}}}(x,y)+\mathscr{V}_{1,2}^{\tilde{P}_{\mathrm{e}}}(x,y)\\ \mathscr{V}_{2,1}^{\tilde{P}_{\mathrm{e}}}(x,y)+\mathscr{V}_{2,2}^{\tilde{P}_{\mathrm{o}}}(x,y) & \mathscr{V}_{2,1}^{x\tilde{P}_{\mathrm{o}}}(x,y)+\mathscr{V}_{2,2}^{\tilde{P}_{\mathrm{e}}}(x,y)\end{bmatrix}
\begin{bmatrix}P_{\mathrm{e}}(x) & -xP_{\mathrm{o}}(x)\\ -P_{\mathrm{o}}(x) & P_{\mathrm{e}}(x)\end{bmatrix},
$$

$$
\begin{bmatrix}\kappa_{1,1}(x,y) & \kappa_{1,2}(x,y)\\ \kappa_{2,1}(x,y) & \kappa_{2,2}(x,y)\end{bmatrix}:=\begin{bmatrix}\tilde{P}_{\mathrm{e}}(y) & y\tilde{P}_{\mathrm{o}}(y)\\ \tilde{P}_{\mathrm{o}}(y) & \tilde{P}_{\mathrm{e}}(y)\end{bmatrix}\begin{bmatrix}(A_1(y))_{[n]}^\top\\ (A_2(y))_{[n]}^\top\end{bmatrix}\left[(W_1(x))_{[n]}\ (W_2(x))_{[n]}\right]
+\begin{bmatrix}\mathscr{P}_{1,1}(x,y) & \mathscr{P}_{1,2}(x,y)\\ \mathscr{P}_{2,1}(x,y) & \mathscr{P}_{2,2}(x,y)\end{bmatrix}.
$$

***Remark 13*** We have

$$W_1^{(n)}(x) = P_{\rm e}(x)C_1^{(n)}(x) - P_{\rm o}(x)C_2^{(n)}(x) - c_1(x)(\pi p)'(x)B^{(n)}(x),$$
$$W_2^{(n)}(x) = -xP_{\rm o}(x)C_1^{(n)}(x) + P_{\rm e}(x)C_2^{(n)}(x) - c_2(x)(\pi p)'(x)B^{(n)}(x).$$

**Lemma 15** *For $a, b \in \{1, 2\}$ and $\zeta \in Z_{\pi p}$ a root of $\pi(z)p(z)$, see (50), we have:*

*(i) The following limit holds*

$$\lim_{z\to\zeta} \pi(z)p(z)\check{C}_a(z) = c_a(\zeta)(\pi p)'(\zeta)\Omega B(\zeta),$$

*(ii) For $n \geq \max(\deg \tilde{P}_{\rm e}, \deg x\tilde{P}_{\rm o})$ we have*

$$(\Omega W_a(\zeta))^{(n)} = 0.$$

***Proof*** (i) Observe that

$$\begin{aligned}
\check{C}_1(z) &= \int_\Delta \frac{\check{B}(x)}{z-x}\check{w}_1(x)\,\mathrm{d}x = \int_\Delta \frac{\check{B}(x)}{z-x}\frac{\tilde{P}_{\rm e}(x)w_1(x) + \tilde{P}_{\rm o}(x)w_2(x)}{\pi(x)p(x)}\,\mathrm{d}x \\
&\quad + \sum_{i=1}^{M} c_1(r_i^2)\frac{\check{B}(r_j^2)}{z-r_i^2} + \sum_{j=1}^{N} c_1(x_j^2)\frac{\check{B}(x_j^2)}{z-x_j^2}, \\
\check{C}_2(z) &= \int_\Delta \frac{\check{B}(x)}{z-x}\check{w}_2(x)\,\mathrm{d}x \\
&= \int_\Delta \frac{\check{B}(x)}{z-x}\frac{x\tilde{P}_{\rm o}(x)w_1(x) + \tilde{P}_{\rm e}(x)w_2(x)}{\pi(x)p(x)}\,\mathrm{d}x \sum_{i=1}^{M} c_2(r_i^2)\frac{\check{B}(r_i^2)}{z-r_i^2} \\
&\quad + \sum_{j=1}^{N} c_2(x_j^2)\frac{\check{B}(x_j^2)}{z-x_j^2}.
\end{aligned}$$

Consequently,

$$\lim_{z\to r_i^2} \pi(z)p(z)\check{C}_a(z) = c_a(r_i^2)\Omega B(r_i^2)\prod_{k\neq i}(r_i^2 - r_k^2) = c_1(r_i^2)\pi(r_i^2)p'(r_i^2)\Omega B(r_i^2),$$
$$\lim_{z\to x_j^2} \pi(z)p(z)\check{C}_a(z) = c_a(x_j^2)\Omega B(x_j^2)\prod_{k\neq j}(x_j^2 - x_k^2) = c_1(x_j^2)\pi'(x_j^2)p(x_j^2)\Omega B(x_j^2).$$

(ii) Use Corollary 5 and the previous result.

□

**Proposition 20** *It holds that*

$$\kappa_{a,b}^{(n)}(\zeta,y)=\begin{bmatrix}\check{A}_a^{(n)}(y)\cdots\cdots\check{A}_a^{(n+\tilde{N}-1)}(y)\end{bmatrix}\Omega[n]^\top\begin{bmatrix}W_b^{(n-\tilde{N})}(\zeta)\\ \vdots\\ W_b^{(n-1)}(\zeta)\end{bmatrix}. \tag{57}$$

***Proof*** From the previous result we find that

$$\lim_{x\to\zeta}\pi(x)p(x)\check{K}_{a,b}^{(n)}(x,y)=c_b(\zeta)(\pi p)'(\zeta)(\Omega B(\zeta))_{[n]}^\top(\check{A}_a(y))_{[n]},\qquad \zeta\in Z_{\pi p},$$

and using Lemma 12 we obtain

$$\begin{aligned}
&\lim_{x\to\zeta}\pi(x)p(x)\check{K}_{1,a}^{(n)}(x,y)\\
&\quad=c_a(\zeta)(\pi p)'(\zeta)(\tilde{P}_{\mathrm{e}}(y)K_1^{(n)}(\zeta,y)+y\tilde{P}_{\mathrm{o}}(y)K_2^{(n)}(\zeta,y))\\
&\qquad-c_a(\zeta)(\pi p)'(\zeta)\begin{bmatrix}\check{A}_1^{(n)}(y)\cdots\cdots\check{A}_1^{(n+\tilde{N}-1)}(y)\end{bmatrix}\Omega[n]^\top\begin{bmatrix}B^{(n-\tilde{N})}(\zeta)\\ \vdots\\ B^{(n-1)}(\zeta)\end{bmatrix},\\
&\lim_{x\to\zeta}\pi(x)p(x)\check{K}_{2,a}^{(n)}(x,y)\\
&\quad=c_a(\zeta)(\pi p)'(\zeta)\big(\tilde{P}_{\mathrm{o}}(y)K_1^{(n)}(\zeta,y)+\tilde{P}_{\mathrm{e}}(y)K_2^{(n)}(\zeta,y)\big)\\
&\qquad-c_a(\zeta)(\pi p)'(\zeta)\begin{bmatrix}\check{A}_2^{(n)}(y)\cdots\cdots\check{A}_2^{(n+\tilde{N}-1)}(y)\end{bmatrix}\Omega[n]^\top\begin{bmatrix}B^{(n-\tilde{N})}(\zeta)\\ \vdots\\ B^{(n-1)}(\zeta)\end{bmatrix}.
\end{aligned}$$

Then, gathering all this together, we get

$$\begin{aligned}
&\lim_{x\to\zeta}\pi(x)p(x)\begin{bmatrix}\check{K}_{1,1}^{(n)}(x,y) & \check{K}_{1,2}^{(n)}(x,y)\\ \check{K}_{2,1}^{(n)}(x,y) & \check{K}_{2,2}^{(n)}(x,y)\end{bmatrix}\\
&\quad=(\pi p)'(\zeta)\begin{bmatrix}\tilde{P}_{\mathrm{e}}(y) & y\tilde{P}_{\mathrm{o}}(y)\\ \tilde{P}_{\mathrm{o}}(y) & \tilde{P}_{\mathrm{e}}(y)\end{bmatrix}\begin{bmatrix}K_1^{(n)}(\zeta,y)\\ K_2^{(n)}(\zeta,y)\end{bmatrix}\begin{bmatrix}c_1(\zeta) & c_2(\zeta)\end{bmatrix}\\
&\qquad-\begin{bmatrix}\check{A}_1^{(n)}(y)\cdots\cdots\check{A}_1^{(n+\tilde{N}-1)}(y)\\ \check{A}_2^{(n)}(y)\cdots\cdots\check{A}_2^{(n+\tilde{N}-1)}(y)\end{bmatrix}\Omega[n]^\top\\
&\qquad\begin{bmatrix}c_1(\zeta)(\pi p)'(\zeta)B^{(n-\tilde{N})}(\zeta) & c_2(\zeta)(\pi p)'(\zeta)B^{(n-\tilde{N})}(\zeta)\\ \vdots & \vdots\\ c_1(\zeta)(\pi p)'(\zeta)B^{(n-1)}(\zeta) & c_2(\zeta)(\pi p)'(\zeta)B^{(n-1)}(\zeta)\end{bmatrix}.
\end{aligned}$$

Therefore, we arrive to the conclusion that, as $x \to \zeta \in Z_{\pi p}$, we get

$$
\begin{aligned}
&(\pi p)'(\zeta)\begin{bmatrix}\tilde{P}_{\mathrm{e}}(y) & y\tilde{P}_{\mathrm{o}}(y)\\ \tilde{P}_{\mathrm{o}}(y) & \tilde{P}_{\mathrm{e}}(y)\end{bmatrix}\begin{bmatrix}K_1^{(n)}(y,\zeta)\\ K_2^{(n)}(y,\zeta)\end{bmatrix}\begin{bmatrix}c_1(\zeta) & c_2(\zeta)\end{bmatrix}\\
&\quad=\begin{bmatrix}\tilde{P}_{\mathrm{e}}(y) & y\tilde{P}_{\mathrm{o}}(y)\\ \tilde{P}_{\mathrm{o}}(y) & \tilde{P}_{\mathrm{e}}(y)\end{bmatrix}\begin{bmatrix}K_{1,1}^{(n)}(\zeta,y) & K_{1,2}^{(n)}(\zeta,y)\\ K_{2,1}^{(n)}(\zeta,y) & K_{2,2}^{(n)}(\zeta,y)\end{bmatrix}\begin{bmatrix}P_{\mathrm{e}}(\zeta) & -\zeta P_{\mathrm{o}}(\zeta)\\ -P_{\mathrm{o}}(\zeta) & P_{\mathrm{e}}(\zeta)\end{bmatrix}\\
&\qquad+\begin{bmatrix}\mathscr{V}_{1,1}^{\tilde{P}_{\mathrm{e}}}(\zeta,y)+\mathscr{V}_{1,2}^{\tilde{P}_{\mathrm{o}}}(\zeta,y) & \mathscr{V}_{1,1}^{x\tilde{P}_{\mathrm{o}}}(\zeta,y)+\mathscr{V}_{1,2}^{\tilde{P}_{\mathrm{e}}}(\zeta,y)\\ \mathscr{V}_{2,1}^{\tilde{P}_{\mathrm{e}}}(\zeta,y)+\mathscr{V}_{2,2}^{\tilde{P}_{\mathrm{o}}}(\zeta,y) & \mathscr{V}_{2,1}^{x\tilde{P}_{\mathrm{o}}}(\zeta,y)+\mathscr{V}_{2,2}^{\tilde{P}_{\mathrm{e}}}(\zeta,y)\end{bmatrix}\begin{bmatrix}P_{\mathrm{e}}(\zeta) & -\zeta P_{\mathrm{o}}(\zeta)\\ -P_{\mathrm{o}}(\zeta) & P_{\mathrm{e}}(\zeta)\end{bmatrix}\\
&\qquad-\begin{bmatrix}\check{A}_1^{(n)}(y)\cdots\cdots\check{A}_1^{(n+N-1)}(y)\\ \check{A}_2^{(n)}(y)\cdots\cdots\check{A}_2^{(n+N-1)}(y)\end{bmatrix}\Omega[n]^{\top}\begin{bmatrix}W_1^{(n-\tilde{N})}(\zeta) & W_2^{(n-\tilde{N})}(\zeta)\\ \vdots & \vdots\\ W_1^{(n-1)}(\zeta) & W_2^{(n-1)}(\zeta)\end{bmatrix}.
\end{aligned}
$$

Now, recalling the expressions of the CD kernels involved, that is

$$
\begin{aligned}
\begin{bmatrix}K_1^{(n)}(x,y)\\ K_2^{(n)}(x,y)\end{bmatrix}&=\begin{bmatrix}(A_1(y))_{[n]}^{\top}\\ (A_2(y))_{[n]}^{\top}\end{bmatrix}(B(x))_{[n]},\\
\begin{bmatrix}K_{1,1}^{(n)}(x,y) & K_{1,2}^{(n)}(x,y)\\ K_{2,1}^{(n)}(x,y) & K_{2,2}^{(n)}(x,y)\end{bmatrix}&=\begin{bmatrix}(A_1(y))_{[n]}^{\top}\\ (A_2(y))_{[n]}^{\top}\end{bmatrix}\begin{bmatrix}(C_1(x))_{[n]} & (C_2(x))_{[n]}\end{bmatrix},
\end{aligned}
$$

we obtain

$$
\begin{aligned}
&\begin{bmatrix}\tilde{P}_{\mathrm{e}}(y) & y\tilde{P}_{\mathrm{o}}(y)\\ \tilde{P}_{\mathrm{o}}(y) & \tilde{P}_{\mathrm{e}}(y)\end{bmatrix}\begin{bmatrix}(A_1(y))_{[n]}^{\top}\\ (A_2(y))_{[n]}^{\top}\end{bmatrix}\\
&\quad\left((\pi p)'(\zeta)(B(\zeta))_{[n]}\begin{bmatrix}c_1(\zeta) & c_2(\zeta)\end{bmatrix}-\begin{bmatrix}(C_1(\zeta))_{[n]} & (C_2(\zeta))_{[n]}\end{bmatrix}\begin{bmatrix}P_{\mathrm{e}}(\zeta) & -\zeta P_{\mathrm{o}}(\zeta)\\ -P_{\mathrm{o}}(\zeta) & P_{\mathrm{e}}(\zeta)\end{bmatrix}\right)\\
&\quad=\begin{bmatrix}\mathscr{V}_{1,1}^{\tilde{P}_{\mathrm{e}}}(\zeta,y)+\mathscr{V}_{1,2}^{\tilde{P}_{\mathrm{o}}}(\zeta,y) & \mathscr{V}_{1,1}^{x\tilde{P}_{\mathrm{o}}}(\zeta,y)+\mathscr{V}_{1,2}^{\tilde{P}_{\mathrm{e}}}(\zeta,y)\\ \mathscr{V}_{2,1}^{\tilde{P}_{\mathrm{e}}}(\zeta,y)+\mathscr{V}_{2,2}^{\tilde{P}_{\mathrm{o}}}(\zeta,y) & \mathscr{V}_{2,1}^{x\tilde{P}_{\mathrm{o}}}(\zeta,y)+\mathscr{V}_{2,2}^{\tilde{P}_{\mathrm{e}}}(\zeta,y)\end{bmatrix}\begin{bmatrix}P_{\mathrm{e}}(\zeta) & -\zeta P_{\mathrm{o}}(\zeta)\\ -P_{\mathrm{o}}(\zeta) & P_{\mathrm{e}}(\zeta)\end{bmatrix}\\
&\qquad-\begin{bmatrix}\check{A}_1^{(n)}(y)\cdots\cdots\check{A}_1^{(n+\tilde{N}-1)}(y)\\ \check{A}_2^{(n)}(y)\cdots\cdots\check{A}_2^{(n+\tilde{N}-1)}(y)\end{bmatrix}\Omega[n]^{\top}\begin{bmatrix}W_1^{(n-\tilde{N})}(\zeta) & W_2^{(n-\tilde{N})}(\zeta)\\ \vdots & \vdots\\ W_1^{(n-1)}(\zeta) & W_2^{(n-1)}(\zeta)\end{bmatrix}
\end{aligned}
$$

that simplifies to

$$
\begin{bmatrix}\tilde{P}_{\mathrm{e}}(y) & y\tilde{P}_{\mathrm{o}}(y)\\ \tilde{P}_{\mathrm{o}}(y) & \tilde{P}_{\mathrm{e}}(y)\end{bmatrix}\begin{bmatrix}(A_1(y))_{[n]}^{\top}\\ (A_2(y))_{[n]}^{\top}\end{bmatrix}\begin{bmatrix}(W_1(\zeta))_{[n]} & (W_2(\zeta))_{[n]}\end{bmatrix}+\begin{bmatrix}\mathscr{P}_{1,1}(\zeta,y) & \mathscr{P}_{1,2}(\zeta,y)\\ \mathscr{P}_{2,1}(\zeta,y) & \mathscr{P}_{2,2}(\zeta,y)\end{bmatrix}
$$

$$= \begin{bmatrix} \check{A}_1^{(n)}(y) \cdots\cdots \check{A}_1^{(n+\tilde{N}-1)}(y) \\ \check{A}_2^{(n)}(y) \cdots\cdots \check{A}_2^{(n+\tilde{N}-1)}(y) \end{bmatrix} \Omega[n]^\top \begin{bmatrix} W_1^{(n-\tilde{N})}(\zeta) & W_2^{(n-\tilde{N})}(\zeta) \\ \vdots & \vdots \\ W_1^{(n-1)}(\zeta) & W_2^{(n-1)}(\zeta) \end{bmatrix}, \tag{58}$$

and the last relation follows immediately. □

**Definition 21** The analogous tau functions for the Geronimus transformations are defined as

$$\tau_n := \begin{vmatrix} W_1^{(n-\tilde{N})}(r_1^2) \cdots\cdots W_1^{(n-\tilde{N})}(r_M^2) & W_2^{(n-\tilde{N})}(r_1^2) \cdots\cdots W_2^{(n-\tilde{N})}(r_M^2) & W_2^{(n-\tilde{N})}(x_1^2) \cdots\cdots W_2^{(n-\tilde{N})}(x_N^2) \\ \vdots \qquad\qquad \vdots & \vdots \qquad\qquad \vdots & \vdots \qquad\qquad \vdots \\ W_1^{(n-1)}(r_1^2) \cdots\cdots W_1^{(n-1)}(r_M^2) & W_2^{(n-1)}(r_1^2) \cdots\cdots W_2^{(n-1)}(r_M^2) & W_2^{(n-1)}(x_1^2) \cdots\cdots W_2^{(n-1)}(x_N^2) \end{vmatrix}.$$

**Proposition 21** *If $\tau_n = 0$ for some $n \geq \tilde{N}$, then there is a nonzero vector there exists a non-zero vector $\begin{bmatrix} c_1 \\ \vdots \\ c_{\tilde{N}} \end{bmatrix}$ such that, for $n \in \mathbb{N}_0$,*

$$\begin{aligned} &c_1\kappa_{1,1}^{(n)}(r_1^2, y) + \cdots + c_M\kappa_{1,1}^{(n)}(r_M^2, y) + c_{M+1}\kappa_{1,2}^{(n)}(r_1^2, y) + \cdots + c_{2M}\kappa_{1,2}^{(n)}(r_M^2, y) \\ &+c_{2M+1}\kappa_{1,2}^{(n)}(x_1^2, y) + \cdots + c_{\tilde{N}}\kappa_{1,2}^{(n)}(x_N^2, y) = 0. \end{aligned} \tag{59}$$

***Proof*** Let us assume that $\tau_n = 0$. Then, as we have

$$\begin{aligned} &\begin{bmatrix} \check{A}_a^{(n)}(y) \cdots\cdots \check{A}_a^{(n+\tilde{N}-1)}(y) \end{bmatrix} \Omega[n]^\top \\ &\times \begin{bmatrix} W_1^{(n-\tilde{N})}(r_1^2) \cdots\cdots W_1^{(n-\tilde{N})}(r_M^2) & W_2^{(n-\tilde{N})}(r_1^2) \cdots\cdots W_2^{(n-\tilde{N})}(r_M^2) & W_2^{(n-\tilde{N})}(x_1^2) \cdots\cdots W_2^{(n-\tilde{N})}(x_N^2) \\ \vdots \qquad\qquad \vdots & \vdots \qquad\qquad \vdots & \vdots \qquad\qquad \vdots \\ W_1^{(n-1)}(r_1^2) \cdots\cdots W_1^{(n-1)}(r_M^2) & W_2^{(n-1)}(r_1^2) \cdots\cdots W_2^{(n-1)}(r_M^2) & W_2^{(n-1)}(x_1^2) \cdots\cdots W_2^{(n-1)}(x_N^2) \end{bmatrix} \\ &= \begin{bmatrix} \kappa_{1,1}^{(n)}(r_1^2, y) \cdots\cdots \kappa_{1,1}^{(n)}(r_M^2, y) & \kappa_{1,2}^{(n)}(r_1^2, y) \cdots\cdots \kappa_{1,2}^{(n)}(r_M^2, y) & \kappa_{1,2}^{(n)}(x_1^2, y) \cdots\cdots \kappa_{1,2}^{(n)}(x_N^2, y) \end{bmatrix}, \end{aligned} \tag{60}$$

there exists a non-zero vector $\begin{bmatrix} c_1 \\ \vdots \\ c_{\tilde{N}} \end{bmatrix}$ such that (59) is satisfied. □

***Remark 14*** The discussion of when $\tau_n \neq 0$ is still open.

Again, to abbreviate notation we will write

$$\begin{bmatrix} W_1^{(n-\tilde{N})}(r_1^2)\cdots\cdots W_1^{(n-\tilde{N})}(r_M^2) & W_2^{(n-\tilde{N})}(r_1^2)\cdots\cdots W_2^{(n-\tilde{N})}(r_M^2) & W_2^{(n-\tilde{N})}(x_1^2)\cdots\cdots W_2^{(n-\tilde{N})}(x_N^2) \\ \vdots & \vdots & \vdots \\ W_1^{(n-1)}(r_1^2)\cdots\cdots W_1^{(n-1)}(r_M^2) & W_2^{(n-1)}(r_1^2)\cdots\cdots W_2^{(n-1)}(r_M^2) & W_2^{(n-1)}(x_1^2)\cdots\cdots W_2^{(n-1)}(x_N^2) \end{bmatrix}$$
$$= \begin{bmatrix} W_1^{(n-\tilde{N})}(r_1^2)\cdots\cdots W_2^{(n-\tilde{N})}(x_N^2) \\ \vdots \\ W_1^{(n-1)}(r_1^2)\cdots\cdots W_2^{(n-1)}(x_N^2) \end{bmatrix}.$$

**Lemma 16** *Let us assume $\tau_n \neq 0$. Then, the following relations hold*

$$\left[\Omega_{n,n-\tilde{N}}\cdots\cdots\Omega_{n,n-1}\right] = -\left[W_1^{(n)}(r_1^2)\cdots\cdots W_2^{(n)}(x_N^2)\right]\begin{bmatrix} W_1^{(n-\tilde{N})}(r_1^2)\cdots\cdots W_2^{(n-\tilde{N})}(x_N^2) \\ \vdots \qquad\qquad \vdots \\ W_1^{(n-1)}(r_1^2)\cdots\cdots W_2^{(n-1)}(x_N^2) \end{bmatrix}^{-1}.$$

**Theorem 13** (Christoffel formulas for Geronimus perturbations) *Let us assume $\tau_n \neq 0$. For $n \geq \tilde{N}$ and $a \in \{1, 2\}$, we have the following Christoffel formulas for the perturbed multiple orthogonal polynomials*

$$\check{A}_a^{(n)}(y) = \frac{\begin{vmatrix} W_1^{(n-\tilde{N}+1)}(r_1^2)\cdots\cdots W_2^{(n-N+1)}(x_N^2) \\ \vdots \qquad\qquad \vdots \\ W_1^{(n-1}(r_1^2)\cdots\cdots\cdots W_2^{(n-1)}(x_N^2) \\ \kappa_{a,1}(r_1^2, y)\cdots\cdots\cdots\kappa_{a,2}(x_N^2, y) \end{vmatrix}}{\tau_{n+1}},$$

$$\check{B}^{(n)}(x) = \frac{\begin{vmatrix} W_1^{(n-\tilde{N})}(r_1^2)\cdots\cdots W_1^{(n-\tilde{N})}(x_N^2) & B^{(n-\tilde{N})}(x) \\ \vdots \qquad\qquad \vdots & \vdots \\ W_1^{(n-1)}(r_1^2)\cdots\cdots W_2^{(n-1)}(x_N^2) & B^{(n-1)}(x) \\ W_1^{(n)}(r_1^2)\cdots\cdots\cdots W_2^{(n)}(x_N^2) & B^{(n)}(x) \end{vmatrix}}{\tau_n}.$$

***Proof*** The type II formula is proven with the aid of Proposition 18 and Lemma 16. For the type I we first notice that for $a \in \{1, 2\}$ Eq. (58) can be written as follows

$$\left[\kappa_{a,1}(r_1^2, y)\cdots\cdots\kappa_{a,1}(x_N^2, y)\right]\begin{bmatrix} W_1^{(n-\tilde{N})}(r_1^2)\cdots\cdots W_1^{(n-\tilde{N})}(x_N^2) \\ \vdots \qquad\qquad \vdots \\ W_1^{(n-1)}(r_1^2)\cdots\cdots W_1^{(n-1)}(x_N^2) \end{bmatrix}^{-1}$$
$$= -\left[\check{A}_a^{(n)}(y)\cdots\cdots\check{A}_a^{(n+\tilde{N}-1)}(y)\right]\Omega[n]^\top,$$

so that

$$\left[\kappa_{a,1}(r_1^2,y)\cdots\cdots\kappa_{a,2}(x_N^2,y)\right]\begin{bmatrix} W_1^{(n-\tilde{N})}(r_1^2)\cdots\cdots W_2^{(n-\tilde{N})}(x_N^2) \\ \vdots \qquad\qquad \vdots \\ W_1^{(n-1)}(r_1^2)\cdots\cdots W_2^{(n-1)}(x_N^2)\end{bmatrix}^{-1}\begin{bmatrix}1\\0\\\vdots\\0\end{bmatrix} = -\check{A}_a^{(n)}(y)\Omega_{n,n-\tilde{N}}$$

Then, from Lemma 16

$$\Omega_{n,n-\tilde{N}} = -\left[W_1^{(n)}(r_1^2)\cdots\cdots W_2^{(n)}(x_N^2)\right]\begin{bmatrix} W_1^{(n-\tilde{N})}(r_1^2)\cdots\cdots W_2^{(n-\tilde{N})}(x_N^2) \\ \vdots \qquad\qquad \vdots \\ W_1^{(n-1)}(r_1^2)\cdots\cdots W_2^{(n-1)}(x_N^2)\end{bmatrix}^{-1}\begin{bmatrix}1\\0\\\vdots\\0\end{bmatrix}.$$

and we finally get the Christoffel formula for the type I multiple orthogonal polynomials. □

### 3.3 Vectors of Markov–Stieltjes functions

We consider here how the vector of Markov–Siteltjes functions

$$\vec{\mathscr{F}}(z) := \int_\Delta \frac{\vec{w}(x)}{z-x}\,\mathrm{d}\,x$$

behaves under the Christoffel and Geronimus transformations discussed above.

We start with the Christoffel transformation (48). Each of the entries of the vector of Markov–Stieltjes functions transform according to

$$\begin{aligned}\hat{\mathscr{F}}_1(z) &= \int_\Delta \frac{p(x)(P_{\mathrm{e}}(x)w_1(x)+P_{\mathrm{o}}(x)w_2(x))}{z-x}\,\mathrm{d}\,x\\ &= -\int_\Delta \frac{p(z)P_{\mathrm{e}}(z)-p(x)P_{\mathrm{e}}(x)}{z-x}w_1(x)\,\mathrm{d}\,x - \int_\Delta \frac{p(z)P_{\mathrm{o}}(z)-p(x)P_0(x)}{z-x}w_2(x)\,\mathrm{d}\,x\\ &\quad+p(z)(P_{\mathrm{e}}(z)F_1(z)+P_0(z)F_2(z))\end{aligned}$$

and

$$\begin{aligned}\hat{\mathscr{F}}_2(z) &= \int_\Delta \frac{p(x)(xP_{\mathrm{o}}(x)w_1(x)+P_{\mathrm{e}}(x)w_2(x))}{z-x}\,\mathrm{d}\,x\\ &= -\int_\Delta \frac{zp(z)P_{\mathrm{o}}(z)-xp(x)P_{\mathrm{o}}(x)}{z-x}w_1(x)\,\mathrm{d}\,x - \int_\Delta \frac{p(z)P_{\mathrm{e}}(z)-p(x)P_{\mathrm{e}}(x)}{z-x}w_2(x)\,\mathrm{d}\,x\\ &\quad+p(z)(zP_{\mathrm{o}}(z)F_1(z)+P_{\mathrm{e}}(z)F_2(z)).\end{aligned}$$

Hence, Christoffel transformations imply the following affine transformations, with polynomials in $z$ variable as coefficients, for the vector of Markov–Stieltjes functions

$$\hat{\vec{\mathscr{F}}}(z) = p(z)\vec{\mathscr{F}}(z)\begin{bmatrix}P_{\mathrm{e}}(z) & zP_{\mathrm{o}}(z)\\ P_{\mathrm{o}}(z) & P_{\mathrm{e}}(z)\end{bmatrix}$$

$$
-\left(\int_\Delta \frac{p(z)P_{\mathrm{e}}(z)-p(x)P_{\mathrm{e}}(x)}{z-x}w_1(x)\,\mathrm{d}\,x+\int_\Delta \frac{p(z)P_{\mathrm{o}}(z)-p(x)P_{\mathrm{o}}(x)}{z-x}w_2(x)\,\mathrm{d}\,x\right.
$$
$$
\left.\int_\Delta \frac{zp(z)P_{\mathrm{o}}(z)-xp(x)P_{\mathrm{o}}(x)}{z-x}w_1(x)\,\mathrm{d}\,x+\int_\Delta \frac{p(z)P_{\mathrm{e}}(z)-p(x)P_{\mathrm{e}}(x)}{z-x}w_2(x)\,\mathrm{d}\,x\right).
$$

The Geronimus case, see (55) and (56), leads to

$$
\begin{aligned}
p(z)\pi(z)\check{\vec{\mathscr{F}}}(z) &= \begin{bmatrix} p(z)\pi(z)\check{\mathscr{F}}_1(z) & p(z)\pi(z)\check{\mathscr{F}}_2(z)\end{bmatrix}\\
&= \begin{bmatrix} p(z)\pi(z)\int_\Delta \frac{\check{w}_1(x)}{z-x}\,\mathrm{d}\,x & p(z)\pi(z)\int_\Delta \frac{\check{w}_2(x)}{z-x}\,\mathrm{d}\,x\end{bmatrix},
\end{aligned}
$$

with

$$
\begin{aligned}
p(z)\pi(z)\check{\mathscr{F}}_1(z) &= \int_\Delta \frac{p(z)\pi(z)-p(x)\pi(x)}{z-x}\check{w}_1(x)\,\mathrm{d}\,x\\
&\quad+\int_\Delta \frac{P_{\mathrm{e}}(x)w_1(x)-P_{\mathrm{o}}(x)w_2(x)}{z-x}\,\mathrm{d}\,x,\\
p(z)\pi(z)\check{\mathscr{F}}_2(z) &= \int_\Delta \frac{p(z)\pi(z)-p(x)\pi(x)}{z-x}\check{w}_2(x)\,\mathrm{d}\,x\\
&\quad+\int_\Delta \frac{-xP_{\mathrm{o}}(x)w_1(x)+P_{\mathrm{e}}(x)w_2(x)}{z-x}\,\mathrm{d}\,x
\end{aligned}
$$

and, from the above arguments, we finally arrive to the following affine transformations, with coefficients being rational functions in $z$, for the vector of Markov–Stieltjes functions

$$
\begin{aligned}
&p(z)\pi(z)\check{\vec{\mathscr{F}}}(z)-\begin{bmatrix}\int_\Delta \frac{p(z)\pi(z)-p(x)\pi(x)}{z-x}\check{w}_1(x)\,\mathrm{d}\,x & \int_\Delta \frac{p(z)\pi(z)-p(x)\pi(x)}{z-x}\check{w}_2(x)\,\mathrm{d}\,x\end{bmatrix}\\
&=\vec{\mathscr{F}}(z)\begin{bmatrix} P_{\mathrm{e}}(z) & -zP_{\mathrm{o}}(z)\\ -P_{\mathrm{o}}(z) & P_{\mathrm{e}}(z)\end{bmatrix}\\
&\quad+\begin{bmatrix}\int_\Delta \frac{P_{\mathrm{e}}(z)-P_{\mathrm{e}}(x)}{z-x}w_1(x)\,\mathrm{d}\,x-\int_\Delta \frac{P_{\mathrm{o}}(z)-P_{\mathrm{o}}(x)}{z-x}w_2(x)\,\mathrm{d}\,x & \int_\Delta \frac{zP_{\mathrm{o}}(z)-xP_{\mathrm{o}}(x)}{z-x}w_1(x)\,\mathrm{d}\,x\\ & -\int_\Delta \frac{P_{\mathrm{e}}(z)-P_{\mathrm{e}}(x)}{z-x}w_2(x)\,\mathrm{d}\,x\end{bmatrix}.
\end{aligned}
$$

The interested reader can compare these results, technically more difficult, with those for the standard case [45].

## 3.4 The even perturbation

We study the strong simplification that supposes to take $P_{\mathrm{e}}=1$, $P_{\mathrm{o}}=0$, thus $P(x)=p(x^2)$ is an even polynomial. In this case the transformation for the vector of weights is very simple

$$
\hat{\vec{w}}(x)=p(x)\vec{w}(x),\quad p(x)\check{\vec{w}}(x)=\vec{w}(x),\quad \check{\vec{w}}(x)=\frac{1}{p(x)}\vec{w}(x)+\vec{c}(x)\sum_{j=1}^{M}\delta(x-r_j^2).
$$

In this case there is a further simplification coming from the bi-Hankel structure of the moment matrix described in (13). This leads to

$$g_{\hat{\vec{w}}} = g_{\vec{w}} p\big((\Lambda^\top)^2\big) = p\big(\Lambda\big) g_{\vec{w}}, \qquad p\big(\Lambda\big) g_{\check{\vec{w}}} = g_{\check{\vec{w}}} p\big((\Lambda^\top)^2\big) = g_{\vec{w}}.$$

Hence, we can consider a banded upper triangular matrix $\omega$ (with $N$ superdiagonals) for both perturbations. For the Christoffel transformations we have

$$\omega := \hat{S} p(\Lambda) S^{-1} = (H^{-1} \tilde{S} \hat{\tilde{S}}^{-1} \hat{H})^\top,$$

while for the Geronimus perturbation we set

$$\omega := S p(\Lambda) \check{S}^{-1} = (\check{H}^{-1} \check{\tilde{S}} \tilde{S}^{-1} H)^\top.$$

The corresponding connection formulas for the Christoffel transformations are

$$\begin{aligned} p(x)\hat{B}(x) &= \omega B(x), \\ \omega^\top \check{A}_a &= A_a, \end{aligned} \tag{61}$$

while for the Geronimus transformations are

$$\begin{aligned} \omega \check{B}(x) &= p(x) B(x), \\ \check{A}_a(x) &= \omega^\top A_a(x), \end{aligned} \tag{62}$$

so that for the linear forms we find

$$p(x)\check{Q}(x) = \omega^\top Q(x). \tag{63}$$

From our previous experience we see that, in the Christoffel case, Eq. (61) will lead to a expression for the entries of $\omega$ in terms of the the type II polynomials evaluated at the zeros $r_i^2$. Moreover, for the Geronimus case a similar reasoning will give the entries of $\omega$ in terms of the type I linear forms evaluated at zeros $r_i^2$. These allow for Christoffel formulas with no use of kernel polynomials of any type.

**Proposition 22** (Alternative Christoffel formulas for the even perturbation)

*(i) For the even Christoffel perturbation we get the following Christoffel formula*

$$\hat{B}^{(n)}(x) = \frac{1}{p(x)} \frac{\begin{vmatrix} B^{(n)}(x_1^2) \cdots\cdots\cdots B^{(n)}(x_N^2) & B^{(n)}(x) \\ B^{(n+1)}(x_1^2) \cdots\cdots B^{(n+1)}(x_N^2) & B^{(n+1)}(x) \\ \vdots \qquad\qquad \vdots & \vdots \\ B^{(n+N)}(x_1^2) \cdots\cdots B^{(n+N)}(x_N^2) & B^{(n+N)}(x) \end{vmatrix}}{\begin{vmatrix} B^{(n)}(x_1^2) \cdots\cdots\cdots B^{(n)}(x_N^2) \\ \vdots \qquad\qquad \vdots \\ B^{(n+N-1)}(x_1^2) \cdots\cdots B^{(n+N-1)}(x_N^2) \end{vmatrix}}.$$

*(ii) For the even Geronimus transformation we have*

$$\check{A}_a(x) = \frac{\begin{vmatrix} Q^{(n-N)}(x_1) \cdots\cdots Q^{(n-N)}(x_N) & A_a^{(n-N)}(x) \\ \vdots \qquad\qquad \vdots & \vdots \\ Q^{(n)}(x_1) \cdots\cdots Q^{(n)}(x_N) & A_a^{(n)}(x) \end{vmatrix}}{\begin{vmatrix} Q^{(n-N)}(x_1) \cdots\cdots Q^{(n-N)}(x_N) \\ \vdots \qquad\qquad \vdots \\ Q^{(n-1)}x_1) \cdots\cdots Q^{(n-1)}(x_N) \end{vmatrix}}.$$

***Proof*** (i) Equation (61) implies that

$$\begin{aligned} &\left[(\omega)_{n,n+1} \cdots\cdots (\omega)_{n,n+N}\right] \\ &= -\left[B^{(n)}(x_1) \cdots\cdots B^{(n)}(x_N)\right] \begin{bmatrix} B^{(n+1)}(x_1) \cdots\cdots B^{(n+1)}(x_N) \\ \vdots \qquad\qquad \vdots \\ B^{(n+N)}(x_1) \cdots\cdots B^{(n+N)}(x_N) \end{bmatrix}^{-1} \end{aligned}$$

so that, from (61) we get

$$\begin{aligned} \check{B}^{(n)}(x) =& B^{(n)}(x) - \left[B^{(n)}(x_1) \cdots\cdots B^{(n)}(x_N)\right] \\ & \begin{bmatrix} B^{(n+1)}(x_1) \cdots\cdots B^{(n+1)}(x_N) \\ \vdots \qquad\qquad \vdots \\ B^{(n+N)}(x_1) \cdots\cdots B^{(n+N)}(x_N) \end{bmatrix}^{-1} \begin{bmatrix} B^{(n+1)}(x) \\ \vdots \\ B^{(n+N)}(x) \end{bmatrix}, \end{aligned}$$

and we get the stated Christoffel formula follows.

(ii) From (63) we get

$$\begin{aligned} &\left[(\omega^\top)_{n,n-N} \cdots\cdots (\omega^\top)_{n,n-1}\right] \\ &= -\left[Q^{(n)}(x_1) \cdots\cdots Q^{(n)}(x_N)\right] \begin{bmatrix} Q^{(n-N)}(x_1) \cdots\cdots Q^{(n-N)}(x_N) \\ \vdots \qquad\qquad \vdots \\ Q^{(n-1)}(x_1) \cdots\cdots Q^{(n-1)}(x_N) \end{bmatrix}^{-1} \end{aligned}$$

and using (62) we get

$$\begin{aligned} \check{A}_a^{(n)}(x) = & A_a^{(n)}(x) - \left[Q^{(n)}(x_1) \cdots\cdots Q^{(n)}(x_N)\right] \\ & \begin{bmatrix} Q^{(n-N)}(x_1) \cdots\cdots Q^{(n-N)}(x_N) \\ \vdots \qquad\qquad \vdots \\ Q^{(n-1)}(x_1) \cdots\cdots Q^{(n-1)}(x_N) \end{bmatrix}^{-1} \begin{bmatrix} A_a^{(n-N)}(x) \\ \vdots \\ A_a^{(n-1)}(x) \end{bmatrix}, \end{aligned}$$

and the Christoffel formula follows.

□

## Conclusions and outlook

Multiple orthogonal polynomials in general are difficult to tackle given the existence of different weights. In this paper we have considered the simple case with only two weights, and sequences of multiple orthogonal polynomials in the step-line. Among the results in this paper some of the findings are, that the recursion matrix of the type II multiple orthogonal polynomials do not fix the system of weights, there is a hidden symmetry, the permuting Christoffel transformations and corresponding Christoffel formulas and the discussion of Pearson equations and the corresponding differential equations for the orthogonal polynomials and linear forms.

The major contribution of this paper is doubtless the finding of Christoffel formulas for Christoffel and Geronimus perturbations. We use the discussion of the permuting and basic Christoffel transformations to introduce the general theory for the perturbation of vector of weights by polynomial or rational functions. This was a big challenge and we think that here it is presented a pretty general answer to the question on whether there is a general theory that includes Christoffel and Geronimus transformations.

A number of open questions are raised by these findings. In the general Christoffel situation we proved that tau functions were non zero, and we could divide by them. Do we have a similar statement for the tau functions of the Geronimus case? What happens with the mixed Christoffel–Geronimus transformation, also known as Uvarov transformations, for multiple orthogonal polynomials? Second, how to extend the theory to more than two weights, in where a further decomposition of the perturbing polynomial will be needed? Finally, there is a question of generality. Is there a more general polynomial perturbation of the vector of weights that admits Christoffel formulas?

**Acknowledgements** AB acknowledges Centro de Matemática da Universidade de Coimbra (CMUC)—UID/MAT/00324/2020, funded by the Portuguese Government through FCT/MEC and co-funded by the European Regional Development Fund through the Partnership Agreement PT2020.
AFM is partially supported by CIDMA Center for Research and Development in Mathematics and Applications (University of Aveiro) and the Portuguese Foundation for Science and Technology (FCT) within project UID/MAT/04106/2020. MM thanks the financial support from the Spanish "Agencia Estatal de Investigación" research project [PGC2018-096504-B-C33], *Ortogonalidad y Aproximación: Teoría y Aplicaciones en Física Matemática* and research project [PID2021- 122154NB-I00], *Ortogonalidad y aproximación con aplicaciones en machine learning y teoría de la probabilidad*. The authors also acknowledge economical support from ICMAT's Severo Ochoa program mobility B.
The authors are grateful for the excellent job of the referees, whose suggestions and remarks improved the final text.

**Funding** Open Access funding provided thanks to the CRUE-CSIC agreement with Springer Nature.

**Data availability** This paper has no associated data.

## Declarations

**Conflict of interest** The authors declare no conflict of interest.